\documentclass [12pt,leqno]{book}
\usepackage{amsfonts}
\catcode`\@=11
\def\resetthefootnote{\renewcommand{\thefootnote}{\@arabic\c@footnote} }

\def\@principiremex#1{\trivlist
  \item[\hskip \labelsep{\bfseries #1\ \thetheo}]\ignorespaces}
\def\opar@principiremex#1[#2]{\trivlist
  \item[\hskip \labelsep{\bfseries #1\ \thetheo\ (#2)}]\ignorespaces}

\newcommand{\newTHEOremrom}[2]{\newenvironment{#1}{\refstepcounter
  {theo}\@ifnextchar[{\opar@principiremex{#2}}{\@principiremex{#2}}
   }{\qedB\endtrivlist}}
\catcode`\@=12
\DeclareMathSymbol{\square}{\mathord}{AMSa}{"03}
\newcommand{\qed}{\nopagebreak\hspace*{\fill}{\vrule width6pt height6pt depth0pt}\par}
\newcommand{\qedB}{\nopagebreak\hspace*{\fill}$\square$\par}

\newif\ifAddress \Addressfalse

\newtheorem {theo} {Theorem} [section]

\newTHEOremrom{defi}{Definition}
\newTHEOremrom{rem}{Remark}
\newTHEOremrom{conj}{Conjeture}
\newTHEOremrom{ex}{Example}
\newTHEOremrom{nota}{Notation}

\newenvironment{trivlistparm}[2]{\trivlist \item[{#1 #2}]}{\endtrivlist}
\newenvironment{proclama}[1]{\trivlistparm{\bfseries}{#1}\itshape}{\endtrivlistparm}
\newenvironment{prooftext}[1]{\trivlistparm{\itshape}{#1.}}{\qed\endtrivlistparm}
\newenvironment{proof}{\trivlistparm{\itshape}{Proof.}}{\qed\endtrivlistparm}
\newenvironment{case}[2][Case]{
  \trivlistparm{\bfseries}{#1 #2.}}{\endtrivlistparm}

\newcommand{\start}[2]{\begin{#1}\label{#2}}

\newcommand{\theoc}[1]{Theorem~\ref{#1}}
\newcommand{\propc}[1]{Proposition~\ref{#1}}
\newcommand{\coryc}[1]{Corollary~\ref{#1}}
\newcommand{\notc}[1]{Notation~\ref{#1}}
\newcommand{\lemc}[1]{Lemma~\ref{#1}}

\newcommand{\remc}[1]{Remark~\ref{#1}}

\newcommand{\figc}[1]{Figure~\ref{#1}}
\newcommand{\tabc}[1]{Table~\ref{#1}}

\newcommand{\refc}[1]{~\ref{#1}}
\newcommand{\refeq}[1]{\rlap(~\ref{#1})}

\newcommand{\paragrafetiq}[2]{%
  \dimen0=\textwidth \advance\dimen0 by -\parindent
  \setbox0=\hbox{#1} \advance\dimen0 by -\wd0
  \par\par\vspace{\topsep}
  \hbox to \textwidth{\box0\hfill\parbox{\dimen0}{\noindent #2}}
  \par\par\vspace{\topsep}}
\newcounter{llistadepth} \newcounter{llistai}

\newenvironment{manlist}[1]{\addtocounter{llistadepth}{1}
      \edef\llistacontador{llista\romannumeral\the\value{llistadepth}}
      \list{({#1{\llistacontador}})}{\usecounter{\llistacontador}
      \def\makelabel##1{\hss\llap{##1}}
      \itemsep=2pt\parsep=0pt\topsep=3pt plus 1pt minus 1 pt}}{\endlist
      \addtocounter{llistadepth}{-1}}

\newenvironment{romlist}{\begin{manlist}{\roman}}{\end{manlist}}
\newenvironment{numlist}{\begin{manlist}{\arabic}}{\end{manlist}}

\def\Int{\mathop\mathrm{Int}}

\def\Sing{\mathop\mathrm{Sing}}
\def\Cl{\mathop\mathrm{Cl}}
\def\id{\mathop\mathrm{Id}}
\def\Card{\mathop\mathrm{Card}}
\def\Per{\mathop\mathrm{Per}}

\def\Ker{\mathrm{Ker}}
\def\ima{\mathrm{Im}(\phi)}
\def\lcm{\mathrm{l.c.m}}
\def\gcd{\mathrm{g.c.d}}
\def\genus{\mathrm{genus}}
\def\bc{\mathrm{bc}}
\def\ggg{\mathrm{g.c.d}}
\def\suma{\mathrm{p}}

\def\Fix{\mathop\mathrm{Fix}}

\def\tr{\mathop\mathrm{trace \/}}
\def\d{\mathop\mathrm{det \/}}

\def\CC{\mathop\mathcal{C}}

\def\tf{\widetilde f}
\def\fs{f_{*1}}
\def\ffs{\widetilde{f}_{*1}}
\def\v{\varphi}
\def\vv{\vartheta}

\def\campo{v}

\def\sp{s_{_P}}

\def\sf{{\sigma_f}}
\def\si{{\sigma_i}}
\def\sfi{\sigma_{f^{n_i}|_{C_i}}}

\def\s{\sigma}
\def\ve{\varepsilon}
\def\uuu{\vartheta}
\def\gg{\gamma}

\def\inn{\mathrm{I}}
\def\deg{\mathrm{deg}}
\def\mmm{\mathrm{m}}
\def\Suu{\mbox{$\Sigma$}}
\def\Su#1{\mbox{$\Sigma_{#1}$}}
\def\Suf#1#2{\mbox{$\Sigma_{#1,#2}$}}
\def\typeop#1#2{\mbox{$[#1;0;\{#2\}]$}}
\def\typeor#1#2#3{\mbox{$[#1;#2;\{#3\}]$}}

\def\hg{\mathop\mathcal{H}\nolimits _{g}}
\def\hgp{\mathop\mathcal{H}\nolimits _{g}^+}
\def\hgr{\mathop\mathcal{H}\nolimits _{g}^-}
\def\hgpx#1{\mathop\mathcal{H}\nolimits _{#1}^+}
\def\hgrx#1{\mathop\mathcal{H}\nolimits _{#1}^-}

\def\hgb{\mathop\mathcal{H}\nolimits  _{g,b}}
\def\hgbp{\mathop\mathcal{H}\nolimits _{g,b}^+}
\def\hgbr{\mathop\mathcal{H}\nolimits _{g,b}^-}
\def\hgbx#1#2{\mathop\mathcal{H}\nolimits _{#1,#2}}
\def\hgbpx#1#2{\mathop\mathcal{H}\nolimits _{#1,#2}^+}
\def\hgbrx#1#2{\mathop\mathcal{H}\nolimits _{#1,#2}^-}

\def\fg{\mathop\mathcal{F}\nolimits _{g}}
\def\fgp{\mathop\mathcal{F}\nolimits _{g}^+}
\def\fgr{\mathop\mathcal{F}\nolimits _{g}^-}
\def\fgpx#1{\mathop\mathcal{F}\nolimits _{#1}^+}
\def\fgrx#1{\mathop\mathcal{F}\nolimits _{#1}^-}

\def\fgb{\mathop\mathcal{F}\nolimits _{g,b}}
\def\fgbp{\mathop\mathcal{F}\nolimits _{g,b}^+}
\def\fgbr{\mathop\mathcal{F}\nolimits _{g,b}^-}
\def\fgbpx#1#2{\mathop\mathcal{F}\nolimits _{#1,#2}^+}
\def\fgbrx#1#2{\mathop\mathcal{F}\nolimits _{#1,#2}^-}

\def\CC{\mathop\mathcal{C}}
\def\ind#1#2{{\rm Ind}_{#1}({#2})}
\def\ii#1#2{\mbox{${#1}~\in~\{1,2,\dots,{#2}\}$}}

\def\map#1#2#3{\mbox{${#1}\colon {#2} \longrightarrow {#3}$}}
\def\Smap#1#2{\mbox{${#1}\colon{#2} \longrightarrow {#2}$}}
\def\union#1#2{\mbox{$\bigcup\limits_{#1}^{#2}$}}

\def\clase#1{\mbox{$[{#1}]_{n}$}}
\def\minrep#1{\mbox{$\widetilde{{#1}}$}}

\newcommand{\N}{\ensuremath{\mathbb{N}}}
\newcommand{\Z}{\ensuremath{\mathbb{Z}}}
\newcommand{\Q}{\ensuremath{\mathbb{Q}}}
\newcommand{\R}{\ensuremath{\mathbb{R}}}
\newcommand{\C}{\ensuremath{\mathbb{C}}}

\newcommand{\SI}{\ensuremath{\mathbb{S}^1}}

\newcommand{\Sn}{\ensuremath{\mathbb{S}^{n}}}
\newcommand{\D}{\ensuremath{\mathbb{D}^2}}

\def\ind#1#2{{\rm Ind}_{#1}({#2})}

\usepackage{amscd}
\usepackage{latexcad}
\usepackage{makeidx}
\usepackage{pstricks}
\makeindex
\oddsidemargin
\evensidemargin
\begin{document}

\pagestyle{plain}
\pagenumbering{roman}

\pagestyle{empty}

\begin{center}\begin{Large}
\begin{bf}
%
%
Minimum periods of  homeomorphisms of orientable surfaces
\\[6cm]
\end{bf}
\end{Large}
Moira Chas \\ 1998

\end{center}

\vspace*{7.5cm}
\begin{flushright}\parbox{7cm}{

Mem\`{o}ria presentada per a aspirar al grau de doctor en Ci\`{e}ncies Matem\`{a}tiques.

Departament de Matematiques \newline Universitat Autonoma de Barcelona \newline
Bellaterra, abril de 1998.

}

\end{flushright}

\newpage
\vspace*{13.5cm}
\newpage
\vspace*{13.5cm}\parbox{10cm}{

Certifiquem que la present mem\`{o}ria ha estat realitzada per la Moira Chas, i dirigida per
nosaltres al Departament de Matem\`{a}tiques de la Universitat Aut\`{o}noma de Barcelona.

\vspace{1cm}

Bellaterra, abril de 1998.

\vspace{3cm}

\begin{flushright}

Dr. Llu\'{i}s Alsed\`{a} i Soler \hfill  Dr. Warren Dicks  i McLay

\end{flushright}}
\newpage
\mbox{}
\newpage\newpage
\vspace*{13.5cm}\parbox{8cm}{

To the beloved memory of Dante Leal, \newline who laughed, struggled and loved so much
\newline in such a short time. }


\chapter*{Acknowledgements}
\pagestyle{empty}

I would like to thank  Prof. Warren Dicks for being my advisor during the final years of
the preparation of this thesis. My knowledge of English (and probably of Spanish!) is not
enough to express how much I owe him. He always had time and energy to talk to me, to
search for answers for my numerous (and often silly) questions, and to encourage me in
hard times. I learned a great deal from him through long hours of exciting mathematical
conversations. I have no doubt that without his advice, patience and generous help this
thesis would not have been possible.

I owe special thanks to Prof. Llu\'{\i}s Alsed\`{a}, who was my advisor during the first
years, for helping me in uncountably many situations and for standing by me. He was the
first (almost heroic) reader of the chaotic early version of this work and his
suggestions contributed a great deal to making it less incomprehensible.

Thanks also are due to Prof. Jaume Llibre who proposed the problem which originated this
thesis as  a course project.

This project was begun with Victor Manosa, and I would like to thank him because his
influence survives in Chapter\ref{S10}, in some results we obtained together.

I discussed the ideas which appear in this thesis with many people and these discussions
opened doors and windows of my mind to new mathematical light. I am particularly grateful
to Profs. Jane Gilman, John Guaschi, Boju Jiang and Dennis Sullivan.

During most of this period I was supported by  the DGICYT through grant number PB93-0860,
thanks to the Dynamical Systems Group of  the Universitat Aut\`{o}noma de Barcelona. I
wish to express my gratitude to them, and also to the Departament de Matem\`{a}tiques of
the Universitat Aut\`{o}noma de Barcelona for treating me so kindly during these years.

I wish to thank everyone at the State University of New York at Stony Brook and the City
University of New York  who made my stay in New York so pleasant and stimulating.

The more elegant drawings were made by Mauricio; I am very grateful to him for the long
time he spent preparing them without losing his smile.

I owe an enormous debt of gratitude to my friends. When this work lead me to mathematical
heavens or hells, or  to labyrinths without apparent exit, there was always one of them
at my side listening to commentaries, complaints and sighs, supporting, encouraging and a
large etcetera. I particularly thank  Marcelo, Fabiana, Carola, Mariana, Olivier, Montse,
Antonia, Blanca, Maria Jose, Marcia, Ana, Anna, Dora, Pilar, Carlos,  Antonio and Marina.
Warm thanks are due to my mother, for reasons which are too profound to be explained
here.

It is not easy to teach mathematics, and it is even more difficult to transmit the love
of it, which includes the absolute concentration, the emotions created by a beautiful
idea, the often obsessive search for an answer, the happiness and the relief at finding
an answer, and the necessary strength and will to continue when there is neither
understanding  nor answers. I learnt all of this from the people of the Departamento de
Matem\'{a}tica of the Universidad de Buenos Aires and from Prof. Jean-Marc Gambaudo and I
wish to thank them for making my first years of mathematics so passionate and
stimulating.

\tableofcontents

\pagenumbering{arabic}
\renewcommand{\thetheo}{\thechapter.\arabic{theo}}

\chapter{Introduction and main results}\label{S1}

Given a set $M$ and a map $f$ from $M$ to itself, recall  the \emph{orbit} \index{orbit} of a point $x \in M$ is the set $\{x, f(x), f^2(x),...\}$. The data $(M,f)$ is  a \emph{discrete dynamical system}, and the aim of the theory of dynamical systems is to understand the structure of the set of all orbits. The simplest kind of orbit  is a \emph{periodic orbit} \index{periodic orbit}, that is,  the orbit of a point $x$ for which there exists a positive integer $n$ such that $f^n(x)=x$. If $x$ is a point contained in a periodic orbit then it is called a \emph{periodic point} \index{periodic point}, and the least $n$ such that $f^n(x)=x$ is defined to be its \emph{period}. 

To obtain information about the structure of the set of orbits, we need to specify the data$(M,f)$. In general,  $M$ has some  structure and one considers maps preserving this structure. This specification can be done in several ways. For instance, $M$ can be a topological space and $f$ a continuous map, or $M$ a differentiable manifold and $f$ a differentiable map, or $M$ equipped with a $\sigma$-algebra structure and a measure and $f$ preserves this measure or at least sets of measure zero.

Periodic orbits have always been an object of special interest in dynamical systems. Not only because of their simplicity but because their existence often has strong consequences for the dynamics of the map.  In this sense, it is sometimes said that the set of orbits is the skeleton of the set of all orbits. Thus, another interesting object is the set of periods of $f$, which is the subset of $\N$ consisting of the periods of all periodic points in $M$. This set is denoted by $\mathrm{Per}(f)$. \index{$\mathrm{Per}(f)$}

It is interesting  to deal with all the questions related to the periodic orbits in the field of \emph{topological dynamics}, i.e., where $M$ is a topological space and $f$ a continuous map. Here, it is often possible to obtain useful information about the structure of the set of orbits. In particular, much is known about the set $\mathrm{Per}(f)$ when $M$ has dimension one; see, for example, \cite{ALM} and references therein. There are also many results for other spaces and classes of maps.

One of the main problems of the theory of dynamical systems is the determination of the existence of periodic orbits and, more generally the structure of $\mathrm{Per}(f)$ We define the \emph{minimum period of $f$} \index{minimum period of $f$} to be the maximum $m$, positive or infinite, such that the iterates $f, f^2, \dots, f^{m-1}$ are fixed point free. We denote this number by $\mathrm{m}(f)$. Observe that the minimum period of $f$ is the greatest lower bound of $\mathrm{Per}(f)$. If $\mathcal{C}$ is a class of maps, the \emph{minimum period of $\mathcal{C}$} \index{minimum period of $\mathcal{C}$} is defined to be the maximum of the minimum periods of the mappings in $\mathcal{C}$, and is denoted by $\mathrm{m}(\mathcal{C})$.

Throughout this thesis we are going to deal with some aspects of the above problem in the case of two-dimensional surfaces. More precisely, we will study the case where $M=\Sigma$, a connected orientable compact surface (possibly with boundary) and the maps of $\Sigma$ we will consider are the orientation-preserving and orientation-reversing homeomorphisms. The reason for studying homeomorphisms of $\Sigma$ and not just continuous maps is that for each surface $\Sigma$ of genus $g$ at least one there exists a map $\Smap{f}{\Sigma}$ with no periodic points, i.e., such that $\mathrm{m}(f)=\infty$. Indeed, there exist a simple closed curve $\gamma \subset \Sigma$ and a map $\map{g}{\Sigma}{\gamma}$ which is the identity on $\gamma$.There exists a map $\Smap{h}{\gamma}$ conjugate to an irrational rotation $\Smap{r}{\SI}$, so $h$ has no periodic points. Clearly, we can view $h \circ g$ as a self-map of $\Sigma$ and since $\mathrm{Per}(h \circ g)=\emptyset$, $\mathrm{m}(h \circ g)=\infty$.

Before going on, let us introduce some notation: The class of all (resp. all orientation-preserving, resp. all orientation-reversing) homeomorphisms of $\Sigma_{g,b}$ will be denoted by $\hgb$ (resp. $\hgbp$, $\hgbr$). Analogously, the class of all  (resp. all orientation-preserving, resp. all orientation-reversing) homeomorphisms of $\Sigma_{g}$ will be denoted by $\hg$ (resp. $\hgp$, 
$\hgr$).

In the case of closed surfaces, the problem of determining the minimum period of the classes of orientation-preserving and orientation-reversing homeomorphisms is completely solved. The aim of this thesis is to study the minimum periods of homeomorphisms of surfaces with non-empty boundary, i.e., $\mmm(\hgbp)$ and $\mmm(\hgbr)$. 

Let us return to the case of closed surfaces. Here, both bounds, $\mmm(\hgp)$ and $\mmm(\hgr)$ can be explicitly expressed as a function of the genus $g$. Their values are summarized in the following formulas.

\begin{equation}
\mmm(\hgp) = \left\{\begin{array}{ll}
        1& \mbox{if $g=0$,}\\
         \infty & \mbox{if $g=1$}, \\
        2g-2 & \mbox{if $g \ge 2$.}\\
        \end{array}
        \right.
        \label{?}
\end{equation}

\begin{equation}
\mmm(\hgr) = \left\{\begin{array}{ll}
        2 & \mbox{if $g=0$,}\\
         \infty & \mbox{if $g=1$}, \\
          4 & \mbox{if $g=2$}, \\
        2g-2 & \mbox{if $g \ge 2$.}\\
        \end{array}
        \right.
        \label{??}
\end{equation}

The first and well-known result in this field appeared in 1910. It is Brouwer's theorem \cite{B} that an orientation-preserving homeomorphism of the sphere always has a fixed point. With the notation we have introduced, this result can be expressed as $\mmm(\hgpx{0})=1$.

It is a simple matter to check that $\mmm(\hgpx{1})$ and $\mmm(\hgrx{1})$ are both infinite by exhibiting examples of orientation-preserving and orientation-reversing homeomorphisms of the torus $\Sigma_1$ with no periodic points. Indeed, view the torus $\Sigma_1$ as $\SI \times \SI$, and consider the homeomorphism $(z,w) \rightarrow (z e^{\alpha i}, w)$, where $\alpha$ is an irrational real number. It is clear that this homeomorphism preserves and does not have fixed points. Similarly, the orientation reversing homeomorphism $(z,w) \rightarrow (z e^{\alpha i}, \overline{w})$ does not have periodic points.

In \cite{N} Nielsen showed that  $\mmm(\hgp)=2g-2$ if $g\ge 3$. He also showed that $\mmm(\hgpx{2}) \in \{2,3\}$. The proof of these results uses the fixed-point theory due to Alexander  \cite{A} and Lefschetz \cite{L1} and \cite{L2} and some elementary algebra.

Later, Wang \cite{W2} showed that $\mmm(\hgr)=2g-2$ if $g \ge 3$ and $\mmm(\hgrx{2})=4$ by using methods analogous to these of Nielsen. 

The problem about the determination of $\mmm(\hgpx{2})$, raised by Nielsen in 1942, remained open until 1996 when Dicks and Llibre \cite{DL} gave an algebraic proof that $\mmm(\hgpx{2})=2$, which completes \refeq{?}.

The only remaining case in \refeq{??} is $\mmm(\hgrx{0})$. Since the antipodal map is an orientation-reversing homeomorphism of $\Sigma_0$ which is fixed point free, $\mmm(\hgrx{0}) \ge 2$. The equality $\mmm(\hgrx{0})=2$ can be deduced from a theorem of Fuller \cite{F}. This theorem states the existence of, and gives a bound for, the minimum periods of classes of homeomorphisms of compact ANRs. (See below for a definition of compact ANRs). In particular, it gives a general bound (and so the finiteness) for $\mmm(\hgb)$ except for two particular cases for which, as we shall see, the minimum period is $\infty$. Before stating this result, we require some terminology.

A subset $A$ of a topological metric space $X$ is called an \emph{compact absolute neighborhood retract} \index{absolute neighborhood retract}\index{ANR} (or, briefly, \emph{compact ANR})  if it has the following property: If $A$ is a subspace of a separable metric space $Y$ and $A$ is homemorphic to $X$, then $A$ is a neighborhood retract of $Y$.

If $K$ is a compact ANR  we denote by $H_k(K;\Q)$ \index{$H_k(K;\Q)$} the $k$-th rational homology group of $K$. For each $k$, the dimension of $H_k(K,\Q)$  is called the \emph{$k$-th Betti number of $K$} \index{$k$-th Betti number of $K$} and denoted by $b_k(K)$.\index{$b_k(K)$} The \emph{Euler characteristic of $K$} \index{Euler characteristic of $K$}  is denoted by $\chi(K)$ \index{$\chi(K)$}  and is defined to be $\sum(-1)^jb_k(K)$, a finite sum. Now we are ready to state Fuller's theorem. We will do it in the slightly more general version given in \cite[Theorem III.E.2]{Bn}.

\begin{theo} Let $K$ be a compact ANR. If $\chi(K) \ne 0$ and $\Smap{T}{K}$ is a homeomorphism then
$$
\mmm(T) \le \max \left\{ \sum_{j\mbox{ odd}}(-1)^jb_k(K) , \sum_{j \mbox{ even} }(-1)^jb_k(K) \right\}
$$
\end{theo}

In particular, surfaces are compact ANR's and their Betti numbers are
\begin{eqnarray*}
b_k(\Sigma_g) = \left\{\begin{array}{ll}
        1 & \mbox{if $k \in \{0,2\}$,}\\
         2g & \mbox{if $k=1$}, \\
        0 & \mbox{if $k \ge 3$,}\\
        \end{array}
        \right.
        \nonumber\end{eqnarray*}
and
\begin{eqnarray*}
b_k(\Sigma_{g,b}) = \left\{\begin{array}{ll}
        1 & \mbox{if $k =0$,}\\
         2g+b-1 & \mbox{if $k=1$}, \\
        0 & \mbox{if $k \ge 2$.}\\
        \end{array}
        \right.
       \nonumber
\end{eqnarray*}
Therefore, $\chi(\Sigma_{g})=2-2g$ and $\chi(\Sigma_{g,b})=2-2g-b$, so Fuller's theorem restricted to surfaces can be written in the following way:

\begin{theo}\label{fuller} (Fuller's Theorem for surfaces) \begin{numlist}
\item If $g \ne1$ then $\mmm(\hg) \le \max\{2,2g\}.$
\item  If $(g,b) \ne (0,2)$ then $\mmm(\hgb) \le \max\{1,2g+b-1\}.$
\end{numlist}
\end{theo}

In particular, this result implies that $\mmm(\hgbp)$ and $\mmm(\hgbr)$ are finite whenever $(g,b) \ne (0,2)$. Also, $\mmm(\hgrx{0})\le 2$, which completes \refeq{??}.

A subclass of homeomorphisms whose minimum period is also known is the class of finite-order maps. A homeomorphism $\Smap{f}{\Sigma}$ is said to be \emph{finite-order} \index{finite-order} if there exists some positive integer $n$ such that $f^n=\id_\Sigma$. 
The class of all (resp. all orientation-preserving, resp. all orientation-reversing) finite-order maps of $\Sigma_g$ is denoted by $\fg$ (resp. $\fgp$, $\fgr$). Nielsen \cite{N} and Wang \cite{W2} determined the minimum period for $\fgp$ and $\fgr$, respectively. These results, together with the simple cases where $g \in \{0,1\}$ are summarized in the following formulas.

\begin{eqnarray*}
\mmm(\fgp) = \left\{\begin{array}{ll}
        1 & \mbox{if $g =0$,}\\
         \infty & \mbox{if $g=1$}, \\
         2 & \mbox{if $g=2$}, \\
        g-1 & \mbox{if $g \ge 3$,}\\
        \end{array}
        \right.
       \nonumber
\end{eqnarray*}

\begin{eqnarray*}
\mmm(\fgr) = \left\{\begin{array}{ll}
        2 & \mbox{if $g =0$,}\\
         \infty & \mbox{if $g=1$}, \\
         4 & \mbox{if $g=2$}, \\
        2g-2 & \mbox{if $g \ge 3$.}\\
        \end{array}
        \right.
       \nonumber
\end{eqnarray*}

\section{Statement of the main results}

Now we discuss the object of our study, surfaces with non-empty boundary. 

It follows from Fuller's theorem \ref{fuller}(2) that $\mmm(\hgbx{0}{1} )=1 $ and 
$\mmm(\hgb) \le 2g+b-1$ when $(g,b) \notin \{(0,1), (0,2)\}$. However, as we will see, except for $(g,b)\in \{(0,1),(1,1)\} $, the bounds given by that theorem are not the best possible. For example, the following gives a strictly smaller bound for $\mmm(\hgb)$ in most cases.

\begin{proclama}{Proposition A} If $2g+b\ge 4$ then $\mmm(\hgb) \le 2g+b-2$.
\end{proclama}

We denote the class of all (resp. all orientation preserving, resp. all orientation reversing) finite-order maps of $\Sigma_{g,b}$ by $\fgb$ (resp. $\fgbp$, $\fgbr$). By exhibiting specific maps we will show that the bound given by Proposition A can be achieved if the pair $(g,b)$ satisfies certain numerical conditions. Moreover, these maps are finite-order, so we have the following two theorems. 

\begin{proclama}{Theorem B}  Let $g \ge 2$. Then $\mmm(\fgbp) = 2g+b-2$ if and only if $b \in\{2,3,4\}$ or there exist positive integers $p_1, p_2, p_3$ such that they are pairwise coprime, each of them divides $2g+b-2$ and $p_1+p_2+p_3=b$. 
\end{proclama}

\begin{proclama}{Theorem C} Let $g \ge 2$. Then $\mmm(\fgbp) = 2g+b-2$ if and only if $b \in\{2,4\}$ or one of the following conditions holds:
\begin{numlist}
\item $g$ is even and there exist positive integers $p_1, p_2$ such that  each of them divides $2g+b-2$, $\mathrm{g.c.d}(p_1,p_2)=2$  and $p_1+p_2=b$. 
\item $g$ is odd, $b$ is even, and $b$ divides $2g-2$.
\end{numlist}
\end{proclama}

By Proposition A, and Theorem B, $\mmm(\hgbp)=2g+b-2$ for certain pairs $(g,b)$. Also, if $g\ge 2$ it can be proved that there exists a homeomorphism $f \in \hgpx{g}{1}$ such that $\mmm(f)=2g-1$ (and, clearly, $f \notin \fgpx{g}{1}$). Hence, we have the following.

\begin{proclama}{Theorem D} Let $g \ge 2$. Then  $\mmm(\hgbp)=2g+b-2$ if one of the conditions holds.
\begin{numlist}
\item There exist positive integers $p_1, p_2, p_3$ such that they are pairwise coprime, each of them divides $2g+b-2$ and $p_1+p_2+p_3=b$.
\item $b-2$ divides $2g$.
\item $b-3$ divides $2g+1$.
\item $b \in \{1,2,3,4,g+2,2g+2, 2g+4\}.$
\end{numlist}
\end{proclama}

Also, using Theorem C, we prove the following.

\begin{proclama}{Theorem E} Let $g \ge 2$. 
\begin{numlist}
\item If $b$ is odd then $\mmm(\hgbr)\le b$ and equality holds if $b \le 2g-2$. 
\item $\mmm(\hgbr)=2g+b-2$ if one of the following conditions holds.
\begin{romlist}
\item $b \in \{2,4\}$.
\item $g$ is odd, $b$ is even and $b$ divides $2g-2$.
\item $g$ is even, and  there exist positive integers $p_1, p_2$ such that  each of them divides $2g+b-2$, $\mathrm{g.c.d}(p_1,p_2)=2$  and $p_1+p_2=b$. 
\item $g$ is even and $b-2 $ divides $2g$
\item $g$ is even and $b-4 $ divides $2g+2$
\item $g$ is even and $b \in \{g+2, 2g+2, 2g+6\}$.
\end{romlist}
\end{numlist}
\end{proclama}

We give values of the minimum periods for orientation-preserving (resp. orientation-reversing) homemorphisms of $\Sigma_{0,b}$, $\Sigma_{1,b}$, and $\Sigma_{2,b}$ in Theorem F (resp. Theorem G). Notice that in these cases, the bounds of Proposition A are achieved if an only if $g=0$ and $b \ge 3$ (resp. $g \ge 2$), $g=1$ and $b \ge 2$ (resp. $b \ge 1$), or $g \ge 2$ and $b$ satisfies one of the conditions listed in Theorem D (resp. Theorem E). 

\begin{proclama}{Theorem F}  \begin{numlist}  
\item \begin{eqnarray*}  
\mmm(\hgbpx{0}{b}) = \left\{\begin{array}{ll}
        1 & \mbox{if $b =1$,}\\
         \infty & \mbox{if $b=2$}, \\
        b-2 & \mbox{if $b \ge 3$.}\\
        \end{array}
        \right.
\end{eqnarray*}
\item \begin{eqnarray*}  
\mmm(\hgbpx{1}{b}) = \left\{\begin{array}{ll}
        2 & \mbox{if $b =1$,}\\
        b & \mbox{if $b \ge 2$.}\\
        \end{array}
        \right.
\end{eqnarray*}
\item Table \ref{ttaabb} shows the values of $\mmm(\hgbpx{2}{b})$.

\begin{table}
\begin{center}
\begin{tabular}{|l|llllllllllllll|}
\hline b &1&2&3&4&5&6&7&8&9&10&11&12&13&14
\\ \hline $\mmm(\hgbpx{2}{b})$&3&4&5&6&3&8&4&10&5&6&6&6&7&8\\ \hline
$\mmm(\hgbrx{2}{b})$&1&4&3&6&4&8&4&4&5&12&6&6&7&8\\ \hline
\end{tabular}
\end{center}
\begin{center}
\begin{tabular}{|l| llllllll |}
\hline b&15 &16&17&18&19&20&21&$b \ge 22$\\ 
\hline $\mmm(\hgbpx{2}{b})$&8 &8&9&10&10&10&10&10\\ \hline
$\mmm(\hgbrx{2}{b})$&8 &8&8&  8&  9&10&11&12\\ \hline
\end{tabular}
\end{center}
\caption{Values of $\mmm(\hgbpx{2}{b})$ and $\mmm(\hgbrx{2}{b})$.}\label{ttaabb}
\end{table}

\end{numlist}
\end{proclama}

\begin{proclama}{Theorem G}  \begin{numlist}  
\item \begin{eqnarray*}  
\mmm(\hgbrx{0}{b}) = \left\{\begin{array}{ll}
        1 & \mbox{if $b =1$,}\\
         \infty & \mbox{if $b=2$}, \\
         2 & \mbox{if $b =3$,}\\
        b-2 & \mbox{if $b \ge 4$.}\\
        \end{array}
        \right.
\end{eqnarray*}
\item $\mmm(\hgbrx{1}{b})=b-2$.
\item Table \ref{ttaabb} shows the values of $\mmm(\hgbrx{2}{b})$.
\end{numlist}

\end{proclama}

By Table \ref{ttaabb}, $\mmm(\hgbrx{2}{b}) \le 10$. We now explain roughly why this happens. (Precise arguments and definitions will be given later). Suppose that $f$ is an orientation preserving homeomorphism of $\Sigma_{2,b}$. If some iterate of the map $\Smap{\widetilde{f}}{\Sigma_2}$ has a fixed point, of index different from one, then the same iterated of the original $f$ has also a fixed point. On the other hand, it can be proved that for any $\Smap{k}{\Sigma_2}$, there exists a positive integer $n$ such that $n \le 10$ and $k^n$ has a fixed point of index different from one. So,  $\mmm(\hgbrx{2}{b}) \le 10$. That equality can be achieved whenever $b\ge 18$ is shown by means of examples. The situation is analogous for any genus larger than or equal to $2$, as is stated in Theorem H. To prove this, besides the ideas of \cite{N}, we also use the Thurston-Nielsen classification of surface homeomorphisms, and Nielsen fixed-point theory

\begin{proclama}{Theorem H} If $g \ge 2$ then $\mmm(\hgbrx{g}{b}) \le 4g+2$. Moreover, if $b \ge 6g+6$, then equality holds.
\end{proclama}

By Theorem H, the values of the minimum periods of orientation-preserving homeomorphisms  of surfaces of genus at least $2$ are bounded by a constant which does not depend on the number of boundary components. This situation is analogous for the classes of orientation reversing homeomorphisms as is stated in the following. 
\begin{proclama}{Theorem I} Let $g \ge 2$. Then $\mmm(\hgbrx{g}{b}) \le 4g+(-1)^g 4$ and equality holds if $b \ge 6g+2+(-1)^g8$. 
\end{proclama}

Obviously, these theorems do not cover all possible cases. Indeed, if $g \ge 3$, the values of $\mmm(\hgbp)$ are not given if $b < 6g+6$ and $b$ does not satisfy the conditions of Theorem D. Similarly, the values of $\mmm(\hgbr)$ are not given if $g \ge 3$, $b < 4g+(-4)^g$, and $b$ does not satisfy the conditions of Theorem E. the case $g=3$ could be solved by completely analogous methods to teh used in Theorems F and G. However, when $g \ge 4$ the quantity of variables makes the calculations too complicated. Also, it is not clear that $\mmm(\hgb)$ expressed as a simple function of $g$ and $b$.

This thesis is organized as follows. In Chapters \ref{S2} and \ref{S3} we give a summary of fixed-point theory and of the Thurston-Nielsen classification of homeomorphisms of surfaces, respectively. In Chapter \ref{S4} we present a standard form for such homeomorphisms. In Chapter \ref{S5} we present some features about planar discontinuous groups, and in Chapter \ref{S6} we apply these results about planar discontinuous groups to determine necessary  and sufficient conditions for the existence of certain finite-order maps of closed surfaces. Chapter \ref{S7} and \ref{S8} are devoted to developing the technical machinery which we will use in Chapters 9, 10 and 11 to prove our main results. 

These thesis has three main branches, which are interconnected. One has to do with the application of fixed-point theory described in Chapter \ref{S2}. All the upper bounds on the mimimum periods except the ones stated in Theorems H and I are consequences of this theory. To obtain the upper bounds of Theorems H and I we apply also the Thurston-Nielsen classification of homeomorphisms and some of its consequences, described in Chapters \ref{S3} and \ref{S4}. This is the second branch. Finally, then third branch has to do with the theory of planar discontinuous groups presented in Chapter \ref{S5}, which will provide us with the tools for constructing examples which will show the existence of lower bounds for minimum periods.

\part{Preliminary results}\label{rios}

\chapter{Fixed-point theory}\label{S2}

Let $X$ be a space and let \Smap{f}{X} be a self map.  Fixed-point theory studies the
nature of the set $\Fix(f)$ in relation to the space $X$ and the map $f$. This study can
be undertaken from different points of view. Since we are dealing with homeomorphisms of
topological manifolds, we focus our attention on the topological setting. Topological
fixed-point theory tries to answer concerning $\Fix(f)$, like what is  the cardinal of
this set, whether is it empty or not, or how does it change under homotopy.

Our aim is to study the existence of fixed points of iterates of continuous maps.
Consequently, fixed-point theory provides very useful tools which are described in this
chapter. More precisely, in Section~\ref{desierto resplandeciente} we give the definition
of the Lefschetz number and state the Lefschetz Fixed-Point Theorem; Section~\ref{it} is
devoted to index theory; and, in Section~~\ref{hecho}, we describe some particular
features of the Lefschetz number for homeomorphisms of surfaces.

\section{The Lefschetz Fixed-Point Theorem}\label{desierto resplandeciente}

Early in the history of fixed-point theory it was discovered that, if $X$ is a polyhedron
(see \cite{Bn} for a definition), and \Smap{f}{X} is a map with only a finite number of
fixed points satisfying an additional technical requirement, it is possible to associate
to each fixed point an integer, called the {\em index} \index{index} which describes the
way in which the map ``winds around'' the point. Furthermore, the sum of all indices was
found to be equal to the {\em Lefschetz number\/},\index{Lefschetz number}
\index{1l@$L(f)$} which is defined  for a continuous self-map $f$ of a polyhedron $X$
(or, more generally, a compact ANR) as the finite sum
\begin{equation}
L(f)= \sum (-1)^k\tr(f_{*k}),\label{literatura}
\end{equation}
where $\tr(f_{*k})$ denotes the trace of the map $f_{*k}$ induced by the action of $f$ on
the $k$-th rational homology group of $X$.

For every homeomorphism $\Smap{f}{X}$, every natural number $k$ and every integer $m$,
$(f_{*k})^m=(f^m)_{*k}$ (see \cite[III.3]{Do}), and we shall write $f_{*k}^m$ to denote
their common value.

\start{rem}{iso} Observe that if $\Smap{f,g}{X}$ are two homotopic maps then for every
natural number $k$, $f_{*k}=g_{*k}$. Therefore $L(f)=L(g)$.
\end{rem}

Perhaps the best known fixed-point theorem in topology is the Lefschetz Fixed-Point
Theorem.

\start{theo}{lef}{\em(Lefschetz (1923), Hopf (1929))} Let $X$ be a compact ANR and let
\Smap{f}{X} be a continuous map. If $L(f) \ne 0$ then every map homotopic to $f$ has a
fixed point.
\end{theo}

The first announcement of this theorem (for a restricted class of polyhedra) was in 1923
\cite{L1} and the details appeared in \cite{L2} and  \cite{L3}. The first proof of the
Lefschetz Fixed-Point Theorem for all polyhedra was given by Hopf \cite{Ho}. Also, a
particular case of this theorem valid for  $s$-to-$1$ maps of surfaces follows  from a
theorem proved by Alexander \cite{A} in 1923.

The Lefschetz number  counts the fixed points ``with multiplicity''.  It is a homotopy
invariant and is easily computable.  As is the case in all the other literature about
minimum periods, the Lefschetz Fixed-Point Theorem will be one of  our main tools.

\section{Index theory}\label{it}

The ``multiplicity'' of a fixed point (and, more generally, of  a fixed-point set which
is open in the whole set of fixed points of the map) is measured by the index. To study
this index we shall develop some theory, following mainly \cite[VII.5]{Do}.

Our first step will be to define the index of certain types of map, firstly for a map
from an open set of the Euclidean space $\R^n$ to $\R^n$; and subsequently replacing
$\R^n$ with any  ENR, (see the  definition below). In order to achieve the former goal,
we need to define a homological object: the fundamental class around a compact set.
Before doing this, we need to recall some notions from algebraic topology.

Let $X, Z$ be spaces, let $Y$ (resp. $W$) be a subspace of $X$ (resp. $Z$),  and let
$\map{f}{(X,Y)}{(Z,W)}$ be a map. The $k$-th integer homology group of the pair $(X,Y)$
is denoted  by $H_k(X,Y; \Z)$,  and the map induced by $f$, from  $H_k(X,Y; \Z)$ to
$H_k(Z,W; \Z)$, is denoted by $H_k(f,\Z)$. (Recall that the map induced by $f$ on
$H_k(X,\Q)$ is denoted by $f_{*k}$).

Let $n$ be a positive integer, let  $V$ be an open subset of $\R^n$, and let $K $ be a
compact subset of $V$. View  $\Sn$ as $\R^n \cup \{\infty\}$. Consider the inclusions
$$\map{i}{(\Sn,\emptyset) }{(\Sn, \Sn \setminus K)}$$ and $$\map{j}{(V, V \setminus K)
}{(\Sn, \Sn \setminus K)}.$$ These maps induce homomorphisms
$$\map{H_{n}(i,\Z)}{H_n(\Sn,\emptyset;\Z) }{H_n(\Sn, \Sn \setminus K;\Z)}$$ and
$$\map{H_{n}(j,\Z)}{H_n(V, V \setminus K;\Z)}{H_n(\Sn, \Sn \setminus K;\Z)}.$$ By the
Excision Lemma (see \cite[Corollary III.7.4]{Do}), $H_{n}(j,\Z)$ is an isomorphism.  On
the other hand, $H_n(\Sn,\emptyset; \Z)$  is isomorphic to $\Z$, so we can fix one of its
generators and denote it by $o$. We define  the {\em fundamental class around
$K$},\index{fundamental class around $K$} denoted by  $o_K$, \index{1o@$\o_K$} as the
element  $H_{n}(j,\Z)^{-1}H_{n}(i,\Z)(o)$ of $H_n(V, V \setminus K;\Z)$.

\begin{rem} $o_K$ is characterized by the property that its image under the map induced by the inclusion
$$\map{H_{n}(i,\Z)}{H_n(V, V \setminus K;\Z)}{H_n(V, V \setminus \{p\};\Z)}$$ agrees with
$o_{\{p\}}$ for every $p \in K$. Roughly speaking, it is an  element of $H_n(V, V
\setminus K;\Z)$ which bounds  $K$ taking into account its orientation.
\end{rem}

Let $V \subset \R^n$ be open and  consider a map $\map{f}{V}{\R^n}$ such that the set of
fixed points, $\Fix(f)$ is compact. Denote this set by $K$ and let the map
$$\map{i-f}{(V,V\setminus K)}{(\R^n,\R^n \setminus \{0\})}$$ be defined by
$(i-f)(x)=x-f(x)$. Consider the homomorphism $$ \map{H_{n}(i-f,\Z)}{H_n(V,V \setminus
K;\Z)}{H_n(\R^n,\R^n \setminus \{0\};\Z)}, $$ and define the {\em index of $f$}
\index{index of a map!of a subset of $\R^n$}  \index{1i@$\inn(f)$} as the integer
$\inn(f)$ such that $$ H_{n}(i-f,\Z)(o_K)=\inn(f) .o_{\{0\}}. $$ (Recall that $o_{\{0\}}$
generates $H_n(\R^n,\R^n \setminus \{0\};\Z)$). This definition does not depend on the
initial choice of the generator of $H_n(\Sn,\emptyset,\Z)$ because
$(-o)_{\{K\}}=-(o_{\{K\}})$ and $(-o)_{\{0\}}=-(o_{\{0\}})$.

Now we extend our definition of index to more general spaces, namely, Euclidean
neighborhood retracts. A topological space $Y$ is said to be a  {\em Euclidean
neighborhood retract} (or, briefly, an {\em ENR})  \index{Euclidean neighborhood retract}
\index{ENR} if a neighborhood retract $X \subset \R^n $ exists and is homeomorphic to
$Y$. The following proposition will allow us to define the index for maps of  ENR's (see
\cite[Proposition and Definition VII.5.10]{Do}).

\start{prop}{definition} If $Y$ is a topological space and $U$ is an open subset of $Y$
which is also an ENR, then every map $\map{f}{U}{Y}$ admits a decomposition $f= \beta
\alpha$ where $\map{\alpha}{U}{V}$, $\map{\beta}{V}{Y}$, and $V$ is open in some
Euclidean space $\R^n$.
\end{prop}

With the notation of the above proposition, consider the restrictions
$$\map{\beta|_{\beta^{-1}(U)}}{\beta^{-1}(U)}{U \cap \beta(V)}$$ and $$ \map{\alpha|_{U
\cap \beta(V)}}{U \cap \beta(V)}{V}.$$ If $\Fix(f)$ is compact we define the {\em index
of  $f$}, \index{index of map!of an ENR} denoted by $\inn(f)$,   \index{1i@$\inn(f)$} as
the index of $$\map{\alpha|_{U \cap \beta(V)} \circ
\beta|_{\beta^{-1}(U)}}{\beta^{-1}(U)}{V}.$$ The index defined in this way is independent
of the decomposition  (see \cite[Proposition and Definition VII.5.10]{Do}).

\start{rem}{pillow book} Observe that for every ENR $X$, and every map $\Smap{f}{X}$, if
$\Fix(f)=\emptyset$ then $H_k(V, V\setminus \Fix(f);\Z)$ is trivial for every natural
number $k$. Hence, $I(f)=0$.
\end{rem}

If  $K$ is an open compact subset of $\Fix(f)$ then we define the {\em index of K with
respect to $f$}, \index{index!of a compact set} denoted by $\ind{f}{K}$,
\index{1i@$\ind{f}{K}$} as the index of $f|_W$ where $W$ is an open subset of $U$ such
that $K=\Fix(f) \cap W$. This value is independent of the choice of $W$ (see
\cite[VII.5.11]{Do})).

In particular, if $x \in X$ is an isolated fixed point then the set $\{x\}$ is compact
and open in $\Fix(f)$,  so we can define its index, called the  {\em  index of $f$ at
$x$}  \index{index of a map!at a fixed point} and denoted by
$\ind{f}{x}$.\index{1i@$\ind{f}{x}$} There is  an alternative way of calculating the
index of a map at an isolated fixed point. Before describing it, we need to define
another important notion of algebraic topology, the degree of a map.

Let $M$ be an $n$-manifold and let $U \subset M$ be an open set. A map $\map{f}{U}{M}$ is
said to be {\em proper} \index{proper map} if $f^{-1}(K)$ is compact for every compact
set $K \subset M$. If the manifold $M$ is connected and orientable, and $\map{f}{U}{M}$
is a proper map, we define the {\em degree of $f$,} \index{degree of a map} denoted
$\deg(f)$, \index{1d@$\deg(f)$} as the integer  which satisfies $$
H_{n}(f,\Z)(o_{f^{-1}(K)})=\deg(f). o_K $$ for some non-empty compact subset $K$ of $M$.
It can be shown that $\deg(f)$ is independent of the choice of $K$ (see
\cite[VIII.4]{Do}).

Now, we can state the promised alternative definition of index of a fixed point.

\start{rem}{isolated} Let $V$ be an open subset of $\R^2$ such that $\Cl(\D) \subset V$.
Assume that $\map{f}{V}{\R^2}$ is a map such that the origin is the only fixed point of
$f$ in $\Cl(\D)$. Define $\Smap{\v_f}{\partial \D}$ by $\v_f(x)=\frac{x-f(x)}{|x-f(x)|}$.
It is a simple exercise in algebraic topology to show that $\ind{f}{x}$ equals the degree
of $\v_f$.
\end{rem}

There are other ways of defining the index of a fixed-point set although some of them
require strong restrictions on the class of maps for which the definitions apply.  In the
general case, it is possible to give an axiomatic definition, (see \cite{Bn}).

In the following remark, we give a property of the degree of a map of the circle.

\start{rem}{grado} If $\Smap{f}{\SI}$ is a map then $L(f)=1-\deg(f)$. In particular, if
$f$ is a homeomorphism,
\begin{eqnarray*}
L(f)&= &\left\{\begin{array}{ll}
    0& \mbox{if $f$ preserves orientation,}\\
    2 & \mbox{if $f$ reverses orientation.}
    \end{array}
    \right.
\end{eqnarray*}
\end{rem}

The following is a direct consequence of \cite[VII.5.13]{Do}.

\start{lem}{elolvidado} If $C \subset \Fix(f)$ is finite and open in $\Fix(f)$ then
$\ind{f}{C }$ is equal to the sum of the indices of the elements of $C$.
\end{lem}

As  mentioned above, the index of $f$ on the whole of $X$ equals $L(f)$, as is stated in
the next theorem  (see \cite[VII.6.13]{Do}).

\start{theo}{suma}   If $X$ is a compact ENR and  $\Smap{f}{X}$ is a map then
$$L(f)=\inn(f).$$ Furthermore, if  $U_i \subset X$,  $\ii{i}{k}$ are open subsets such
that $X=\bigcup\limits_{i=1}^{k}U_i$ and for each pair $i \ne j$, $U_i \cap U_j \cap
\Fix(f)=\emptyset$ then $$L(f)=\sum_{i=1}^{k}\inn(f|_{U_i}). $$ Consequently, if
$\Fix(f)$ is finite then $L(f)$ equals the sum of the indices of the fixed points.
\end{theo}

By \remc{pillow book}, the Lefschetz Fixed-Point \theoc{lef} is a consequence of
\theoc{suma}. Observe that the definition of $\inn(f)$ uses integer homology groups
whereas the definition of $L(f)$ uses rational homology groups.

For each homeomorphism $\Smap{f}{\Suu}$, unless we specify the contrary, when we speak
about the index of a fixed point of $f^i$ (or a pointwise fixed set), we mean the index
with respect to $f^i$.

In $\Fix(f)$ we define the following relation: Given $x, y \in \Fix(f)$, we say that  $x$
and $y$ are {\em $f$-equivalent\/}  \index{equivalence, $f$-} if there is a path $\alpha$
from $x$ to $y$ such that $\alpha$ and $f \circ \alpha$ are homotopic keeping the endpoints fixed. 
It is easy to see that the relation of $f$-equivalence is an equivalence relation.
The equivalence classes are called {\em fixed-point classes of $f$\/}. \index{fixed-point
classes} It can be proved (see \cite{Bn}, \cite{J}) that a fixed-point class is compact
and open in $\Fix(f)$, so its index is defined. A fixed-point class is called {\em
essential \/}\index{essential fixed-point class} if its index is different from $0$.
Essential fixed-point  classes will be important for us because they survive (preserving
their index) under isotopy, as  is stated in the following theorem, which is a corollary
of \cite[Theorem  VI.E.3]{Bn}.

\start{theo}{jiang} Let $i$ be an integer different from $0$, let $X$ be a compact
polyhedron, and let $\Smap{f}{X}$ be  a map. If $f$ has a fixed-point class of index $i$
and   $\Smap{g}{X}$  is homotopic to $f$, then $g$ has a fixed-point class of index $i$.
\end{theo}

Another important property of fixed-point classes is given in the following theorem
(see, for instance, \cite[Theorem I.4.3 and Theorem I.4.4]{J}).

\start{theo}{essential} The number of essential  fixed-point classes  is finite and the
sum of the indices of all (essential) fixed-point classes of $f$ equals $L(f)$.
\end{theo}

We end this section  by stating the following proposition (see \cite[Exercise
VII.6.25.2]{Do}).

\start{prop}{collapse} Let $X, A$ be compact ENR's such that $A \subset X$, and let
$\Smap{f}{(X,A)}$ be a map. Then $$ L(f)+1=L(\bar{f})+L(f|_A), $$ where
$\Smap{\bar{f}}{X/A}$ denotes the map induced by $f$ on the quotient space $X/A$.
\end{prop}

\section{Lefschetz numbers of maps of surfaces}\label{hecho}

If  $\Smap{f}{\Suu}$ is a homeomorphism, the Lefschetz number of $f$, $L(f)$, takes a
particular form. Here $\tr(f_{*0})=1$ and
\begin{eqnarray*}
\tr(f_{*2}) = \left\{\begin{array}{ll}
        1 & \mbox{if $f$ preserves  orientation and $\partial \Suu=\emptyset$,}\\
        -1 & \mbox{if $f$  reverses orientation and
          $\partial \Suu= \emptyset$}, \\
        0 & \mbox{if $\partial \Suu \ne \emptyset$,}\\
        \end{array}
        \right.
        \nonumber
\end{eqnarray*}
so \refeq{literatura} can be rewritten as
\begin{equation}
L(f) = \left\{\begin{array}{ll}
        2-\tr (\fs) & \mbox{if $f$ preserves  orientation and $\partial \Suu=\emptyset$,}\\
         -\tr (\fs) & \mbox{if $f$  reverses orientation and
          $\partial \Suu= \emptyset$}, \\
        1-\tr(\fs) & \mbox{if $\partial \Suu \ne \emptyset$.}\\
        \end{array}
        \right.
        \label{anne}
\end{equation}

Let $\lambda_1, \lambda_2, \ldots, \lambda_n$ be the eigenvalues of $\fs$. For each
positive integer  $i$, $L(f^i)=\tr(f_{*0}^i)+\tr(f_{*2}^i)-\sum_{j=1}^n\lambda_j^i$.
Therefore,  the polynomials $$p_i=p_i(\lambda_1, \lambda_2, \ldots,
\lambda_n)=\sum_{j=1}^n\lambda_j^i$$ play a key role in the computation of the numbers
$L(f^i)$. Indeed, by \refeq{anne},
\begin{equation}
L(f^i) = \left\{\begin{array}{ll}
        2-p_i & \mbox{if $f^i$ is orientation-preserving and $\partial \Suu=\emptyset$,}\\
                   -p_i & \mbox{if $f^i$  is orientation reversing and $\partial \Suu=\emptyset$,} \\
        1-p_i & \mbox{if $\partial \Suu \ne \emptyset$.}
                              \end{array}
        \right.
        \label{numerar}
\end{equation}

Clearly, the sequence of $p_i$'s determines the sequence of $L(f^i)$'s. For this reason,
we will study in more detail some of its properties.  Namely, we will study Newton's
equations.

Let $k$ be a positive integer and let $M \in GL_k(\Z)$. Denote by $\lambda_1, \lambda_2,
\dots, \lambda_k$ the eigenvalues of $M$. If we write the characteristic polynomial of
$M$ as \index{characteristic polynomial!of a matrix}

$$P(x)=\d(xI_k-M)=x^k+s_1x^{k-1}+s_2x^{k-2}+ \cdots+s_{k-1}x+s_k,$$ then, for each  $j
\in \{1,2,\dots,k\}$, $$s_j\, =\, (-1)^j \sum \limits_{1\leq i_1<\cdots<i_j\leq k}
\lambda_{i_1}\cdots\lambda_{i_j}, \mbox{and}$$ \index{Newton's equations}\label{newton}

\renewcommand{\theequation}{N.\arabic{equation}}
\setcounter{equation}{0}
\begin{eqnarray}
p_1+s_1 &=& 0, \label{N1}\\ p_2+s_1 p_1+2s_2&=& 0,\label{N2}\\ \vdots & & \nonumber
\end{eqnarray}
\renewcommand{\theequation}{N.k}
\begin{equation}
p_k+s_1 p_{k-1}+s_2 p_{k-2}+\dots + k s_k = 0,\mbox{ and }\label{Nk}
\end{equation}
\renewcommand{\theequation}{N.(k+l)}
\begin{equation}
p_{k+l}+s_1 p_{k+l-1}+s_2 p_{k+l-2}+\dots + s_k p_l =0 \mbox{ for $l \ge 1$.}
\label{Nk+l}
\end{equation}
\renewcommand\theequation{\thechapter.\arabic{equation}}
\setcounter{equation}{3}

\noindent See,  for example, \cite[Exercise 2, Section 29]{V}.

If $\Smap{f}{\Suu}$ is a homeomorphism then the characteristic polynomial of $\fs$ will
be called the {\em characteristic polynomial  of $f$}.\index{characteristic polynomial!of
a map}

Now, in order to study one of the properties of $\fs$, we shall introduce more notation.
For each positive integer $g$ denote by $J_g$ the element  \index{1j@$J_g$} $$ \left
(\begin{array}{ll}
         0&I_g\\
         -I_g & 0
    \end{array}
    \right )
$$ of $GL_{2g}(\Z)$. We say that $M \in GL_{2g}(\Z)$ is {\em proper symplectic} (resp.
{\em improper symplectic}) \index{proper symplectic matrix} \index{improper symplectic
matrix} if $M^t J_g M=  J_g$ (resp. $M^t J_g M=  -J_g$). A matrix $M$ is {\em symplectic}
\index{symplectic matrix} if it is either proper symplectic or improper symplectic.

A property of proper simplectic matrices is stated in the following; see \cite{Si}.

\start{prop}{gracias} If $M \in GL_{2g}(\Z)$ is proper symplectic then $\d(M)=1$.
\end{prop}

The next result  is a consequence of  \cite[Theorem 3.6.7]{Zie}.

\start{prop}{efectos} With respect to a certain basis, for every $f \in \hgp$ (resp.
$\hgr$), the matrix of $\fs$ is  proper symplectic (resp. improper symplectic).
\end{prop}

\chapter{The Thurston-Nielsen classification}\label{S3}

The goal of this chapter is to give  a brief introduction to the Thurston-Nielsen
classification of isotopy classes of surface homeomorphisms. This is, undoubtedly, the
most important tool in the topological theory of surface dynamics. It can be viewed as a
prime decomposition theorem: it gives the existence in each isotopy class of a
homeomorphism that is constructed by gluing together homeomorphisms of two types,
pseudo-Anosov and finite-order. The theory has numerous applications and implications for
many diverse areas of mathematics,  but we will focus on some of its dynamical aspects.
The reader is referred to \cite{T} for the original proof, to \cite{FLP} or \cite{HT} for
proofs of the theorem for orientable surfaces, and to \cite{Wu} for a proof for
non-orientable surfaces. Also, an algorithmic proof 
can be found in \cite{BH}.

This chapter is organized in the following way: Sections~\ref{fo} and  \refc{pam} are
devoted to finite-order maps and  to pseudo-Anosov maps respectively, and, in
Section~\ref{tct} we define reducible maps, give the Thurston-Nielsen  classification and
state some properties of finite-order maps.

\section{Finite-order maps}\label{fo}

The simplest types of maps used in the construction of a Thurston representative are the
finite-order maps defined in Chapter~\ref{S1}.  They are dynamically very simple: the
period of each orbit equals the order $n$ of the map, except for a finite number of orbits
whose period is a divisor of $n$. If $\Smap{f}{\Suu}$ is an isometry with respect to  a
hyperbolic metric, then it is standard that $f$ is finite-order; see \cite[Expos\'{e} 3,
Th\'{e}oreme 18]{FLP}. Conversely, when $f$ is finite-order on a surface of negative Euler
characteristic, it is an isometry with respect to some hyperbolic metric;  see
\cite[Theorem 2.8]{E}.

If $\Smap{f}{\Suu}$ is a homeomorphism, then the {\em order of $f$}, \index{order!of a
finite-order map} denoted by $\sf$, \index{1s@$\sf$} is the least positive integer
$\sigma$ such that $f^\sigma=\id_{_{\Suu}}$, or $\infty$ if no such $\sigma$ exists. The
{\em order} of a non-empty class $\CC$ of  homeomorphisms \index{order!of a class of
finite-order homeomorphism} is defined to be the supremum of  the $\sf$, $f \in \CC$.
Wiman \cite{Wi} determined the order of $\fgp$  and Wang \cite{W1} determined the order
of $\fgr$; see also \cite{Ha}. These results are summarized in the following.

\start{theo}{www} If $g \ge 2$ then the order of  $\fgp$ is $4g+2$ and the order of
$\fgr$ is $4g+(-1)^g4$.
\end{theo}

\start{rem}{soniar} Observe that the order of $\fgp$ (resp. $\fgr$) coincides with the
upper bound given for $\mmm(\hgbp)$ (resp. $\mmm(\hgbr)$) by Theorem~H (resp. Theorem~I).
\end{rem}

As well as the order of finite-order maps, we require some information about the
connected components of $\Fix(f)$. This is stated in the following lemma which is a
consequence of \cite[Lemma 1.1]{JG}; see also \cite[Theorem 2.8]{E}.

\start{lem}{jiji} If $\Smap{f}{\Suu}$ is an orientation-preserving finite-order
homeomorphism, and $A$ is a connected component of $\Fix(f)$, then either $A=\Suu$, or
$A$ contains only one point and it has  a neighborhood homeomorphic to $\D$  where $f$
acts as a rotation. In the former case,  $\ind{f}{A}=\chi(\Suu)$, and in the latter case,
$\ind{f}{A}=1$.
\end{lem}

Observe that the trace of  $(\id_{_{\Suu}})_{*k}$  equals the $k$-th Betti number of
$\Suu$, i.e., $\tr ((\id_{_{\Suu}})_{*k})=b_k(\Suu)$ for each $k \in \N$.  In particular,
by  \refeq{anne}, $L(\id_{_{\Suu}})=\chi(\Suu)$. By \lemc{jiji} and \theoc{suma}, we have
the following.

\start{lem}{nw} If $\Smap{f}{\Suu}$ is an orientation-preserving finite-order map then

\begin{equation}
L(f) = \left\{\begin{array}{ll}
        \Card(\Fix(f)) & \mbox{if $f \ne \id_{_{\Suu}}$,}\\
                   \chi(\Suu) & \mbox{if $f=\id_{_{\Suu}}$.}
                              \end{array}
        \right.
        \nonumber
\end{equation}
\end{lem}

Now we state a result analogous to \lemc{jiji} for the orientation-reversing case; see
\cite[Theorem 2.8]{E} for the proof.

\start{lem}{jijior} Let $f \in \fgr$ and $A$ be a connected component of $\Fix(f)$. Then
$A$ is a simple closed curve, with a neighborhood $U$ homeomorphic to $\SI \times (0,1)$
where $f$ acts as the reflection $(z,t) \mapsto (z, 1-t)$. Moreover, $\ind{f}{A}=0$.
\end{lem}

Applying Lemmas~\ref{jiji} and \refc{jijior} we obtain the following.

\start{lem}{parageo} Let $f \in \fgr$. If there exists $\ii{i}{\sf-1}$ such that
$f^i|_C=\id_C$ for some simple closed curve $C \subset \Su{g}$, then $\sf \equiv_4 2$
and $i=\sf/2$.
\end{lem}

\section{Pseudo-Anosov maps}\label{pam}

A detailed description of pseudo-Anosov maps of surfaces without boundary can be found in
\cite{KG}, but this author could not find a good description for  pseudo-Anosov maps of
surfaces with boundary. For this reason we give a complete definition here, although  for
our purposes it would be enough to give the description of a standard form defined in
Chapter~\ref{S4}.

A {\em singular foliation $F$ of $\Sigma$} \index{singular foliation} is a partition of
$\Suu\setminus \{x_1, x_2,\ldots,x_m\}$, for some finite subset $ \{x_1,
x_2,\ldots,x_m\}$ of $\Suu$,  into a disjoint union of one-dimensional manifolds, called
{\em leaves}, \index{leaves} such that there exists a  finite  $C^\infty$ atlas with
charts  $$(\varphi_i,U_i)_{1 \le i \le l},$$ where $\map{\varphi_i}{U_i}{\C}$, with $m
\le l$ and $\cup_{1 \le i \le l}U_i=\Suu$,   and a finite sequence of integers  $p_i \ge
3$, $i \in \{1, 2, \ldots, m\}$, and  $s \in \{0,1,2,\dots,m\}$, such that the following
statements hold:

\begin{case}{1}If $1 \le i \le s$ then
\end{case}

\begin{numlist}
\item $\varphi_i(U_i)=D_{a_i}  \cap \{z \in \C  \,\,:\, \Re(z) \ge 0\}$ for some $a_i >0$;
\item $x_i \in U_i$ and $\varphi_i(x_i)=0$;
\item For each leaf $L$ of $F$, if $K$ is a component of  $L \cap U_i$, then there exists $k \in [0,\infty)$ such that $K$ is mapped bijectively to a component of $\{z \in \C  \,\,:\,  \Im(z^{p_i-1})=k \} \cap \varphi_i(U_i).$
\end{numlist}

\begin{case}{2}If $s+1 \le i \le m$ then
\end{case}
\begin{numlist}
\item $\varphi_i(U_i)=D_{a_i}$ for some $a_i >0$;
\item $x_i \in U_i$ and $\varphi_i(x_i)=0$;
\item For each leaf $L$ of $F$, if $K$ is a component of  $L \cap U_i$, then there exists $k \in [0,\infty)$ such that $K$ is mapped bijectively to a component of $\{z \in \C  \,\,:\, \Im(z^{p_i/2})=k \} \cap \varphi_i(U_i).$
\end{numlist}

\begin{case}{3}If $m+1 \le i \le l$ then
\end{case}
\begin{numlist}
\item $\varphi_i(U_i)=(0,b_i)\times (0,c_i)$ or $(0,b_i)\times [0,c_i)$ for some $b_i ,c_i>0$;
\item   For each leaf $L$ of $F$, if $K$ is a component of  $L \cap U_i$, then there exists $k \in [0,\infty)$ such that $K$ is mapped bijectively to a component of $\{z \in \C  \,\,:\,  \Im (z) =k \} \cap \varphi_i(U_i).$
\end{numlist}

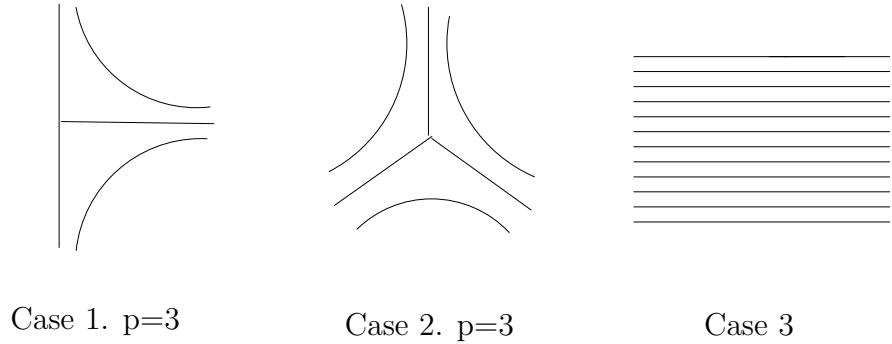
\placedrawing{pseudo.lp}{Charts for singular foliations}{figfol}

The points $x_i$,  $i \in \{1, 2, \ldots, m\}$, are called the {\em singularities}
\index{singularities} of the foliation and the points in $\Suu \setminus\{x_1,x_2,
\ldots, x_m\}$ are called {\em regular points}. \index{regular points} Observe that each
boundary component of $\Suu$  is a finite union of leaves and singular points.

For each $i \in \{1,2, \ldots, m\}$ there are $p_i$ leaves emanating from $x_i$. These
leaves are called {\em prongs}  \index{prongs} and we refer to $x_i$ as a {\em
$p$-pronged singularity} \index{ppronged sin@$p$-pronged!singularity}. If $x$ is a
regular point then it is contained in a single leaf $L$. However for convenience, we will
refer to the oriented components of $L \setminus \{x\}$ as the {\em prongs emanating from
$x$}  \index{prongs!emanating from a regular point} and to $x$ as a {\em $2$-pronged
point}.  \index{two@$2$-pronged point} A leaf emanating from a boundary component $B$ but
not contained in $B$ is called  a {\em prong of $B$} \index{prong!of a boundary
component} and $B$ is called a {\em $p$-pronged boundary component}  \index{p
pronged@$p$-pronged!boundary component} if there are exactly $p$ prongs emanating from
it.

Let $x$ be a  singularity of a singular foliation $F$ and let $(\varphi_x,U_x)$ be a
chart such that $x \in U_x$. We say that an arc $\alpha \subset \Suu$ such that $x \in
\alpha$ is {\em transverse to $F$ in $x$} \index{arc transverse!to a foliation in a
point} if there exists $r>0$ such that for every $\overline{r}<r$,  $(\alpha \setminus
\{x\}) \cap \varphi_x^{-1}(D_{\overline{r}})$ has exactly two  connected components and
each  of them is included in one connected component of
$\varphi_x^{-1}(D_{\overline{r}}\cap (\varphi_x(U_x)\setminus \{z \in \C  \,\,:\, \Im
(z^{p/2})=0\}))$; see \figc{labarctrans}.

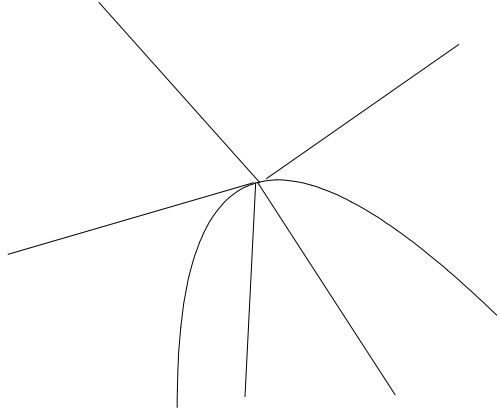
\placedrawing{arctrans.lp}{An arc transverse to a foliation in a
singularity}{labarctrans}

An arc $\alpha$ is  {\em transverse to a foliation $F$} \index{arc transverse!to a
foliation $F$} if it is transverse to the  leaves of $F$ in the usual sense in the
regular points of $F$ and transverse to $F$ in the singularities.

A {\em transverse measure $\mu$ to a foliation with singularities $F$} \index{transverse
measure to a foliation with singularities} is a map which assigns to each  arc $\alpha$
transverse to $F$ a non-negative Borel measure $\mu_\alpha$ on $\alpha$,   with the
following properties:
\begin{numlist}
\item If $\beta$ is a subarc of $\alpha$ then the measure $\mu_\beta$ is the restriction of $\mu_\alpha$ to $\beta$.
\item If  $\alpha_0$, $\alpha_1$  are two arcs in $\Suu$ transverse to $F$ related by a homotopy $\map{\alpha}{[0,1]\times[0,1]}{\Suu}$ such that $\alpha([0,1]\times \{0\})=\alpha_0$ and $\alpha([0,1]\times \{1\})=\alpha_1$ and
$\alpha( [0,1] \times \{a\})$ is contained in a leaf of $F$ for each $a \in [0,1]$, then
$\mu_{\alpha_0}=\mu_{\alpha_1}$.
\end{numlist}

If $\Smap{k}{\Suu}$ is a homeomorphism and $F$ is a foliation of $\Suu$ we say that $F$
is {\em invariant under  $k$}  \index{invariant foliation under  a map} if the images
under $k$ of the leaves of $F$ are leaves of $F$. If $F$ is a foliation invariant under
$k$ and $\mu$ is a measure transverse to $F$, we define the {\em image measure}
\index{image measure} $k(\mu)$ as the transverse measure to $F$ such that if $\alpha$ and
$\beta$ are arcs transverse to $F$ then
$k(\mu)_\beta(\gamma)=\mu_{k^{-1}(\beta)}(k^{-1}(\gamma))$ for every Borel set $\gamma
\subset \alpha$.

A pair $(F,\mu)$ is a {\em measured foliation} \index{measured foliation} if $F$ is a
singular foliation and $\mu$ is a transverse measure to $F$.

Let $F_{1}$ and $F_{2}$ be two singular foliations and assume that $x$ is an interior
$p$-pronged singularity for $F_{1}$ and $F_{2}$.  We say that  $F_{1}$ and $F_{2}$ are
{\em transverse in $x$} \index{transverse foliations!at a point} if there exists a
neighborhood $U$ of $x$, $r >0$, and a homeomorphism $\map{\v}{U}{D_r}$ such that

 \begin{numlist}
\item $\v(x)=0$;
\item For each leaf $L$ of $F_1$, if $K$ is a component of  $L \cap U$, then there exists $k \in [0,\infty)$ such that $K$ is mapped bijectively to a component of $\{z \in \C  \,\,:\,   \Re(z^{p/2})=k \} \cap \v(U)$;  see \figc{labtransfol}
\item For each leaf $L$ of $F_2$, if $K$ is a component of  $L \cap U$, then there exists $k \in [0,\infty)$ such that $K$ is mapped bijectively to a component of $\{z \in \C  \,\,:\,   \Im(z^{p/2})=k \} \cap \v(U)$;  see \figc{labtransfol}
\end{numlist}

Let $B$ be a boundary component of $\Suu$, and  let  $F_1$ and $F_{2}$ be two singular
foliations on $\Suu$. Collapsing $B$ to a point $x$, the foliations $F_{1}$ and $F_{2}$
induce two (not necessarily singular) foliations $\overline{F_{1}}$ and $\overline{F_{2}}$ respectively in
the resultant surface $\overline{\Suu}$. We say that $F_{1}$ and $F_{2}$ are {\em
transverse in $B$} \index{transverse foliations!at a boundary component} if
$\overline{F_{1}}$ and $\overline{F_{2}}$ are transverse in   $x$;  (see
\figc{labtransfol}). Two singular foliations are {\em transverse } \index{transverse
foliations} if they have the same interior singularities, and they are transverse in
interior singularities, boundary components and regular points (see \figc{labtransfol}).

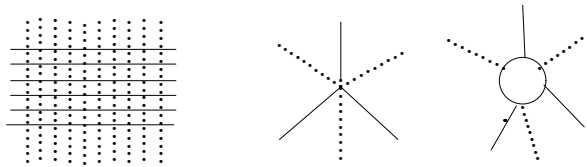
\placedrawing{transv.lp}{Charts for transverse singular foliations}{labtransfol}

A map $\Smap{f}{\Suu}$ is called {\em pseudo-Anosov}\index{pseudo-Anosov! homeomorphism}
if there exist two measured transverse foliations $(F^s, \mu^s)$ and $(F^u, \mu^u)$ which
are invariant under $f$, and such that each boundary component contains at least a
singular point,  and  a real constant $\lambda > 1$ such that
$f(\mu^s)=\lambda^{-1}\mu^s$ and $f(\mu^u)=\lambda \mu^u$. (This  is usually denoted as
$f(F^s,\mu^s)=(F^s, \lambda^{-1} \mu^s)$ and $f(F^u,\mu^u)=(F^u, \lambda \mu^u)$).  This
number $\lambda$ is called the {\em expansion constant for $f$}.\index{expansion constant
for a pseudo-Anosov homeomorphism}

\start{rem}{lambda}  If $\alpha \subset \Suu$ is an arc included in a leaf of $F^u$
(resp. $F^s$) and it is not included in the boundary of $\Suu$ then $\alpha$ is
transverse to $F^s$ (resp. $F^u$). Thus, $\mu^s(\alpha)$ (resp. $\mu^u(\alpha)$) is
defined. Furthermore,  $\mu^s(f(\alpha))=\lambda^{-1} \mu^s(\alpha)$ (resp.
$\mu^u(f(\alpha))=\lambda\mu^u(\alpha)$).

On the other hand, if $\alpha$ is an arc included in $\partial \Suu$ then neither
$\mu^s(\alpha)$ nor $\mu^u(\alpha)$ are defined.
\end{rem}

>From now on we will  work with only  the unstable foliation $F^u$. For simplicity, we
will refer to it as {\em the\/} foliation.

Denote by  $\Sing(\Suu)$  the set of singularities of the foliation. A   very useful
equation which relates the Euler characteristic of $\Suu$ with the foliation is the {\em
Euler-Poincar\'{e} Formula}, \index{Euler-Poincar\'{e} Formula}
\begin{equation}
\sum_{s \in \Sing(\Suu)} (2-p_s)=2 \chi(\Suu),\label{ep}
\end{equation}
where $p_s$ is the number of prongs emanating from $s$ for each $s \in \Sing(\Suu)$; see
\cite[expos\'{e} 5]{FLP} for a proof.

\start{rem}{pprong} If $B$ is a $p$-pronged boundary component, then $$\sum_{s \in
\Sing(B)} (2-p_s)=-p.$$
\end{rem}

\section{The classification theorem}\label{tct}

A {\em system of reducing curves for a surface $\Suu$}  \index{system of reducing
curves!for a surface} is a finite  (possibly empty) set of pairwise disjoint simple
closed curves $\Gamma=\{\Gamma_1,\Gamma_2,\ldots,\Gamma_n\} \subset \Suu$ such that each
connected component of $\Suu \setminus \Gamma$ has negative Euler characteristic. If
\Smap{f}{\Suu} is a homeomorphism, then a {\em system of reducing curves for $f$},
\index{system of reducing curves!for a homeomorphism} or an {\em $f$-system of reducing
curves}, is a system of reducing curves $\Gamma$ for $\Suu$ which is $f$-invariant and
has an $f$-invariant tubular neighborhood  $N(\Gamma)$ of $\Gamma$, called an {\em
$f$-tubular neighborhood.} A homeomorphism $\Smap{f}{\Suu}$ is said to be {\em reducible}
\index{reducible homeomorphism} if there exists a non-empty system of reducible curves
$\Gamma$ for $f$  and, for each connected component $R$ of $\Suu \setminus N(\Gamma)$,
there exists  a positive integer $n$ such that:
\begin{numlist}
\item $f^n(R)=R$.
\item $\Smap{f^n|_R}{R}$ is either finite-order or pseudo-Anosov.
\end{numlist}
The subsurface $R$ is called an {\em $f$-component} \index{component, $f$-} or {\em
component of $f$}.  \index{component of a homeomorphism}The least positive integer $n$
which satisfies $f^n(R)=R$ is called the {\em period of $R$} \index{period!of a
component} or {\em $f$-period of $R$}. \index{period of $R$, $f$-} We say that a $R$ is
a {\em pseudo-Anosov component of $f$}  \index{pseudo-Anosov!component of a
homeomorphism} or a  {\em finite-order component of $f$}  \index{finite-order component
of a homeomorphism} if  the  homeomorphism $\Smap{f^n|_R}{R}$ is pseudo-Anosov or
finite-order, respectively.

\start{theo}{th}{\em (Thurston-Nielsen)} If $\Smap{f}{\Suu}$ is a
homeomorphism then $f$ is isotopic to a homeomorphism $f'$ which is finite-order,
pseudo-Anosov or reducible.
\end{theo}

Here we say that  $f'$ is in {\em Thurston canonical form}. \index{Thurston canonical
form} Notice that $f'$ is not uniquely determined by $f$. In general, one says that an
isotopy class  is {\em finite-order}, \index{finite-order class}{\em pseudo-Anosov}
\index{pseudo-Anosov!class} or {\em reducible}  \index{reducible class} if  an element in
Thurston canonical form has the corresponding property. If  an isotopy class is
finite-order then any complicated behavior of the maps in this class  can be isotoped
away. However, complicated behavior of maps in a pseudo-Anosov class persists under
isotopy.

\chapter{The standard form}\label{S4}

Any  isotopy class contains infinitely many homeomorphisms as  Thurston canonical forms.
However, while finite-order isotopy classes have, roughly speaking, a unique
representative, for pseudo-Anosov homeomorphisms of surfaces with boundary,  there  is a
certain amount of choice involved in the structure of the foliation and the dynamics of
the representatives of a given class.  There are a number of papers that have used this
freedom to refine the Thurston canonical  form for dynamical purposes. Such a  refinement
was used to prove the existence of a dynamically minimal representative for Nielsen
classes of fixed points in the category of surface homeomorphisms. This result was
sketched by Jiang in \cite{J2} and given in full detail by Jiang and Guo in \cite{JG}.
They isotope the Thurston canonical form in different steps. We are not interested in the
final refinement but in a specific type of Thurston canonical form obtained during the
process. They call this particular type {\em the standard form}.  \index{standard form!of
a map}

In order to prove our main results we need to prove the existence of fixed-points classes
of index different from one  for iterates of a homeomorphism. Since essential fixed-point
classes (and their indices) are preserved under isotopy, it will be enough to prove that
each map in an isotopy class has a fixed-point class of index $0$, or that one
representative of the isotopy class has a fixed-point class of negative index.  We will
do this for maps in standard form because this class of  maps has very useful properties:
iterates of a map in standard form are also in standard form, fixed-point classes are
connected,  their structure is well understood, and their indices are easily calculable.

Throughout  this chapter we will follow \cite{JG}, in order to give  a description of the
standard form and to state some of its properties. This description will be done in the
spirit of the Thurston-Nielsen Theorem: In Section~\ref{ppaa} we describe the standard
form for pseudo-Anosov maps and its fixed-point classes; in  Section~\ref{maquinista} we
define a reducible map in  standard form to be a reducible map in Thurston canonical
form.  The main result of this chapter is \propc{types}, where we describe the structure
of the fixed-point classes of a map in standard form.

\section{Standard form for pseudo-Anosov maps}\label{ppaa}

In view of the need to glue finite-order pieces and pseudo-Anosov pieces together, the
standard form  is required to be finite-order on $\partial \Suu$. By  definition, a map
in  standard form  is  smooth  but we are not going to take this property into account
here, because  it will not be necessary  for our purposes.

We begin by defining some auxiliary maps. Fix a real number $\lambda>1$ and a natural
number $p$.  Consider  the correspondence ${\Phi_p}\colon{\C \setminus\{0\}} \rightarrow
{\C \setminus\{0\}}$, $z \mapsto z^{p/2}$, and the diffeomorphism  $\map{\Psi}{\C
\setminus \D}{\C \setminus \{0\}}$ defined as $$\Psi(z)=z-\frac{z}{|z|}.$$ Let
$\overline{V}$ be  the vector field on $\C$  obtained by ``slowing down'' the vector
field $$V(s)=\lambda \overline{s},$$ with a non-decreasing smooth function
$\map{\alpha}{\R}{[0,1]}$ such that $\alpha(0)=0$ and $\alpha(1/\lambda^2)=1$; that is,
$\overline{V}(s)=\alpha(|s|)V(s)$.

Recall that the time-one map \index{time-one map} associated to a vector field $W$ is
defined as $\v(x,1)$, where $\v(x,t)$ is the solution of the equation
\begin{eqnarray*}
 \cases{\frac{d\v(x,t)}{dt}=W(\v(x,t))\cr
                \v(x,0)=x,}\cr
\end{eqnarray*}
and recall that the vector field $\Phi_*(\overline{V})$ is defined in a point $q$ as the
product of the matrix $d\Phi^{-1}_{\phi(q)}$ with the vector $\overline{V}(q)$.

Let $\Smap{f}{\C}$ and $\Smap{f'}{\C \setminus\D}$ be  the time-one maps associated to
the vector fields
\begin{eqnarray*}
\campo_p(z) & = & \cases{\Phi_{p*}^{-1}\overline{V}(z) & \mbox{if  $z \ne 0$},\cr
                                            0 & \mbox{if $z=0$},}\cr
\end{eqnarray*}

\begin{eqnarray*}
\campo'_p(z) & = & \cases{\Psi_*^{-1}\campo_p(z) & \mbox{if  $z \in \C \setminus
\Cl(\D)$},\cr
                                            0 & \mbox{if $z \in \partial \D$}.}\cr
\end{eqnarray*}
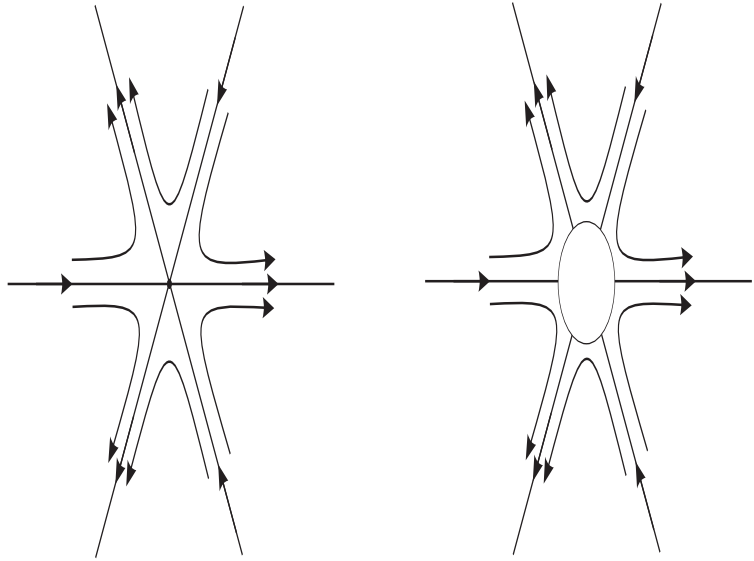
See \figc{labflows}. \placedrawing{cstand.lp}{The flows of $v_p$ and $v'_p$
($p=3$).}{labflows}

Now, consider  $\alpha \in \R$ and define $\Smap{r_\alpha^{+}}{\C}$  and
$\Smap{r_\alpha^{-}}{\C}$ by $$ r^{+}_{\alpha}(z)=z e^{ 2i\alpha\pi} \mbox{ and }
r^{-}_{\alpha}(z)=\overline{z} e^{-2i\alpha\pi}. $$

We say that a pseudo-Anosov map $\Smap{f}{\Sigma}$ is in {\em standard form} if there is
a finite smooth atlas $U$ of $\Suu$, consisting of  one chart for each interior
singularity, one chart for each boundary component and charts at interior regular points
such that: \index{standard form!of a pseudo-Anosov map}
\begin{numlist}
\item If $x$ is an interior $p$-pronged point (possibly regular)  and  $(U_x,u_x)$ is the chart for $x$ then the measures $\mu^s$ and $\mu^u$ are mapped by $u_x$ to the measures $|\Re\,\,  d \Phi_p(z)|$ and $|\Im \,\,d \Phi_p(z)|$ on $\C$ respectively. The leaves of $F^s$ and $F^u$ are mapped to the lines $\{z \in \C \,\,:\,\,\Re \Phi_p(z)=\mbox{constant}\}$ and $\{z \in \C \,\,:\,\, \Im \Phi_p(z)=\mbox{constant} \}$ respectively.
\item If $A$ is a $p$-pronged boundary component and $(U_A, u_A)$ is its chart, then the measures $\mu^s$ and $\mu^u$ on $U_A$ get mapped to the measures
$|\Re\,\, d \Phi_p \Psi(z)|$ and $|\Im\,\, d \Psi \Phi_p(z)|$ on $\C \setminus \D$
respectively. The leaves of $F_s$ and $F_u$ get mapped to the lines $\{z \in
\C\,\,:\,\,\Re \Phi_p\Psi(z)=\mbox{constant}\}$ and $\{\Im
\Phi_p\Psi(z)=\mbox{constant}\}$ respectively.
\item For each chart $(U_x, u_x)$ at an interior singularity $x$, $u_x(U_x)$ contains the closed disk $\Cl(\D)$. Moreover $\Cl(\D) \cap U_y =\emptyset$ for every $y \ne x$ such that $(U_y, u_i)$ is a chart of the atlas. Similarly for the closed annulus $\{z \,\,:\,\,1 \le |z|\le2\}$ in each chart $(U_A, u_A)$.
\item If $x$ is a $p$-pronged point or $A$ is a $p$-pronged boundary component then
for each chart $(U_x,u_x)$ and $(U_A,u_A)$, there exists  $0 \le k <p$ such that the
following diagrams commute.
\end{numlist}
$$
\begin{CD}
(U_x,\{ x\}) @>{f}>>  (U_{\v(x)},\{\v(x)\})     \\ @VV{ u_x}V        @VV{u_{\v(x)}}V \\
(\C, \{0\}) @>r^\epsilon_{_{{k/p}}} f>> (\C,\{0\})
\end{CD}
\hspace{1cm}
\begin{CD}
(U_A, A)    @>{f}>>(U_{\v(A)},\v(A))\\ @VV{u_A}V        @VV{u_{\v(A)}}V\\ (\C \setminus
\D, \{0\}) @>r^\epsilon_{_{{k/p}}}   f'>> (\C \setminus \D,\{0\})
\end{CD}
$$ where $\epsilon$ is $+$ or $-$ when the map $f$ is orientation-preserving or
reversing, respectively.

We say that a fixed point $x$ or an invariant boundary component $A$ is {\em  of type
$(p,k)^+$} \index{type of!a fixed point of a pseudo-Anosov map}\index{type of!an
invariant boundary component of a pseudo-Anosov map} (resp. {\em of type $(p,k)^-$}) if
$f$ preserves orientation (resp. reverses orientation) and  one of the following
statements holds:
\begin{numlist}
\item $x$ is a singularity or $A$ is a boundary component and $(p,k)^+$ (resp.   $(p,k)^-$) is as in part $(4)$ of the definition of standard form.
\item $x$ is a regular point, $(p,k)=(2,0)$  and the two prongs in $F^u$ emanating from $x$ remain fixed under the action of $f$.
\item $x$ is a regular point, $(p,k)=(2,1)$  and the two prongs in $F^u$ emanating from $x$ are interchanged by  $f$.
\end{numlist}

Roughly speaking, we can say that a pseudo-Anosov map acts as a rotation of angle
$2k\pi/p$  on the prongs of a fixed point of type $(p,k)^+$ and as a reflection on the
prongs of a fixed point of type $(p,k)^-$. So we have the following  remark.

\start{prop}{dash} Let $x$ (resp. A) be a fixed point (resp. an invariant boundary
component) of type $(p,k)^+$ of  an orientation-preserving  pseudo-Anosov homeomorphism
$f$.  Then, for each positive integer $n$,  $x$ (resp. $A$),
 is of type $(p,\widetilde{nk}^p)^+$  for $f^n$. In particular, $x$ (resp. $A$) is of type $(p,0)^+$ for $f^p$. In other words, all prongs emanating from $x$  (resp. $A$) remain fixed under the action of $f^p$.

On the other hand, if $x$  (resp. $A$) is a fixed point  (resp. an invariant boundary
component) of an orientation-reversing pseudo-Anosov map $f$ of type $(p,k)^-$  then $x$
(resp. $A$) is of type $(p,0)^+$ for $f^2$. Hence, all prongs emanating from $x$  (resp.
$A$) remain fixed under the action of $f^2$.
\end{prop}

\start{prop}{pa} Let $\Smap{f}{\Suu}$ be a pseudo-Anosov homeomorphism in  standard form.
\begin{numlist}
\item If $f$ preserves orientation then the fixed-point classes of $f$ are either  interior fixed points or invariant boundary components.
\item If $f$ reverses orientation then the fixed-point classes of $f$ are fixed points.
\end{numlist}
Furthermore, the indices of fixed-point classes depend on their types as described in
Tables ~\ref{indicespA} and \ref{indicesbdpA}.
\end{prop}
\begin{table}
\begin{center}
\begin{tabular}{|l| lllll|}
\hline Type of $x$ &  $(p,0)^+$ &  $(p,k)^+$ & $(2n,0)^-$& $(2n,1)^-$ & $(2n+1,0)^-$\\
\hline $\ind{f}{x}$ & $1-p$ & $1$ & $-1$ & $1$ & $0$ \\ \hline
\end{tabular}
\end{center}
\caption{Indices of  fixed points of a pseudo-Anosov map}\label{indicespA}

\end{table}

\begin{table}
\begin{center}
\begin{tabular}{|l| lllll|}
\hline Type of $A$ &  $(p,0)^+$ &  $(p,k)^+$ & $(2n,0)^-$ & $(2n,1)^-$ & $(2n+1,0)^-$ \\
\hline $\Fix(f) \cap A$ & $a$ & $\emptyset$ & $2$ points & $2$ points & $2$ points \\
$\ind{f}{A}$ & $-p$ & $0$ & $0+0$ & $1+1$ & $1+0$ \\ \hline
\end{tabular}
\end{center}
\caption{Indices of invariant boundary components of a pseudo-Anosov
map}\label{indicesbdpA}
\end{table}

\section{Definition and properties of the standard form}\label{maquinista}

We say that  a homeomorphism $\Smap{f}{\Suu}$ in  Thurston canonical form is {\em in
standard form } \label{standardpage} if one of the following holds
\begin{romlist}
\item $f$ is finite-order
\item $f$ is a pseudo-Anosov map in standard form.
\item $f$ is reducible, every component of $f$ satisfies (i) or (ii), and for each connected component $A$ of $N(\Gamma)$, if $n$ is a positive integer such that $f^n(A)=A$, then $f^n|_{\Cl(A)}$  is conjugate to one of the following  maps of
$\SI \times I$.

\begin{numlist}
\item  $(z,t) \mapsto (z e^{2(a+bt)\pi i},t)$,  where $a$ and $b$ are rational numbers.
\item  $(z,t) \mapsto (\overline{z}e^{a(1-2t)\pi i},1-t)$, where $a$ is rational.
\item $(z,t) \mapsto (z e^{2a\pi i},1-t)$ where $a$ is a rational number.
\item  $(z,t) \mapsto (\overline{z},t)$;  where $\overline{z}$ denotes the  conjugate of $z$.
\end{numlist}
\end{romlist}
For each connected component $A$ of $N(\Gamma)$, the minimum  positive integer $m$ such
that $f^m(A)=A$ is called the {\em $f$-period of $A$}. \index{period!of a component of
$N(\Gamma)$}

Now we state properties of the standard form.

\start{lem}{class} Each isotopy class of homeomorphisms of a surface contains a map in
standard form.
\end{lem}

The following can  easily be deduced from the definition of standard form.

\start{lem}{iterates} Iterates of a map in standard form are also in  standard form.
\end{lem}

\start{rem}{rotation} If $\Smap{f}{\Suu}$ is an orientation-preserving map in standard
form and   $B$ is an $f$-invariant boundary component of $\Suu$ then $f|_B$ is a
rotation. Furthermore, if $B$ belongs to a pseudo-Anosov component of $f$ and  $p$ is the
number of prongs of $B$ then $f^p|_B=\id_B$  and all the prongs emanating from $B$ remain
fixed under the action of $f^p$.
\end{rem}

A description of the fixed-point classes of a map in standard form is given in the next
result, which  plays a key role in the proofs of Theorems~H and I.

\start{prop}{types} {\em (\cite[Lemma 3.6]{JG})} Let $\Smap{f}{\Suu}$ be a map in
standard form. Assume that $C$ is a fixed-point class of $f$. Then one of the following
holds.

\begin{manlist}{\Alph}
\item  $C=\{x\}$ where $x$ is an isolated fixed point, $f$ preserves orientation and one of the following holds:
\begin{manlist}{\Alph{llistai}.\arabic}
\item   $x \in \Int(\Suu)$ and $f$ is conjugate to a rotation in a neighborhood of $x$ and $\ind{f}{x}=1$.
\item $x \in \Int(\Suu)$ is a fixed point of a connected component  $A$  of $N(\Gamma)$, $f|_{\Cl(A)}$ is conjugate to a map of the form $(z,t) \mapsto (\overline{z}e^{a(1-2t)\pi i},1-t)$ and $\ind{f}{x}=1$.
\item  $x \in \Int(\Suu)$  is a type $(p,k)^+$ interior fixed point of  a pseudo-Anosov component and  $\ind{f}{x}=1-p$ if $k=0$ and $\ind{f}{x}=1$ otherwise.
\end{manlist}

\item $C=\{x\}$ where $x$ is an isolated fixed point and $f$ reverses orientation and  one of the following  holds:
\begin{manlist}{\Alph{llistai}.\arabic}
\item $x \in \Int(\Suu)$ is an interior fixed point of a pseudo-Anosov
component and $\ind{f}{x} \in \{1, -1, 0\}$.
\item  $x \in \partial \Suu$ is in a type $(p,k)^-$ invariant boundary component of some pseudo-Anosov component and $\ind{f}{x}=1$ or $0$.
\end{manlist}

\item   $C$ is a simple closed curve,  $f$ is orientation-preserving and one of the following  holds:

\begin{manlist}{\Alph{llistai}.\arabic}
\item $C \subset \Int(\Suu)$ and $C \subset A$ for some connected component  $A$ of $N(\Gamma)$,  $f_{\Cl(A)}$ is conjugate to a map of the form $(z,t) \mapsto (z e^{2(a+bt)\pi i},t)$ and $\ind{f}{C}=0$
\item  $C \subset \Int(\Suu)$, $C$ is a type $(p,0)^+$ boundary component of a pseudo-Anosov component of $f$, and $\ind{f}{C}=-p$.
\item $C \subset \partial \Suu$,  $C$ is a type $(p,0)^+$ boundary component of a pseudo-Anosov component of $f$ and $\ind{f}{C}=-p$.
\end{manlist}
\item   $C$ is a simple closed curve and  $f$ is orientation-reversing, $C \subset \Int(\Suu)$, in a neighborhood of $C$ $f$ is conjugate to the reflection $(z,t) \mapsto (z,1-t)$ and $\ind{f}{C}=0$.

\item $f$ is orientation-reversing, and $C$ is a fixed arc.  $C$ is contained in a subsurface $B$ on which $f$ acts as an involution (i.e., $f^2=\id$). Moreover, if $x$ is an endpoint of $C$ such that $x \in \Int(\Suu)$ then it is in a boundary component of a pseudo-Anosov component. Also, $\ind{f}{x} \in \{1, -1,0\}$.

\item $f$ is orientation-preserving and $C$ is a fixed subsurface  of $\Suu$  with $\chi(C) \le 0$. If  $B$ is a boundary component of $C$ such that $B \subset  \Int(\Suu)$ then $B$ is also a boundary component of either a  component of $N(\Gamma)$ or a component of a pseudo-Anosov component of $f$. In the latter case, $B$ is of type
$(p,0)^+$. Moreover $\ind{f}{C}=\chi(C)- \sum p_B <0$ where the summation is over the
components $B$ of $\partial C$ which are also boundary components of a pseudo-Anosov
component of $f$ and $p_B$ is the number of prongs emanating from $B$.
\end{manlist}
\end{prop}

\chapter{Planar discontinuous groups}\label{S5}
\renewcommand{\thetheo}{\thechapter.\arabic{theo}}

The goal of this chapter is to give a brief exposition of some features of the theory of
planar discontinuous groups. This theory will provide us the necessary tools  for the
construction and study of  finite-order homeomorphisms of closed surfaces in
Chapter~\ref{S6}.

Let us begin by introducing some notation and definitions.  We shall denote by $\uuu$ an
element of the set $\{+,-\}$. A {\it signature}\index{signature} consists of a sign $+$
or $-$ and an ordered sequence of integers with certain subsequences bracketed together
in the following manner: $$ (\uuu, T, [m_1, m_2, \dots,  m_{_R}],
\{(m_{1,1},m_{1,2},\dots, m_{1,M_1}), $$$$
  (m_{2,1},m_{2,2},\dots, m_{2,M_2}),\dots
  (m_{B,1},m_{B,2},\dots, m_{B,M_{_B}})\})
$$ and  satisfying:
\begin{numlist}
\item   $T \ge 0$, $R \ge 0$, $B \ge 0$, and $M_i \ge 0$ for each $i \in \{1,2,\dots, B\}$. Moreover, if the sign of the signature is $+$, $T$ is even;
\item If $T=B=0$ then the sign of the signature is $+$;
\item  For every $\ii{i}{R}$, $m_i \ge 2$ and $m_i \le m_j$ if $i \le j$;
\item  For every $\ii{i}{B}, \ii{j}{M_{_i}}$, $m_{i,j} \ge 2$.
\end{numlist}
Observe that  $R$, $M_i$ for some $i \in \{1,2,\dots,  B\}$ or $B$ can be $0$. In such
cases, the signature will be written with the brackets inserted, but with no symbols
between them. In fact, the signatures  which we will consider in  \lemc{dihedral}, and
afterwards, satisfy $M_i=0$  for every $\ii{i}{B}$. To abbreviate the notation we shall
denote these signatures by $$ (\uuu, T, [m_1, m_2, \dots,  m_{_R}],B). $$ If $$
\Psi=(\uuu, T, [m_1, m_2, \dots,  m_{_R}], \{(m_{1,1},m_{1,2},\dots, m_{1,M_1}), $$$$
  (m_{2,1},m_{2,2},\dots, m_{2,M_2}),\dots
  (m_{B,1},m_{B,2},\dots, m_{B,M_{_B}})\})
$$ is a signature, the {\it Euler characteristic} of $\Psi$  \index{Euler
characteristic!of a signature} is defined as $$
\mu(\Psi)=2-T-B-\sum_{i=1}^{R}\left(1-\frac{1}{m_i}\right)
-\frac{1}{2}\sum_{i=1}^{B}\sum_{j=1}^{M_i}\left(1-\frac{1}{m_{i,j}}\right). $$ We say
that a group $G$ is a {\em (cocompact)  planar discontinuous group }\index{cocompact
planar discontinuous group}  if there exists a signature $\Psi$ as in the above paragraph
such that $\mu(\Psi) \le 0$ and $G$ has a  presentation with generators:
\begin{numlist}
\item $\s_i$, $i \in \{1,2, \dots,  R\},$
\item $\tau_i$, $i \in \{1, 2,\dots,  T\},$
\item $\pi_i$, $i \in \{1, 2, \dots,  B\},$
\item $\rho_{i,j}$, $i \in \{1,2, \dots,  B\}$, $j \in \{1,\dots,  M_i+1\}$;
\end{numlist}
and defining relations
\begin{numlist}
\item $\sigma_i^{m_i}$, $i \in \{1,2, \dots,  R\},$
\item $\rho_{i,j}^2$,  $i \in \{1,2,\dots,  B\}$, $j \in \{1,2,3, \dots, M_i+1\}$,
\item $(\rho_{i,j}\rho_{i,j+1})^{m_{i,j}}$, $i \in \{1,2, \dots,  B\}$, $j \in \{1,2,\dots,  M_i\}$,
\item $\pi_{i}^{-1}\rho_{i,M_i+1} \pi_i\rho_{i,1}$, $\ii{i}{B}$;
\item
\begin{romlist}
\item  If  $\uuu=+$ then 
$$\pi_1 \pi_2, \dots, \pi_{_B} \s_1 \s_2 \dots, \s_{_R} [\tau_1,\tau_2][\tau_3,
\tau_4]\dots[\tau_{T-1},\tau_{T}]$$ (where  $[a,b]$ denotes $a b a^{-1} b^{-1}$).
\item  If $\uuu=-$ then 
$$\pi_1 \pi_2, \dots, \pi_{_B} \s_1 \s_2 \dots, \s_{_R}\tau_1^2 \tau_2^2 \dots,
\tau_T^2.$$
\end{romlist}
\end{numlist}

It can be proved that a planar discontinuous group $G$ determines  its signature  up to
certain permutations, see \cite[Theorems 1a, 2a and 3]{Ma}.  We require a weaker version
of this result:

\start{theo}{unicidad} If the signatures $$ (\uuu, T, [m_1, \dots,
m_{_R}],\{(m_{1,1},\dots, m_{1,M_1}),\dots
  (m_{B,1},\dots, m_{B,M_{_B}})\}) \mbox{ and }
$$ $$ (\widetilde{\uuu}, \widetilde{T}, [\widetilde{m_1}, \dots,
\widetilde{m}_{_R}],\{(\widetilde{m}_{1,1},\dots, \widetilde{m}_{1,M_1}),\dots
 (\widetilde{m}_{B,1},\dots, \widetilde{m}_{B,M_{_B}})\})
$$ are associated to a planar discontinuous group $G$ then $\uuu=\widetilde{\uuu}$,
$T=\widetilde{T}$, $R=\widetilde{R}$, $m_i=\widetilde{m}_i$ for each $\ii{i}{R}$,
$M_i=\widetilde{M}_i$ for each $\ii{i}{B}$ and $$ \{m_{i,1}, m_{i,2},\dots,
m_{i,M_i}\}=\{\widetilde{m}_{i,1}, \widetilde{m}_{i,2},\dots,  \widetilde{m}_{i,M_i}\} $$
for each  $\ii{i}{B}$.

In particular, the order of each of the generators is determined by the relations $(1)$
and $(2)$, and, for each $i \in \{1,2, \dots,  B\}$, $j \in \{1,2,\dots,  M_i\}$, the
order of $(\rho_{i,j}\rho_{i,j+1})$ is $m_{i,j}$.
\end{theo}

\begin{rem} The requirement that $m_i \le m_j$ if $i \le j $ in  the definition of a signature could have been dropped. In this case, the $m_i$ are uniquely determined up to permutation.
\end{rem}

Let $G$ be a planar discontinuous group with signature $\Psi$. We  define the {\it Euler
characteristic} of $G$ \index{Euler characteristic!of a group} as $$ \mu(G)=\mu(\Psi). $$
By \theoc{unicidad}, the Euler characteristic of a planar discontinuous group is well
defined. A planar discontinuous group is said to be {\em non-Euclidean}
\index{non-Euclidean group} (resp. {\em Euclidean}) \index{Euclidean group} if $\mu(G)<0$
(resp. $\mu(G)=0$).

Given a pla\-nar dis\-continuous group $G$, we define a homo\-mor\-phism
$$\map{\ve}{G}{\{-1,1\}},$$ called the {\em orientation map}, \index{orientation map}
which acts on  the generators in the following way: $$
\begin{array}{ll}
        \ve(\s_{i})=1&  \mbox{for each $\ii{i}{R}$,}\\
        \ve(\pi_i)= 1 & \mbox{for each $\ii{i}{B}$,}\\
        \ve(\rho_{ij})=-1&  \mbox{for each $\ii{i}{B}$, $\ii{j}{M_i}$},\\
        \ve(\tau_{i})=\lefteqn{
            \left\{
            \begin{array}{ll}
            1 & \mbox{if  $\uuu=+$,}\\
                         -1  & \mbox{otherwise.}
                      \end{array}
                      \right.
                      }&
    \end{array}
$$ A planar discontinuous group $G$ is said to be {\em orientable} if $\ve(G)=\{1\}$ and
{\em non-orientable} if $\ve(G)=\{-1,1\}$. By \theoc{unicidad} this definition is
independent of the signature of $G$.

A {\em surface group} \index{surface group} is a planar discontinuous group for which
$B=R=0$.

\start{rem}{fundamental} Notice that if $G$ is a surface group, then it is isomorphic to
the fundamental group of a surface. More precisely, if $G$ is orientable (resp.
non-orientable)  and has signature $(\uuu,T,[\,\,], \{\,\,\})$, then it is isomorphic to
$\pi_1(\Su{\frac{T}{2}})$ (resp. $\pi_1(N_T)$).  Here,  $\pi_1(\,)$
\index{1p@$\pi_1(\cdot)$} denotes  the fundamental group, and $N_T$  \index{1n@$N_T$} a
non-orientable compact connected closed surface of genus $T$. Observe that if $G$ is
orientable then $T$ is even and $$\mu(G)=\chi(\Su{\frac{T}{2}}),$$ and if $G$ is
non-orientable then $$\mu(G)=\chi(N_{T}).$$
\end{rem}

Now we state the well-known Riemann-Hurwicz Formula (see, for instance, \cite[Theorem 4.14.22]{Zie}

\start{theo}{RH} Let $G$ be a planar discontinuous group and let $H$ be a subgroup of $G$
of finite index.  Then $H$ is a planar discontinuous group and $$ \mu(H)=[G:H]\mu(G). $$
 Moreover, the orientation of $H$ is the restriction of the orientation of $G$.
\end{theo}
Here,  $[G:H]$ denotes the index of $H$ in $G$.

\start{rem}{palabras} Observe that, by \theoc{RH}, if $H$ is a subgroup of finite index
of a non-Euclidean (resp. Euclidean)  planar discontinuous group $G$ then $H$ is a
non-Euclidean (resp. Euclidean) planar discontinuous group. Further, if $G$ is orientable
so is $H$.
\end{rem}

As usual, if $G$ and $H$ are groups and $\map{\phi}{G}{H}$ is a homomorphism,
$\Ker(\phi)$ denotes the kernel of $\phi$. If $H$ is finite,  it follows from \theoc{RH}
that $\Ker(\phi)$  is also a planar discontinuous group.  In the following lemma we
determine necessary and sufficient conditions for $\Ker(\phi)$  be a surface group. The
proof   is elementary and can be found in \cite{Ha} for the particular case of
non-Euclidean planar discontinuous orientable groups. The general case is analogous and
for this reason, we do not include it in this exposition.

\start{lem}{preservation} Let  $\map{\phi}{G}{H}$ be a homo\-mor\-phism where $G$ is a
pla\-nar dis\-con\-tinuous group and $H$ is a finite group. Then $\phi$  preserves the
orders of the elements of finite order in $G$ if and only if  $\Ker(\phi)$ is a surface
group.
\end{lem}

In the following lemma, we state a property which must be satisfied for planar
discontinuous groups $G$ which admit epimorphisms  $\map{\phi}{G}{\Z_n}$ preserving the
orders of its finite-order elements.

\start{lem}{dihedral} Let $G$ be a planar discontinuous group with signature $\Psi$ and
let $n$ be a positive integer. If  there exists an epimorphism  $\map{\phi}{G}{\Z_n}$
which preserves the orders of finite-order elements of $G$ then  $M_i=0$ for each
$\ii{i}{B}$.
\end{lem}
\begin{proof}
Assume that $M_i \ge 1$ for some $\ii{i}{B}$. Since $\rho_{i,1}^2=\rho_{i,2}^2=1$, $n$
must be even and  $\phi(\rho_{i,1})=\phi(\rho_{i,2})=\clase{\frac{n}{2}}$. Hence,
$\rho_{i,1}\rho_{i,2} \in \Ker(\phi)$. On the other hand, by \theoc{unicidad}, the order
of  $\rho_{i,1}\rho_{i,2}$ is $m_{i,1} \ge 2$, a contradiction.
\end{proof}

>From  now on, we are going to consider only planar discontinuous groups $G$  which admit
epimorphisms  $\map{\phi}{G}{\Z_n}$ for some positive integer $n$ preserving the orders
of finite-order elements of $G$. By \lemc{dihedral}, the signature of such groups can be
denoted by $(\uuu,T,[m_1,m_2,\dots, m_{_R}],B)$.

Let $n, B, R$ and $1 \le p_1, p_2, \dots,  p_{_R}  < n$ be elements of $\N$. We say that
a  finite-order homeomorphism $\Smap{f}{\Suu}$ is of {\em type} \index{type of a
finite-order homeomorphism} $\typeor{n}{B}{p_1,p_2,\dots, p_{_R}}$ if the following
holds.
\begin{numlist}
\item  $n=\sf$;
\item For each $\ii{i}{R}$ there exists a periodic $f$-orbit $O_i$ of period $p_i$. Furthermore, $O_i  \cap O_j = \emptyset$ if $i \ne j$;
\item There exist $B$ simple closed curves $C_1, C_2, \dots,  C_B$ in $\Suu$ such that
\begin{romlist}
\item  For each $1 \le i,j \le n/2$, and for each $1\le l,m \le B $,
$$f^i(C_l) \cap f^j(C_m) = \emptyset, \mbox{ if } l \ne m \mbox{ or } i \ne j.$$
\item  $f^{\frac{n}{2}}|_{C_i}=\id_{C_i}$ for each $\ii{i}{B}$;
\end{romlist}
\item Every point $x \in \Suu \setminus \left(\bigcup\limits_{i=1}^{B}\bigcup\limits_{j=1}^{\frac{n}{2}} f^j(C_i) \cup \bigcup\limits_{i=1}^{R}O_i\right)$
has period $n$.
\end{numlist}

The following result will be crucial for our arguments, because, combined by Harvey's
\theoc{Harvey} and \theoc{g-b}, will give necessary and sufficient conditions for the
existence of a finite-order map of a given type. The proof can be found in \cite{G1} for
the particular case of  orientable planar discontinuous groups and orientable
homeomorphisms of prime order. This proof can be easily generalized for orientable groups
and homeomorphisms of any order. The more general statement we give here is a consequence
of some results of  \cite[Chapter 4]{Zie}.

\start{theo}{signatures} Let $n$ be a positive integer, let $G$ be a planar discontinuous
group with signature $$(\uuu,T,[m_1,m_2,\dots, m_{_R}],B)$$ and let $\map{\phi}{G}{\Z_n}$
be an epimorphism such that $\Ker(\phi)$ is an orientable surface group. Then there exists
a  finite-order homeomorphism $\Smap{f}{\Suu}$ of type
$$\typeor{n}{B}{\frac{n}{m_1},\frac{n}{m_2},\dots, \frac{n}{m_{_R}}}.$$ where $\Suu$ is a
closed surface such that $\pi_1(\Suu)$ is isomorphic to $\Ker(\phi)$. Moreover, the genus
of $\Suu$ equals $$ 1+\frac{n(T+B+R-2)-\sum_{i=1}^Rn/m_i}{2}. $$

Conversely, let $g$ be a positive integer and let $\Smap{f}{\Su{g}}$ be a finite-order
homeomorphism  of type $\typeor{n}{B}{p_1,p_2,\dots, p_{_R}}$. Then
$$T=2-B-R+\frac{2g-2+\sum_{i=1}^Rp_i}{n}$$ is a non-negative integer and  $f$ determines
an epimorphism from a planar discontinuous group $G$ of signature
$(\uuu,T,[\frac{n}{p_1},\frac{n}{p_2},\dots, \frac{n}{p_{_R}}],B)$ to $\Z_n$ such that
$\Ker(\phi)$ is isomorphic to $\pi_1(\Su{g})$.

In both cases, $f$ is orientation-preserving if and only if $G$ is orientable.
\end{theo}

By \theoc{signatures} (or by \lemc{jiji}) we have the following.

\start{cory}{B} If there exists $f \in \fgp$ of type $$\typeor{n}{B}{ p_1, p_2, \dots,
p_{_R}}$$ then $B=0$.
\end{cory}

The next theorem, due to Harvey \cite{Ha}, determines when, given a signature $\Psi$ and
an integer $n \ge 2$, there exists an epimorphism from an orientable group $G$ with
signature $\Psi$ to $\Z_n$ such that its kernel is a surface group.

\start{theo}{Harvey} {\em (Harvey)}  Let $n \ge 2$. Suppose that  $G$ is an orientable
planar discontinuous group of signature $$(+,T,[m_1,m_2,\dots, m_{_R}],B)$$ and let
$$M=\lcm(m_1,m_2,\dots, m_{_R}).$$ Then there exists an epimorphism $\map{\phi}{G}{\Z_n}$
such that $\Ker(\phi)$ is a surface group if and only if the following conditions are
satisfied:
\begin{numlist}
\item $\lcm(m_1,m_2,\dots, \hat{m_i},\dots, m_{_R})=M$ for each $\ii{i}{R}$ where $\hat{m_i}$ denotes the omission of $m_i$;
\item $M$ divides $n$, and, if $T=0$, $M=n$;
\item $R \ne 1$ and, if $T=0$, $R \ge 3$;
\item if $2^m$ divides $M$, and $2^{m+1}$ does not divide $M$ for some  positive integer $m$, the cardinal of the set $\{\ii{i}{R} \,\,:\,\, 2^m \mbox{ divides }m_i\}$ is even.
\end{numlist}
\end{theo}

The following proposition is the analogue of \theoc{Harvey} for orientable Euclidean
groups.

\start{lem}{paratoro} Let $G$ be a Euclidean planar discontinuous orientable group. Then
there exists an epimorphism $\map{\phi}{G}{\Z_n}$ such that $\Ker(\phi)$ is a surface
group if and only if one of the following holds,
\begin{numlist}
\item The signature of $G$ is $(+,0,[2,2,2,2],0)$ and $n=2$.
\item The signature of $G$ is $(+,0,[3,3,3],0)$ and $n=3$.
\item The signature of $G$ is $(+,0,[2,4,4],0)$ and $n=4$.
\item The signature of $G$ is $(+,0,[2,3,6],0)$ and $n=6$.
\item The signature of $G$ is $(+,2,[\,\,],0)$.
\end{numlist}
\end{lem}
\begin{proof} A simple calculation shows that if $\mu(G)=0$ then the signature of $G$ is one of those listed in $(1)$-$(5)$; see \cite[Table I]{Ma}. Thus, the result follows from \lemc{preservation}.
\end{proof}

\part{Development of the tools}\label{secretos}

\chapter{Finite-order maps of closed surfaces}\label{S6}

This chapter is dedicated to the study of homeomorphisms of finite order of closed
surfaces. This study has two goals: Firstly, to give a characterization of the possible
types of maps in $\fgp$ and in $\fgr$ which will allow us to determine, in
Chapter~\ref{S9}, $\mmm(\fgbp)$ and $\mmm(\fgbr)$. The second goal is the
construction of  maps with  ``large" minimum periods in Chapters~\ref{S9} and
\ref{S10}.

The organization of this chapter is as follows: in Sections~\ref{finiteop} and
\ref{lenguaje} we determine necessary and sufficient conditions for the existence of a
map with a given type in $\fgp$ and $\fgr$,  respectively. In Section~\ref{maniana}  we
construct some special maps in $\fgp$ and $\fgr$.

\section{The orientation-preserving case}\label{finiteop}

The following result is a consequence of Theorem~\ref{signatures} and Harvey's
\theoc{Harvey}.

\start{cory}{har} Let $g \ge 2$, $n \ge 2$, and $1 \le p_1, p_2, \dots, p_R < n$, be
positive integers. Set $$d=\gcd(p_1,p_2,\dots,p_R) \mbox{ and }
T=\frac{2g-2+\sum_{i=1}^Rp_i}{n}-R+2.$$ Then there exists $f \in \fgp$ of type
$\typeop{n}{p_1,p_2,\dots,p_R}$  if and only if the following conditions hold.
\begin{numlist}
\item $T$ is a non-negative even integer;
\item $\gcd(p_1,p_2,\dots,\hat{p_i},\dots,p_R)=d$ for each $\ii{i}{R}$;
\item For each $\ii{i}{R}$, $p_i$ divides $n$;
\item If $T=0$, $d=1$;
\item $R \ne 1$ and,  if $T=0$, $R \ge 3$.
\end{numlist}
\end{cory}

A necessary and sufficient condition for the existence of a finite-order map of certain
type is stated in the following result, which is an immediate consequence of
\coryc{har}.

\start{cory}{triangle} Let $p_1, p_2, p_3$ and $g$ be positive integers such that $g \ge
2$. There exists $f \in \fgp$  of type $\typeop{2g-2+p_1+p_2+p_3}{p_1,p_2,p_3}$ if and
only if    $p_1, p_2, p_3$ are pairwise coprime and $p_i$ divides $2g-2+p_1+p_2+p_3,$ for
each $i \in \{1,2,3\}$.

In particular, for each $g \ge 2$ and each $k$ dividing $g$, there exist
orientation-preserving  finite-order maps of $\Su{g}$ of the following types: $$
\typeop{2g+k}{1,1,k},  \typeop{4g}{1,1,2g }, \typeop{4g+2}{1,2,2g+1 }, $$ $$
\typeop{2g}{1,1,g }, \typeop{2g+1}{1,1,1 }, \typeop{2g+2}{1,1,2 }. $$
\end{cory}

We close this section with the  following consequence of \theoc{signatures} and
\lemc{paratoro}.

\start{lem}{ordertoro} The types of the elements of $\fgpx{1}$ are precisely
$$\typeop{2}{1,1,1,1}, \typeop{3}{1,1,1}, \typeop{4}{1,1,2},\typeop{6}{1,2,3} \mbox{ and
}$$ $$ \typeop{n}{\,\,}, n \ge 1.$$
\end{lem}

\section{The orientation-reversing case}\label{lenguaje}

The aim of this section is to state and give a proof  of the analogue of \coryc{har} for
the orientation-reversing case.

We need to prove the analogue of Harvey's \theoc{Harvey} for non-orientable planar
discontinuous groups. In this case,  we need to determine, not only the existence of an
epimorphism $\map{\phi}{G}{\Z_n}$ with kernel a surface group but also whether
$\Ker(\phi)$ is  orientable. The following lemma states sufficient and necessary
conditions for this.

\start{lem}{oli} Let $G$ be a non-orientable planar discontinuous group  and let
$\map{\phi}{G}{\Z_n}$ be an epimorphism. Then $\Ker(\phi)$ is orientable  if and only if
the following conditions hold:
\begin{numlist}
\item $n$ is even;
\item For every generator $x$ of $G$, $\ve(x)=-1$ if and only if $\minrep{\phi(x)}$ is odd.
\end{numlist}
\end{lem}
\begin{proof}  Assume that $\Ker(\phi)$ is orientable, so $\Ker(\phi) \subset \Ker(\ve)$. Since $\phi$ and $\ve$ are surjective by hypothesis, there exists an epimorphism $\Z_n \mapsto \{-1,1\}$ such that composed with $\phi$ it gives $\ve$. Hence (1) holds.  Since there is then a unique epimorphism $\Z_n \mapsto \{-1,1\}$ and it sends odd elements to $-1$ and even elements to $1$, (2) holds.

Reverse reasoning proves the converse.
\end{proof}

We  now need to introduce some notation. If $G$ is a planar discontinuous group of
signature $(\uuu,T,[m_1,m_2,\dots,m_R],B)$ and $n$ is a positive integer divisible by
each $m_i$, we write $\ggg(G,n)=\gcd(\frac{n}{m_1},\frac{n}{m_2},\dots,\frac{n}{m_R})$
and  $\suma(G,n)=\frac{1}{2}\sum_{i=1}^{R}\frac{n}{m_i}$

\start{theo}{g-b} Let $n$ be a positive integer and let $G$ be a non-orientable planar
discontinuous group of signature $$(\uuu,T,[m_1,m_2,\dots,m_R],B).$$ Then there exists an
epimorphism $\map{\phi}{G}{\Z_n}$ such that $\Ker(\phi)$ is an orientable surface group
if and only if the following conditions hold:
\begin{numlist}
\item $n$ is even;
\item each $m_i$ divides $n$;
\item $\ggg(G,n)$ is even;
\item  If $B \ge 1$ or  $\suma(G,n) \equiv_2 T+1$ then $\frac{n}{2}$ is odd;
\item If  $T+B=1$ then $\ggg(G,n)=2$.
\end{numlist}
\end{theo}
\begin{proof} We prove the ``only if'' direction first. Let $\map{\phi}{G}{\Z_n}$ be an epimorphism such that $\Ker(\phi)$ is an orientable surface group. By \lemc{oli}, (1) holds. By \lemc{preservation}, $\phi(\s_i)$ has order $m_i$ for each $\ii{i}{R}$, so  (2) holds and
\begin{equation}\label{minrep}
\minrep{\phi(\s_i)}=\frac{k_in}{m_i}
\end{equation}
for some positive integer $k_i$ such that $1 \le k_i <m_i$, and  $(k_i,m_i)=1$. By
\lemc{oli},  $\frac{k_in}{m_i}$ is even. Clearly, if $k_i$ is odd, then $\frac{n}{m_i}$
is even. If $k_i$ is even, then $m_i$ is odd, so by (1) $\frac{n}{m_i}$ is even. Hence,
(3) holds.

Let us see  (4). Assume first that $B \ge 1$. By \lemc{preservation},
$\phi(\rho_{11})=\clase{\frac{n}{2}}$ and by  \lemc{oli},
$\minrep{\phi(\rho_{11})}=\frac{n}{2}$ is odd. So, we may assume that $B=0$. Since $G$ is
non-orientable, $\uuu=-$. Suppose now, that $\frac{n}{2}$ is even. Here, by
\refeq{minrep}, for each $\ii{i}{R}$, $\frac{\minrep{\phi(\s_i)}}{2}\equiv_2
\frac{n}{2m_i}$. By \lemc{oli}, $\minrep{\phi(\tau_i)} \equiv_2 1$ for each $\ii{i}{T}$.
Then, since $\phi(\s_1\s_2\dots \s_R\tau_1^2\tau_2^2\dots\tau_R^2)=\clase{0}$, $$ 0
\equiv_2 \frac{1}{2}\sum_{i=1}^R\minrep{\phi(\s_i)}+\sum_{i=1}^T\minrep{\phi(\tau_i)}
\equiv_2 \suma(G,n)+T $$ which proves (4).

To see (5), observe that, by  \refeq{minrep}, $\ggg(G,n)$ divides
$\sum_{i=1}^R\minrep{\phi(\sigma_i)}$.  Assume that $T+B=1$. Then either $T=1$, $B=0$ and
$\uuu=-$ or $T=0$, $B=1$ and $\uuu=+$. In the former case, since $\phi(\s_1 \s_2 \dots
\s_R\tau_1^2)=\clase{0}$, $\ggg(G,n)$ divides $2\minrep{\phi(\tau_1)}$. On the other
hand, $\ima$ is generated by $\clase{\ggg(G,n)}$ and $\phi(\tau_1)$. Since $\phi$ is
surjective, $\gcd(\ggg(G,n),\minrep{\phi(\tau_1)})=1$. Then $\ggg(G,n)$ divides $2$ and
since it is even the proof of this case is complete. In the latter case, that is, when
$T=0$, $B=1$, and $\uuu=+$,  we have $\sum_{i=1}^R\phi(\s_i)+\phi(\pi_1)=\clase{0}$. Then
$\ggg(G,n)$ divides $\minrep{\phi(\pi_1)}$, so $\ima$ is generated by $\clase{\ggg(G,n)}$
and $\phi(\rho_{11})=\clase{\frac{n}{2}}$. Therefore, $\gcd(\ggg(G,n),\frac{n}{2})=1$. On
the other hand, by (4), $\frac{n}{2}$ is odd, so, since $\ggg(G,n)$ divides $n$,
$\gcd(\ggg(G,n),\frac{n}{2})=\frac{\ggg(G,n)}{2}$. Thus, (5) holds.

We now prove the reverse implication. Consider first the case where $\uuu=-$. If $T$ is
odd, define  $\phi$ on the generators in the following way: $$
\begin{array}{ll}
    \phi(\s_{i})=\clase{\frac{n}{m_i}} & \mbox{for $\ii{i}{R}$},\\
    \phi(\pi_{i})=\clase{2}& \mbox{for \ii{i}{B}},\\
                \phi(\rho_{i1})=\clase{\frac{n}{2}}  & \mbox{for $\ii{i}{B}$,}\\
                \phi(\tau_i)= \clase{(-1)^i} & \mbox{for $i \in \{2,3,\dots,T\}$,}\\
    \phi(\tau_{1})=\lefteqn{
                  \left\{
                                          \begin{array}{ll}
                  \clase{-\suma(G,n)-2B} & \mbox{if $\suma(G,n)$ is odd,}\\
                   \clase{-\suma(G,n)-2B+ \frac{n}{2}}  & \mbox{otherwise.}
                    \end{array}
                      \right.
                      }&
\end{array}
$$ We claim that $\phi$ is surjective. If $B \ge 1$, by (4), $\frac{n}{2}$ is odd. Then
$\clase{1} \in \ima$ since $\clase{2}, \clase{\frac{n}{2}} \in \ima$ and
$\gcd(2,\frac{n}{2})=1$. If $B=0$ and $T \ge 3$, $\phi$ is surjective because $\clase{1}
\in \ima$ by definition. If $B=0$ and $T=1$, by (5) $\ggg(G,n)=2$. Then $\clase{2} \in
\ima$. On the other hand, by definition,  $\minrep{\phi(\tau_1)}$ is odd. So,the claim
follows directly.

Suppose now that  $T$ is even. Since $\uuu=-$, $T \ne 0$. So, $T \ge 2$. In this case we
define $$
\begin{array}{ll}
    \phi(\s_{i})=\clase{\frac{n}{m_i}} & \mbox{for $\ii{i}{R}$},\\
    \phi(\pi_{i})=\clase{0}& \mbox{for \ii{i}{B}},\\
                \phi(\rho_{i1})=\clase{\frac{n}{2}}  & \mbox{for $\ii{i}{B}$,}\\
                \phi(\tau_i)= \clase{(-1)^i} & \mbox{for $i \in \{2,3,\dots,T\}$,}\\
    \phi(\tau_{1})=\lefteqn{
                  \left\{
                                          \begin{array}{ll}
                  \clase{-\frac{\suma(G,n)}{2}-1} & \mbox{for $\frac{\suma(G,n)}{2}$ is even,}\\
                   \clase{-\frac{\suma(G,n)}{2}-1+ \frac{n}{2}}  & \mbox{otherwise}.
                    \end{array}
                      \right.
                      }&
\end{array}
$$ Here, $\phi$ is surjective because $\clase{1} \in \ima$.

If  $\uuu=+$  we define $$
\begin{array}{ll}
        \phi(\s_{i})=\clase{\frac{n}{m_i}} & \mbox{for $\ii{i}{R}$},\\
                \phi(\tau_{i})=\clase{2} & \mbox{for $\ii{i}{T}$},\\
                    \phi(\rho_{i1})=\clase{\frac{n}{2}}  & \mbox{for $\ii{i}{B}$},\\
                    \phi(\pi_1)=\clase{2-2B-2\suma(G,n) }& \\
                    \phi(\pi_{i})=\clase{2}& \mbox{for $i \in\{2,3,\dots,B\}$}.\\
    \end{array}
$$

Notice that $B \ge 1$ because $G$ is non-orientable. Therefore, by (4), $\frac{n}{2}$ is
odd. Clearly,  if $\clase{2}, \clase{\frac{n}{2}} \in \ima$, $\phi$ is surjective. It is
trivial that this occurs when $B+T \ge 2$. Otherwise, since $B \ge 1$, $B+T=1$. Then, by
(5),  $\gcd(G,n)=2$. Hence,  in this case also, $\clase{2}, \clase{\frac{n}{2}} \in
\ima$.

That $\phi$ is an orientable surface group follows from  \lemc{oli} and
\lemc{preservation}.
\end{proof}

We now come to the main result of this section.

\start{cory}{anterior} Let $B,  R \ge 0$, $g \ge 2$, $n \ge 1$ and $1 \le p_1, p_2,
\dots, p_R <n$ be integers. Set $$ T=2-R-B+\frac{2g-2+\sum_{i=1}^Rp_i}{n}. $$ There exist
$f \in \fgr$ of type $\typeor{n}{B}{p_1,p_2, \dots, p_R}$ if and only if the following
conditions hold:
\begin{numlist}
\item For each $\ii{i}{R}$, $p_i$ divides $n$;
\item $n$ is even;
\item $T$ is a non-negative integer;
\item If $B=0$ then $T \ge 1$;
\item $\gcd(p_1,p_2, \dots, p_R)$ is even;
\item If $B+T=1$ then $\gcd(p_1,p_2, \dots, p_R)=2$;
\item  If $B \ge 1$ or $\frac{1}{2}\sum_{i=1}^R p_i \equiv_2 T+1$ then $\frac{n}{2}$ is odd.
\end{numlist}
\end{cory}
\begin{proof}
Assume that there exists such an $f$. Clearly, (1) holds. By \theoc{signatures}, there
exists a non-orientable planar discontinuous group $G$ with signature
$$\left(\uuu,T,\left[\frac{n}{p_1},\frac{n}{p_2},\dots,\frac{n}{p_R}\right],B\right)$$
and an epimorphism $\map{\phi}{G}{\Z_n}$ such that $\Ker(\phi)$ is isomorphic to
$\pi_1(\Sigma_g)$. Then (2), (3) and (4) hold, and (5), (6) and (7) follow from
\theoc{g-b}.

To see the converse, consider a non-Euclidean planar discontinuous group $G$ with
signature $$\left(\uuu, T,\left[\frac{n}{p_1},\frac{n}{p_2},\dots,\frac{n}{p_R}\right],
B\right)$$ such that if $B=0$ then  $\uuu=-$. Hence, $G$ is non-orientable and we can
apply  \theoc{g-b} to conclude the existence of an epimorphism $\map{\phi}{G}{\Z_n}$ such
that $\Ker(\phi)$ is an orientable surface group. By \theoc{signatures}, there exists  a
finite-order orientation-reversing homeomorphism $\Smap{f}{\Su{g}}$ of type
$\typeor{n}{B}{p_1,p_2,\dots,p_R}$. Now, we complete the proof, by observing that, by
\theoc{RH}, the signature of $\Ker(\phi)$ is $(+,g,[\,],0)$.
\end{proof}

Conditions (2) and (5) of \coryc{anterior} implies the following.

\start{lem}{suenio} If $\fgbr \ne \emptyset$ then $b$ is even.
\end{lem}

An application of  \coryc{anterior} yields the following result, which will be used in
the proof of Theorem~C.

\start{lem}{arpa} Let  $p, p_1, p_2$ and $g$ be positive integers such that $g \ge 2$.
\begin{numlist}
\item There exist $f \in \fgr$ of type $\typeop{2g-2+p}{p}$ if and only if $g$ is odd, $p$ is even, and $p$ divides $2g-2$.
\item There exist $f \in \fgr$ of type $\typeop{2g-2+p_1+p_2}{p_1,p_2}$ if and only if $g$ is even, $\gcd(p_1,p_2)=2$, and, for each $i \in \{1,2\}$, $p_i$ divides $2g-2+p_1+p_2$.
\end{numlist}
\end{lem}
\begin{proof}
The reverse implication in  (1) as well as the reverse implication in  (2) follow easily
from \coryc{anterior}. So, it only remains to prove both direct implications.

We start with (1). Set $T=2-1-0+\frac{2g-2+p}{2g-2+p}=2$. Assume that there exists $f \in
\fgr$ of type $\typeop{2g-2+p}{p}$. From \coryc{anterior}(5), it follows that $p$  is
even. By \coryc{anterior}(1), $p$ divides $2g-2+p$, so, $p$ divides $2g-2$. If
$\frac{p}{2}$ is even, since $\frac{p}{2}$ divides $g-1$ we have that $g$ is odd. If
$\frac{p}{2}$ is odd,   $\frac{p}{2} \equiv_2 T+1=3$. Hence  by \coryc{anterior}(7),
$\frac{n}{2}=g-1+p$ is odd. So, $g$ is odd.

Now, we prove  (2).   Set $T=2-2-0+\frac{2g-2+p_1+p_2}{2g-2+p_1+p_2}=1$. Assume that
there exists $f \in \fgr$ of type $$\typeop{2g-2+p_1+p_2}{p_1,p_2}.$$ By
\coryc{anterior}(1), $p_i$ divides $2g-2+p_1+p_2$ for each $i \in \{1,2\}$. Since
$B+T=0+1=1$, by \coryc{anterior}(6), $\gcd(p_1,p_2)=2$. Then, both $\frac{p_1}{2},
\frac{p_2}{2}$ are odd or one of them is even and the other odd. In the former case, we
have that $$\frac{p_1+p_2}{2} \equiv_2 0 \equiv_2 T+1,$$ and  by \coryc{anterior}(7),
$\frac{n}{2}=g-1+\frac{p_1+p_2}{2}$ is odd. So, $g$ is even. To study the latter case we
can assume without loss of generality that $\frac{p_1}{2}$ is even and $\frac{p_2}{2}$ is
odd. Since $p_1$ divides $2g-2+p_1+p_2$, $\frac{p_1}{2}$ divides $g-1+\frac{p_2}{2}$.
Hence, $g-1$ is odd. Therefore, $g$ is even and the proof  is complete.
\end{proof}

\section{Examples}\label{maniana}

The purpose of this section is to construct  finite-order orientation-reversing
homeomorphisms of closed surfaces. We shall do this  by ``gluing together'' finite-order
orientation-preserving ones. Some of the ideas of these constructions are based on an
example given in \cite{W1}.  These finite-order maps will be used in Chapter~\ref{S9} to construct maps with a given minimum period.

Before proving these result, let us introduce some notation. For each $\alpha \in \R$,
define $\Smap{R_\alpha}{\SI}$ by $z \mapsto z e^{2 \pi \alpha i}$.

\start{lem}{8g+4} Let  $g$ be a positive even integer. Then there exists $f \in \fgr$  of
type $$\typeor{4g+4}{0}{4,2g+2},$$ and a closed annulus $A \subset \Su{g}$  such that
$f|_A$ is conjugate to the map $(z, t) \mapsto (R_{ \frac{1 }{4g+4}}(z),1-t)$ on $\SI
\times [0,1]$.
\end{lem}
\begin{proof} Set $l =\frac{g}{2}$. By \lemc{ordertoro} and \coryc{triangle}, there exists $h \in \hgpx{l}$ of type $\typeop{4l+2}{1,2,2l+1}$.  We can assume without loss of generality that $\Su{l} \subset \R^3$. Let $D$ be an open disk centered at the fixed point of $h$ such that $h(D)=D$, and let $R$ be a plane which does not intersect $\Su{l}$.
Denote by $\Smap{s_{_R}}{\R^3}$ the reflection with respect to this plane.  Now define a
map $$\Smap{k}{(\Su{l} \setminus D) \cup s_{_R}(\Su{l} \setminus D)}$$
 in the following way:
\begin{eqnarray*}
k(x) &=& \left\{\begin{array}{ll}
        s_{_R}(x) & \mbox{if $x \in \Su{l} \setminus D$,}\\
        h(s_{_R}(x)) & \mbox{if $x \in s_{_R}(\Su{l} \setminus D)$.}
        \end{array}
        \right.
\end{eqnarray*}
Since $h$ has order $4l+2$, there exists $q$ coprime to $4l+2$ and a parametrization
$\map{\omega}{\SI}{\partial D}$ such that $\omega^{-1}\circ  h \circ
\omega=R_{\frac{q}{4l+2}}$. By taking a power of $h$ if necessary, we can assume that
$q=1$. Now define  an equivalence relation $\sim$ on $(\Su{l} \setminus D) \cup
s_{_R}(\Su{l} \setminus D)$ as follows: Let $x, y \in (\Su{l} \setminus D) \cup
s_{_R}(\Su{l} \setminus D)$. Then $x \sim y$ if and only if one of the following
statements holds:

\begin{numlist}
\item $x=y$;
\item $x \in \partial D$, $y \in s_{_R}(\partial D)$, $x=\omega(z)$ and $y=\sp(\omega(R_{-\frac{1}{8l+4}}(z)))$;
\item $y \in \partial D$, $x \in s_{_R}(\partial D)$,  $y=\omega(z)$ and $x=\sp(\omega(R_{-\frac{1}{8l+4}}(z)))$.
\end{numlist}
Observe that $$k(\sp(\omega(R_{-\frac{1}{8l+4}}(z))))=\omega(R_{\frac{1}{8l+4}}(z))
 \sim s_{_R}(\omega (z))=k(\omega(z)).$$
>From this we can conclude that $x \sim y$ implies $k(x) \sim k(y)$. On the other hand,
$(\Su{l} \setminus D) \cup s_{_R}(\Su{l} \setminus D)/\sim$  is a closed surface of genus
$g$. Hence,  $k$ induces $f \in \fgr$. Clearly, the type of $k$ is
$\typeor{4g+4}{0}{4,2g+2}$.

Denoting by $\map{p}{(\Su{l} \setminus D) \cup s_{_R}(\Su{l} \setminus D)}{\Su{g}}$ the
natural projection, it is not hard to see  that there exists an annulus $A \subset
\Su{g}$ invariant under $f$ such that $p(\partial D) \subset A$ and  $f|_{A}$ is
conjugate to the map  $(z, t) \mapsto (R_{\frac{1 }{4g+4}}(z),1-t)$ on $\SI \times
[0,1]$.
\end{proof}

\start{lem}{8g}  Let $g$ be a positive odd integer such that $g \ge 3$. Then there exists
$f \in \fgr$ of type $$\typeor{4g-4}{0}{2g-2},$$ and a closed annulus $A \subset \Su{g}$
such that $f|_A$ is conjugate to the map  $(z, t) \mapsto (R_{\frac{1 }{4g-4}}(z),1-t)$
on $\SI \times [0,1]$.
\end{lem}
\begin{proof} By \coryc{anterior}, there exists $h \in \mathop\mathcal{F}_{g-1}^r$ of type $$\typeor{4g-4}{0}{2,2g-2)}:$$
Let $D_1$, $D_2$ be  open disks  centered at  points of the $2$-periodic orbit of $h$
such that $h(D_1)=D_2$.  Set $$\eta=h|_{\Su{g-1}\setminus (D_1 \cup D_2)}.$$

Taking a power  of $\eta$ if necessary, we can assume that $\eta^2|_{\partial D_1}$ is
conjugate to rotation through an angle of $\frac{2 \pi}{4g-4}$. Then there exists a
parametrization $\map{\omega}{\SI}{\partial D_1}$ such that $\omega^{-1} \eta^2
\omega=R_{\frac{1}{4g-4}}$.

Let $\sim$ be the smallest equivalence relation defined on $\Su{g-1}\setminus (D_1 \cup
D_2)$ containing all the pairs of the form $$(x,\omega
R_{\frac{-1}{8g-8}}\omega^{-1}\eta(x))$$ for each $x \in \partial D_2$. Set $y=\omega
R_{\frac{-1}{8g-8}}\omega^{-1}\eta(x)$. Since $y \in \partial D_1$, $\eta(y) \in \partial
D_2$. Therefore, $$\eta(y) \sim \omega R_{\frac{-1}{8g-8}}\omega^{-1}(\eta^2(y))=\omega
R_{\frac{1}{8g-8}}\omega^{-1}(y)=\eta(x).$$ From this it follows that $x \sim y$ implies
that $\eta(x) \sim \eta(y)$ for every $x, y \in \Su{g-1}\setminus (D_1 \cup D_2)$. On the
other hand, $\Su{g-1}\setminus (D_1 \cup D_2)/\sim$ is an orientable closed surface of
genus $g$. It is easy to check that $\eta$ induces an element of $\fgr$ with the required
properties.
\end{proof}

\start{lem}{rosa} Let $g$ be a positive even integer. Then there exists $f \in \hgr$ of
type $$\typeor{2g-2}{0}{\,\,},$$ and a closed annulus $A \subset \Su{g}$   such that
$f|_A$ is conjugate to the map  $(z, t) \mapsto (R_{\frac{1}{2g-2}}(z),1-t)$ on $\SI
\times [0,1]$.
\end{lem}
\begin{proof} By \coryc{anterior} there exists $h \in \fgrx{g-1}$ of type $\typeor{2g-2}{0}{2}$. Let
$D_1$, $D_2$ be  open disks  centered at  points of the $2$-periodic orbit of $h$ such
that $h(D_1)=D_2$. Now we can proceed as in \lemc{8g}, considering
$S=\Su{g-1}\setminus(D_1 \cup D_2)$, $\eta=h|_S$, ``gluing'' the boundary components in
an appropriate way and taking the map induced by $\eta$ in the quotient surface.
\end{proof}

\start{lem}{intimos dones} Let $g$ be a positive integer even. Then there exists $f \in
\fgr$ of type $\typeop{4g}{2,2g}$ and a  closed annulus $A \subset \Su{g}$   such that
$f|_A$ is conjugate to the map  $(z, t) \mapsto (R_{\frac{1}{4g}}(z),1-t)$ on $\SI \times
[0,1]$.
\end{lem}
\begin{proof} Set $k=\frac{g}{2}$. By \lemc{ordertoro} and \coryc{triangle} there exists $h \in \fgpx{k}$ of type $\typeop{4k}{1,1,2k}$. By taking an $h$-invariant disk around one of the fixed points of $h$, it is easy to check that we can complete the proof as in \lemc{8g+4}.
\end{proof}

\start{rem}{pan} Dropping the requirement of the existence of the invariant annulus $A$
in Lemmas~\ref{8g+4}, \ref{8g}, \ref{rosa} and \ref{intimos dones},  these results follow
directly from Corollaries~\ref{har} and \ref{anterior}. Nevertheless, the existence of
such an annulus will be fundamental for the construction of examples, as we will see in
Chapter~\ref{S9}\end{rem}

\chapter{Consequences of fixed-point theory}\label{S7}

This chapter is devoted to stating some facts concerning the Lefschetz numbers of
homeomorphisms of surfaces. In Section~\ref{key} we prove a result which contains one of
the main ideas of this thesis, that is, a relationship between homeomorphisms of surfaces
with boundary and homeomorphisms of surfaces without boundary.  Some algebraic tools are
described in  Section~\ref{algebra}. In Section~\ref{felix} we prove some results about
sequences of Lefschetz numbers which will be used in Chapter~\ref{S9} and
\ref{S10}. Finally, Section~\ref{lasubseccionqueaterra} is devoted to the (tedious)
study of three particular classes of homeomorphisms of closed surfaces which will be used
in Chapter~\ref{S11}.

\section{The induced map}\label{key}

We begin this section by introducing a definition which we will use frequently. Let $f
\in \hgb$ and let $T$ be an $f$-invariant subset of $\Suf{g}{b}$ consisting of $k$
boundary components of $\Suf{g}{b}$.  In $\Suf{g}{b}$ we consider the equivalence
relation where two points in $\Suf{g}{b}$ are equivalent if they are equal or they belong
to the same boundary component in $T$.  The quotient space is  a surface of genus $g$
with $b-k$ boundary components. Since $T$ is $f$-invariant, $f$ induces a homeomorphism
$\tf^T \in \hgbx{g}{b-k}$. In particular, we can take $T$ to be the boundary of
$\Suf{g}{b}$.  In this case $b=k$, so $f$ induces a homeomorphism of $\Su{g}$ which will
be denoted by $\tf$ and called the {\em homeomorphism induced by $f$}.  \index{induced
homeomorphism} \index{1f@$\tf$}

The aim of this section is to prove a topological result, \propc{index}, which will be
fundamental for (almost all) our arguments and states a relationship between the fixed
points of a map $\Smap{f}{\Suf{g}{b}}$ and the index of the fixed points of the induced
map. In this way, we will be able to obtain information on the fixed-point set of maps of
surfaces with boundary by studying the (simpler) case of maps of surfaces without
boundary.

\start{lem}{paraindex} Let   $f \in \hgb$. Assume that there exists an $f$-invariant
boundary component $B$ of $\Suf{g}{b}$. Set $T=\{B\}$ and consider
$$\Smap{\tf^T}{\Suf{g}{b-1}},$$ the map induced by $f$ on $\Suf{g}{b-1}$ by collapsing
$B$ to a point, $q$. If $f|_B$ is fixed-point free  then $\ind{\tf^T}{q}=1$.
\end{lem}
\begin{proof} Suppose $f|_B$ is fixed-point free. Here, there exists an open set $U \subset \Suf{g}{b}$ such that $B \subset U$ and $f|_U$ is also fixed-point free. Set $V= \Suf{g}{b} \setminus B$. Clearly, $V$ is an open set, $V \cup U=\Suf{g}{b}$ and $U \cap V \cap \Fix(f)=\emptyset$. Hence, we can apply \theoc{suma} to conclude that
$$ L(f)=\inn(f|_U)+\inn(f|_V). $$ Since  $f|_U$ is fixed-point free,  by \remc{pillow
book},  $\inn(f|_U)=0$,  so
\begin{equation}\label{U}
L(f)=\inn(f|_V).
\end{equation}
Let $p$ denote the projection from $\Suf{g}{b}$ to the quotient space $\Suf{g}{b-1}$.
Then $p(U)$ and $p(V)$ are open sets contained in $\Suf{g}{b-1}$ and satisfy $p(U) \cup
p(V)=\Suf{g}{b-1}$, $p(U) \cap p(V) \cap \Fix(\tf^T)= \emptyset$. Thus, by \theoc{suma},
\begin{equation}\label{hat}
L(\tf^T)=\inn(\tf^T|_{p(U)})+\inn(\tf^T|_{p(V)}).
\end{equation}
On the other hand, setting $A=B$ in  \propc{collapse} we obtain $$
L(f)+1=L(\tf^T)+L(f|_B). $$ Since $f|_B$ is fixed-point free, by \theoc{lef} ,
$L(f|_B)=0$. By \refeq{U} and \refeq{hat} we get

\begin{equation}\label{ring}
\inn(f|_V)+1=\inn(\tf^T|_{p(U)})+\inn(\tf^T|_{p(V)}).
\end{equation}
>From the definition of $\tf^T$ it follows that $\inn(\tf^T|_{p(V)})=\inn(f|_V)$ so,
\refeq{ring} becomes $$ \inn(\tf^T|_{p(U)})=1. $$ Clearly
$\ind{\tf^T}{q}=\inn(\tf^T|_{p(V)})$, so our proof is complete.
\end{proof}

\start{prop}{index} Let $f \in \hgb$ and let $\Smap{\tf}{\Su{g}}$ be the map induced by
$f$. If,  for some $i$, $f^i$ is fixed-point free, then all the fixed points of
$\widetilde f^i$ are isolated and have local index  equal to $1$ with respect to
$\widetilde f^i$.
\end{prop}

\begin{proof} Let $x$ be a fixed point of $\tf^i$. Since $f^i$ is fixed-point free, $x$ is a collapsed boundary component. Since the set of collapsed boundary components is finite, $x$ is isolated. That the index is  $1$ follows from \lemc{paraindex}.
\end{proof}

\section{Basic algebraic tools}\label{algebra}

The goal of this section is to present and prove some   elementary algebraic facts
concerning the map $\fs$ induced by a map $\Smap{f}{\Suu}$. The reason for doing this is
that by \refeq{numerar} there exists a close relationship between the map $\fs$ and the
sequence of Lefschetz numbers, $L(f), L(f^2), \dots$.

Recall that, by \propc{efectos}, there exists a basis such that for each  $f \in \hg$,
the matrix representing $\fs$ is symplectic. In the following lemma we state a property
of such matrices. The proof of the second statement  is due to Fenchel and can be found
in \cite{N}, but we include it here for convenience. A proof of this statement for the
improper symplectic case can also be found in \cite{FL}.

\start{lem}{causas} Suppose $M \in GL_{2g}(\Z)$ is proper  (resp. improper) symplectic
and let $P(x)=\d(x\id_{2g}-M)$ be its characteristic polynomial. Then
\begin{numlist}
\item $\d(M)=1$ (resp. $\d(M)=(-1)^g$).
\item  $P(x)=x^{2g}P(1/x)$ (resp. $P(x)=(-1)^g x^{2g}P(-1/x)$).
\end{numlist}
\end{lem}
\begin{proof} If $M$ is proper symplectic then (1) holds by \propc{gracias}. Now, define    $U_g \in GL_{2g}(\Z)$  as
$$ \left (\begin{array}{ll}
         0&I_g\\
         I_g &0
    \end{array}
    \right ).
$$ It is not hard to see that $U_g$ is improper symplectic and $\d(U_g)=(-1)^g$. Let $M
\in GL_{2g}(\Z)$ be an improper symplectic matrix. Since $U_gM$ is proper symplectic,
then, by \propc{gracias}, $(-1)^g\d(U_g)=\d(M U_g)=1$. Hence, (1) holds.

To see (2), assume that $M \in GL_{2g}(\Z)$ is proper symplectic. (The argument where $M$
is improper symplectic is similar.) Then, $$ M^t J_g(\lambda I_{2g}-M)=\lambda M^t
J_g-M^t J_gM=(\lambda M-I_{2g})^t J_g$$ $$ =-\lambda
I_{2g}\left(\frac{1}{\lambda}I_{2g}-M\right)^t J_g, $$ and (2) follows by taking
determinants of the two extremes of the above chain of equalities.
\end{proof}

We need to study the sequence  of integers $$L(f), L(f^2), \dots$$ which, by
\refeq{numerar},  is equivalent to studying the sequence $$p_1, p_2, \ldots.$$ By
Newton's equations (page \pageref{newton}), this infinite sequence is closely related to
the sequence $$s_1, s_2, \ldots, s_{2g}.$$ For this reason, in the next result, which
follows directly from \propc{efectos} and \lemc{causas}, we state some relations
satisfied by the $s_i$'s.

\start{lem}{symmetry} If $f \in \hgp$ then $\d(\fs)=1$ and  $s_{h}=s_{2g-h}$ for all $h
\in \{1,\dots,2g\}$. If $f \in \hgr$ then $\d(\fs)=(-1)^g$ and $s_{h}=(-1)^{g+h}s_{2g-h}$
for all $h \in \{1,\dots,2g\}$.
\end{lem}

\start{rem}{torre} By \lemc{symmetry} and Newton's equations (page \pageref{newton}),
for each $f \in \hg$, the infinite sequence $L(f), L(f^2), \dots$ is determined by the
finite sequence  $L(f), L(f^2), \ldots, L(f^g)$. In other words, if we are given $$L(f),
L(f^2), \dots, L(f^g)$$ then we can compute $L(f^i)$ for each positive integer $i$.
Hence, if  $f, h \in \hgp$ (resp. $f,h \in \hgr$) and $L(f^i)=L(h^i)$ for each
$\ii{i}{g}$ then $L(f^i)=L(h^i)$ for each positive integer $i$.
\end{rem}

The arguments to prove the next result, from which Proposition~A follows as a corollary,
are based on Nielsen  \cite{N}.

\start{lem}{paratraza} Let  $A \in GL_k(\R)$  where $k\geq 3$. If $\tr (A^i)=1$ for $i
\in \{-1, 1, 2, \ldots, k-2\}$ then $\tr (A^{k-1})\ne 1$.
\end{lem}

\begin{proof} The assumptions imply that $p_i=1$ for $i=1,2, \ldots, k-2$. Then by Newton's equations (page \pageref{newton}),

\begin{eqnarray*}
1+s_1 &=& 0,\\ 1+s_1 +2s_2&=& 0,\\ 1+s_1  +s_2  +3 s_3 &=& 0,\\ \dots & &\\ 1+s_1 +s_2
+s_3 +\dots + (k-2) s_{k-2} &=& 0,\\ p_{k-1}+s_1 +s_2 +s_3 +\dots + (k-1) s_{k-1} &=&
0.\\
\end{eqnarray*}
Hence, $s_1=-1$, $s_i=0$ for each $i \in \{2,3,\dots, k-2\}$ and
$(k-1)s_{k-1}=1-p_{k-1}$. Therefore,  $p_{k-1}=1$ if and only if $s_{k-1}=0$. On the
other hand, since the eigenvalues of $A$ are the inverses of the eigenvalues of $A^{-1}$,
$$s_{k-1}=(-1)^{k-1}\d(A)\tr(A^{-1})=(-1)^{k-1}\d(A) \neq 0.$$ Therefore, $p_{k-1} \ne 1$
and the proof is complete.
\end{proof}

To prove our main results we need to prove the existence of fixed-point classes of $f^m$
with negative index, for certain $m$. In view of \lemc{essential}, one way to do this is
to show that $L(f^m)<0$. For this reason we are going to study some conditions on $f$
which imply that inequality.

\start{lem}{potency} Suppose that $g \ge 2$ and that \Smap{f}{\Su{g}} is a homeomorphism
such that, for some positive integer $m$, the $m$-th power of each  of the eigenvalues of
$\fs$ is equal to $1$.  Then $L(f^m)<0$.
\end{lem}
\begin{proof} By hypothesis, $\tr(\fs^m)=\sum_{i=1}^{2g}\lambda_i^m=2g$. Then, by \refeq {anne}, if $f^m$
preserves orientation  then $L(f^m)=2-\tr(\fs^m)=2-2g <0$ and if $f^m$ reverses
orientation  then $L(f^m)=-\tr(\fs^m)=-2g<0$.
\end{proof}

\section{Sequences of Lefschetz numbers I}\label{felix}

The results of this section, as well as the ones of the next, are consequences of
\refeq{numerar}, Newton's equations (page \pageref{newton}) and \lemc{symmetry}.

\start{lem}{nortumbria} Let $f \in \hgp$ be such that $L(f^i)=1$ for each $\ii{i}{g}$.
Then $L(f^{g+1})=-g$.
\end{lem}
\begin{proof} By the assumptions and \refeq{numerar}, $p_i=1$ for each $\ii{i}{g}$, and, by Newton's equations (N.1), (N.2) (page \pageref{newton}), $\dots$, (N.g), $s_1=-1$  and $s_i=0$ for each $i \in \{2,3,\dots,g\}$. By \lemc{symmetry}, $s_{g+1}=s_1=-1$. Hence, we can replace all these values in Newton's equation (N.(g+1)) (page \pageref{newton}) to obtain $p_{g+1}=g+2$. Thus, by \refeq{numerar}, $L(f^{g+1})=2-p_{g+1}=-g$, as desired.
\end{proof}

\start{lem}{amistad} If $g \ge 2$ then there does not exist $f \in \hgp$ such that
$L(f^i)~\ge~4$ for each $\ii{i}{2g}$.
\end{lem}
\begin{proof} Assume that there does exist such an $f$.
By \refeq{numerar},  $p_i=2-L(f^i) < -1$, for each $\ii{i}{2g}$. We claim that $s_i >1$
for each $\ii{i}{2g}$. We prove this claim  by induction. By Newton's equation \refeq{N1}
(page \pageref{newton}),  $s_1=-p_1 >1$. Hence, the claim holds for $i=1$. Now, assume $j
\le 2g$ is such that $s_i >1$ for each $\ii{i}{j-1}$. Observe that $p_is_{j-i} <-1$ for
each $\ii{i}{j-1}$. By Newton's equations (page \pageref{newton}), $$
-js_j=p_1s_{j-1}+p_2s_{j-2}+\dots+p_{j-1}s_1+p_j<-j, $$ so the claim holds for $j$ and we
are done. In particular, we have proved that $s_{2g}>1$. On the other hand, by
\lemc{symmetry}, $s_{2g}=\d(\fs)=1$, which is impossible.
\end{proof}

\start{lem}{verlaine} If $g \ge 2$ then there does not exist $f \in \hgr$ such that
$L(f^{2i})>4$  and $L(f^{2i-1})=0$ for each $\ii{i}{g}$.
\end{lem}
\begin{proof} Assume that there does exist such an $f$. Let $\ii{i}{2g}$.
By \refeq{numerar},
\begin{equation}\label{ay}
p_{2i}=2-L(f^{2i}) <-2 \mbox{ and  } p_{2i-1}=-L(f^{2i-1})=0.
\end{equation}
We claim that  $s_{i}>1$ if $i$ is even and $s_{i}=0$ if $i$ is odd. We prove this claim
by induction. By Newton's equation \refeq{N1} (page \pageref{newton}),  $s_1=-p_1=0$.
Hence, the claim holds for $i=1$. Now, assume $j \le 2g$ is such that the inductive
hypothesis holds for each $\ii{i}{j-1}$. Suppose that $j$ is even. Observe that by the
inductive hypothesis and \refeq{ay}, $p_is_{j-i}<-2$ for each even $\ii{i}{j-1}$. By
Newton's equations (page \pageref{newton}) and \refeq{ay} $$
-js_j=p_1s_{j-1}+p_2s_{j-2}+\dots+p_{j-1}s_1+p_j<-j, $$ so the claim holds  in this case.
If  $j$ is odd, for each $\ii{i}{j}$, $p_is_{j-1}=0$ because, by the inductive hypothesis
and \refeq{ay}, either $p_i=0$ or $s_{j-i}=0$. Hence, $$
-js_j=p_1s_{j-1}+p_2s_{j-2}+\dots+p_{j-1}s_1+p_j=0, $$
 so  the claim is proved.
In particular, we have proved that $s_{2g}>1$. On the other hand, by \lemc{symmetry},
$s_{2g}=\d(\fs)=(-1)^g$, which is impossible.
\end{proof}

\section{Sequences of Lefschetz numbers II}\label{lasubseccionqueaterra}

In Section\refc{S11} we will need to prove the existence of  fixed-point classes of
negative index for  pseudo-Anosov maps in standard form. As we will see, the
pseudo-Anosov maps with more than two orbits of singularities  can be ``controlled'' with
the help of the Euler-Poincar\'{e} Formula \refeq{ep}.  Some of the homeomorphisms with one
or two orbits of singularities do not offer many difficulties and can be studied by using
the Euler-Poincar\'{e} Formula or \lemc{potency}. However, there are some particular cases of
pseudo-Anosov homeomorphisms with one or two orbits of singularities which need  special
treatment. The results of this section will allow us to deal with these cases.

\subsection{The orientation-preserving case}
We start
 with two technical lemmas. The first one will be used to prove the second one which, in turn, will be used in the proof of \propc{cas-part-op}.

\start{lem}{reufa} Let $A \in GL_k(\Z)$ be such that, for some positive integer $n$ and
some non-negative integer $j$,  the coefficients $s_i$ of the characteristic polynomial
of $A$ satisfy $$
\begin{array}{lll}
s_i & = & \cases{-2 & if\,  $i \equiv_n  1$, \cr
                  1 & if\,  $i \equiv_n 0, 2$, \cr
                  0 &   \mbox{otherwise},}
\end{array}
$$ for non-negative $i$ such that $i \le nj$. Then, for each  positive integer $t$  such
that $nj+t \le k$ $$ p_{nj+t}+\sum_{h=0}^{j-1} p_{nh+t}-2\sum_{h=1}^{j}p_{nh+t-1}
+\sum_{h=1}^{j}p_{nh+t-2} $$$$ +\sum_{h=1}^{t-1}s_{nj+h}p_{t-h}+(nj+t)s_{nj+t}=0. $$
\end{lem}
\begin{proof}
Since $nj+t \le k$, the result follows directly from  Newton's equation N.(nj+t) (page
\pageref{newton})
\begin{equation}
p_{nj+t}+s_1 p_{nj+t-1}+s_2 p_{nj+t-2}+\dots + (nj+t) s_{nj+t} = 0,
\end{equation}
by replacing  the $s_h$'s by their values for  $h \in \{1, 2, \ldots, nj\}$ and grouping
the $p_h$'s for $h \in \{t, t+1, \ldots, nj+t-1 \}$ according to the value of $h$ modulo
$n$.
\end{proof}

\start{lem}{polinor} Let $g \ge 2$, let $f \in \hgp$, let $v$, $n$ and $i$  be  positive
integers  such that $v  \le 2g$ and $n \ge 3$.   If for each $i \le v$ $$
\begin{array}{lll}
L(f^i) & = & \cases{n & if\,  $n$ divides $i$, \cr
                    0 &   \mbox{otherwise},}\cr
\end{array}
$$ then,  for each $i \le v$ $$
\begin{array}{lll}
s_i & = & \cases{-2 & if\,  $i \equiv_n  1$, \cr
                  1 & if\,  $i \equiv_n 0, 2 $, \cr
                  0 &   \mbox{otherwise}.}
\end{array}
$$ Consequently, if $j,t $ are positive integers such that $nj \le v$  and $nj+t \le 2g$
then $$ p_{nj+t}+\sum_{h=0}^{j-1} p_{nh+t}-2\sum_{h=1}^{j}p_{nh+t-1}
+\sum_{h=1}^{j}p_{nh+t-2} $$ $$ +\sum_{h=1}^{t-1}s_{nj+h}p_{t-h}+(nj+t)s_{nj+t}=0. $$
\end{lem}

\begin{proof}
Assume that we have proved the first statement of the lemma. Then the second one follows
immediately by  \lemc{reufa}. So, let us prove the formula for the $s_i$'s.

By  \refeq{numerar},
\begin{eqnarray*}
p_i & = & \cases{2-n & if\,  $n$ divides $i$,\cr
                 2   &   \mbox{otherwise,}}\cr
\end{eqnarray*}
for $i \leq v$. We shall proceed by  induction. By Newton's equation \refeq{N1}
 (page \pageref{newton}), the equality holds for $i=1$. Consider now $1 \le i \le v-1$ and suppose that
the statement holds for $s_1,s_ 2, \cdots,s_ i$. Let  $j$ and $r$ be the non-negative
integers such that $i=nj+r$ and $0 \leq r < n$. Recall that a matrix representing $\fs$
has size $2g \times 2g$ and observe that $2 \le i+1=nj+r+1 \le v \le 2g$. By  the
induction hypothesis and \lemc{reufa} with $t=p+1$, $$ p_{i+1}+\sum_{h=0}^{j-1}
p_{nh+r+1}-2\sum_{h=1}^{j}p_{nh+r}+\sum_{h=1}^{j}p_{nh+r-1} $$
\begin{equation}\label{gatitosblancos}
 +s_{nj+1}p_r+s_{nj+2}p_{r-1}+\cdots+s_{nj+r}p_1+(i+1)s_{i+1}=0
\end{equation}

Considering  separately the four cases, $r=0, 1, n-1$  and $r \notin \{ 0, 1, n-1\}$,  we
show that the values of the $p_l$'s for $l \in \{1, 2, \ldots, i+1\}$ and $s_l$'s for $l
\in \{nj+1,nj+2, \ldots, nj+r \}$  can be replaced to find $s_{i+1}$. Indeed, if $r=0$,
then \refeq{gatitosblancos} can be written as

$$ 2+2j-2(2-n)j+2j+(-2)(2-n)+1.2+0+\cdots+ (nj+1)s_{i+1}=0, $$ which gives $s_{i+1}=-2$.

If $r=1$ then  \refeq{gatitosblancos} gives $$2 +2j-4j+(2-n)j-4+(nj+2)s_{i+1}=0, $$ so
$s_{i+1}=1$.

If $r=n-1$ then  \refeq{gatitosblancos} becomes $$ 2-n+(2-n)j-2\cdot 2 j+2j-2\cdot 2
+2+(nj+n)s_{i+1}=0, $$ and hence, $s_{i+1}=1$.

Finally, suppose that $r \notin \{ 0, 1, n-1\}$. Then, by \refeq{gatitosblancos}, $$
2+2j-4j+2j-2\cdot 2+2+(i+1)s_{i+1}=0, $$ which gives $s_{i+1}=0$.
\end{proof}

The following result will be used to prove \propc{pA}.

\start{prop}{cas-part-op} Let $g \ge 2$ and let $f \in \hgp$  satisfy the following
conditions:
\begin{numlist}
\item $1, 2, g+1 \notin \Per(f)$.
\item There exists a positive integer $n$ such that $\Per(f) \cap \{3,4, \ldots, g\}=\{n\}$. Moreover, there is only one periodic orbit of period $n$ and each of its points has index $1$ for $f^{hn}$ for  each positive integer $h$ such that  $hn \le g$.
\end{numlist}
Then either there exists $m$ such that  $1 \leq m \leq 3g-3$ and   $L(f^m) <0$, or there
exists $l \in \Per(f)$ such that $g+2 \le l \le \frac{4}{3}(g-1)$
\end{prop}
\begin{proof} Here there are only finitely many periodic points of $f$ of period less than
or equal to $2g$, so they are isolated. Therefore, by \lemc{suma}, if $1 \leq i \leq 2g$,
then the Lefschetz number of $f^i$ is the sum of the local indices of its fixed points,
so, for each $i \leq g+1$,
\begin{eqnarray*}
L(f^i) & = & \cases{n & if\,  $n$ divides $i$,\cr
                    0 &   \mbox{otherwise}.}\cr
\end{eqnarray*}
Consequently, by  \refeq{numerar},
\begin{equation}
p_i = \left\{\begin{array}{ll}
        2-n & \mbox{if $n$ divides $i$},\\
                   2& \mbox{otherwise,} \\
                              \end{array}
        \right.
        \label{pedeperrito}
\end{equation}
for $i \le g+1$.

Let $r, j$ be the non-negative integers such that $g=nj+r$ and $0 \le r \le n-1$. We
claim that $r \notin \{0,2\}$. Since $n \ge 3$, by Lemmas \refc{polinor} and
\ref{symmetry}, if $r=0$,  $-2=s_{g+1}=s_{g-1} \in \{0,1\}$, which is a contradiction.
Similarly, if $r=2$ then $-2=s_{g-1}=s_{g+1} \in \{0,1\}$. This completes the proof of
the claim. Now we split the proof into two cases.

\begin{case}{1}{ $r=1$.}\end{case}
By Lemmas \refc{polinor} and \ref{symmetry},
$$P(x)=(x-1)^2(x^{2nj}+x^{n(2j-1)}+\ldots+x^n+1)=\frac{x^{2nj+n}-1}{x^n-1}(x-1) ^2.$$
Then $(x^n-1)P(x)=(x^{2nj+n}-1)(x-1)^2$. Hence, if $\lambda \ne 1$ is a root of $P(x)$
then $0=(\lambda^n-1)P(\lambda)=(\lambda^{2nj+n}-1)(\lambda-1)^2$. Therefore the
$(2nj+n)$-th power of the roots of $P(x)$ are $1$. Since these roots are the eigenvalues
of $\fs$, \lemc{potency} completes the proof of this case because $2nj+n \le 2g-2+g-1\le
3g-3$.

\begin{case}{2}{ $r \ge 3$.} \end{case}
We  start by showing by induction on $t$ that if $r+1 \le t <\min\{n, 2r-2\}$ then
$p_{nj+t}=2$. Since $nj+r+1=g+1$ and $0 \le r+1 \le t <n$, $n$ does not divide $g+1$, so,
by \refeq{pedeperrito}, the result holds for $t=r+1$.  Now fix $r+1 < t <\min\{n, 2r-2\}$
and assume that $p_{nj+s}=2$ for every $r+1 \le s < t$.

Recall that a matrix representing $\fs$ has size $2g \times 2g$ and observe that $1 \le
nj+t  < nj+n <2g$. Then,  by \lemc{polinor} (taking $v$ as $g+1$), $$
p_{nj+t}+\sum_{h=0}^{j-1} p_{nh+t}-2\sum_{h=1}^{j}p_{nh+t-1} +\sum_{h=1}^{j}p_{nh+t-2} $$
\begin{equation}\label{holala}
+\sum_{h=1}^{t-1}s_{nj+h}p_{t-h}+(nj+t)s_{nj+t}=0.
\end{equation}
Since $4 \le r+1 \le t <n$,   $t-2, t-1, t \not\equiv_n 0$. By Lemmas \refc{polinor} and
\ref{symmetry},  $s_{nj+1}=-2$, $s_{nj+2}=1$ and $s_{nj+s}=0$ if $0 \le s <t$. Further,
our induction hypothesis implies that $p_{nj+s}=2$ if $r+1 \le s<t$. Now, substituting in
\refeq{holala} these values and the ones given by \refeq{pedeperrito}, we obtain
\begin{equation}
p_{nj+t}+2j-4j+2j-2 \cdot 2 +2 +(nj+t)s_{nj+t}=0.\label{finnn}
\end{equation}
Since $r+1 \le t <2r-2$, we have  $2 < 2r-t<r$. Therefore, by Lemmas\refc{polinor} and
\ref{symmetry},
 $s_{nj+t}=s_{2g-(nj+t)}=s_{nj+2r-t}=0$. Hence,  \refeq{finnn}  gives $p_{nj+t}=2$. This completes the induction step. Now we divide the proof of this case into three subcases.

\begin{case}[Subcase]{2.1}{$2r-2>n$.}\end{case}
By \lemc{polinor}, $$ p_{nj+n}+\sum_{h=0}^{j-1} p_{nh+n}-2\sum_{h=1}^{j}p_{nh+n-1}
+\sum_{h=1}^{j}p_{nh+n-2} $$ $$ +\sum_{h=1}^{n-1}s_{nj+h}p_{n-h}+(nj+n)s_{nj+n}=0. $$ We
have proved that $p_{nj+s}=2$ for $r+1 \le  n-1$. Therefore, $$ p_{nj+n}+(2-n)j-4j+2j-2
\cdot 2 +2+(nj+n)s_{nj+n}=0. $$ Observe that $2 <2r-n<2n-n=n$, so, by \lemc{symmetry},
$s_{nj+n}=s_{2g-(nj+n)}=s_{nj+2r-n}=0$. Hence $p_{nj+n}=2+nj$. By  \refeq{numerar},
$L(f^{nj+n})=-nj<0$. Since $nj+n\le g-3+g-3 < 3g-3$ the desired conclusion holds in this
case.

\begin{case}[Subcase]{2.2}{$2r-2=n$.}\end{case}
By Lemmas  \refc{polinor} and \ref{symmetry},
$$P(x)=(x-1)^2(x^{n(2j+1)}+x^{2nj}+\ldots+x^n+1)= \frac{x^{2nj+2n}-1}{x^n-1}(x-1)^2.$$
Observe that   $2nj+2n \le 2g-6 +g-3 < 3g-3$. Since, as above,  we can see that the
$2(nj+n)$-th powers of the eigenvalues of $\fs$ are $1$, from \lemc{potency} we can
deduce that $L(f^{2nj+2n})<0$. Taking $m=2nj+2n$, we can complete the proof.

 \begin{case}[Subcase]{2.3}{$2r-2 <n$.} \end{case}
 By \lemc{polinor} (taking $v=g+1$),
 $$
p_{nj+2r-2}+\sum_{h=0}^{j-1} p_{nh+2r-2}-2\sum_{h=1}^{j}p_{nh+2r-3}
+\sum_{h=1}^{j}p_{nh+2r-4} $$ $$
+\sum_{h=1}^{2r-3}s_{nj+h}p_{2r-2-h}+(nj+2r-2)s_{nj+2r-2}=0. $$ Since $2 \le 2r-4
<2r-3<2r-2<n$,   $p_{nh+2r-3}=p_{nh+2r-4}=2$ if $0 \le h \le j$ and $p_{nh+2r-2}=2$ if $0
\le h <j$. Further, $s_{nj+1}=-2$ and $s_{nj+2}=1$. Hence, $$
p_{nj+2r-2}+2j-4j+2j-4+2+(nj+2r-2)1=0. $$ So $p_{nj+2r-2}=2-(nj+2r-2)$ and
$L(f^{nj+2r-2})=nj+2r-2 \ne 0$. This assures the existence of a periodic orbit whose
period divides $l=nj+2r-2$. Since $n+4 \le nj+2r-2<nj+n$, $n$ does not divide $nj+2r-2$,
so the period of the orbit is $nj+2r-2$.  Now observe that  $r \ge 3$, so  $nj+2r-2 \ge
g+1$. Moreover, there are no periodic orbits of period $g+1$. Therefore, $nj+2r-2 \ge
g+2$. Since $j \ge 1$, $$3(nj+2r-2)=3nj+3r+2r-2+r-4 < 3g+g-r+r-4=4g-4.$$ Thus, taking $l$
as $nj+2r-2$ the desired conclusion holds for this subcase.
\end{proof}

\subsection{The orientation-reversing case}

The following four lemmas show  results analogous to Lemmas~\ref{reufa}, \refc{polinor}
and \propc{cas-part-op} for the orientation-reversing case. Since, by \lemc{symmetry},
the equation relating the values of the $s_i$'s corresponding to an orientation-reversing
map of $\Su{g}$ depends on the parity of $g$, the analogue of \propc{cas-part-op} for the
orientation-reversing case, splits into two results, Propositions\refc{cas-part-or-even}
and \ref{cas-part-or-odd}. The proof of the following lemma is analogous of the proof of
\lemc{reufa}.

\start{lem}{reufaor} Let $A \in GL_k(\Z)$ be such that, for some positive integer $n$ and
some non-negative integer $j$,  the coefficients $s_i$ of the characteristic polynomial
of $A$ satisfy,

$$
\begin{array}{lll}
s_i & = & \cases{1& if\,  $i \equiv_n  0$, \cr
                  -1 & if\,  $i \equiv_n  2 $, \cr
                  0 &   \mbox{otherwise},}
\end{array}
$$ for $i \le nj$. Then, for each  positive integer $t$  such that $nj+t \le k$, $$
p_{nj+t}+\sum_{h=0}^{j-1} p_{nh+t}-\sum_{h=1}^{j}p_{nh+t-2}+
\sum_{h=1}^{t-1}s_{nj+h}p_{t-h}+(nj+t)s_{nj+t}=0. $$
\end{lem}

\start{lem}{pol-or} Let $g \ge 2$ and  let $f \in \hgr$. If there exist positive integers
$n$ and $v$ such that  $n$ is even, $n \ge 3$ and $v  \le 2g$  such that, for every $i
\le v$, $$
\begin{array}{lll}
L(f^i) & = & \cases{n & if\,  $n$ divides $i$, \cr
                    0 &   \mbox{otherwise},}\cr
\end{array}
$$ then, for $i \le v$, $$
\begin{array}{lll}
s_i & = & \cases{-1 & if\,  $i \equiv_n  2$, \cr
                  1 & if\,  $i \equiv_n 0  $, \cr
                  0 &   \mbox{otherwise}.}
\end{array}
$$ Consequently, if $j, t $ are positive integers such that $nj \le v$  and $nj+t \le 2g$
then $$ p_{nj+t}+\sum_{h=0}^{j-1} p_{nh+t}-\sum_{h=1}^{j}p_{nh+t-2}+
\sum_{h=1}^{t-1}s_{nj+h}p_{t-h}+(nj+t)s_{nj+t}=0. $$
\end{lem}

\begin{proof}
Assume that we have proved the first statement of the lemma. Then, the second one follows
by  \lemc{reufaor}. Now, let us prove the formula for the $s_i$'s.

By  \refeq{numerar},
\begin{equation}
p_i = \left\{\begin{array}{ll}
        2-n & \mbox{if $n$ divides $i$},\\
        0 &  \mbox{if $i$ is odd,}\\
                    2 & \mbox{otherwise,} \\
                              \end{array}
        \right.
        \label{aranias}
\end{equation}
for $i \leq v$.

Now  we  proceed by  induction. By  Newton's equation \refeq{N1} (page \pageref{newton}),
\refeq{aranias} holds for $i=1$. Consider now $1 \le i \le v-1$ and suppose that the
statement holds for $s_1,s_ 2, \cdots,s_ i$. Let $j$ and $r$ by  the non-negative
integers such that $i=nj+r$ and $0 \leq r < n$. Recall that a matrix representing $\fs$
has size $2g \times 2g$ and observe that $2 \le i+1=nj+r+1 \le v \le 2g$. By  the
induction hypothesis and \lemc{reufaor},
\begin{equation}\label{leoncitosblancos}
p_{i+1}+\sum_{h=0}^{j-1} p_{nh+r+1}-\sum_{h=1}^{j}p_{nh+r-1}
\end{equation}
$$ \hspace{-2cm}+\sum_{h=1}^{r}s_{nj+h}p_{t-h}+(i+1)s_{i+1}=0. $$

Since $n$ is even, if $r$ is even, then $nj+r-1$ and $nj+r+1$ are odd for every positive
integer $j$. Therefore, all the terms of the above sum except $(i+1)s_{i+1}$ are equal to
$0$. So, $s_{i+1}=0$. So we can assume that $r$ is odd. Considering  separately the three
cases, $r=1, n-1$  and $r \notin \{ 1, n-1\}$, we show that the values of the $p_l$'s for
$l \in \{1, 2, \ldots, i+1\}$ and $s_l$'s for $l \in \{nj+1,nj+2, \ldots, nj+r \}$  can
be replaced to find $s_{i+1}$. Indeed, if $r=1$, then \refeq{leoncitosblancos} can be
written as $$ p_{i+1}+2j-j(2-n)+(nj+2)s_{i+1}=0, $$ which gives $s_{i+1}=-1$.

If $r=n-1$ then  \refeq{leoncitosblancos} gives $$ p_{i+1}+j(2-n)-2j+2+(nj+1)s_{i+1}=0,$$
so, $s_{i+1}=1$. Finally, if $r$ is odd and different from $1$ and $n-1$, by
\refeq{leoncitosblancos}, $$2+2j-2j-2+(i+1)s_{i+1}=0,$$ which gives  $s_{i+1}=0$.
\end{proof}

The following proposition will be used in the proof of \propc{pAgeven}.

\start{prop}{cas-part-or-even} Let $g$ be a positive even integer such that $g \ge 6$,
and let $f \in \hgr$ satisfying the following conditions:
\begin{numlist}
\item $\Per(f) \cap \{1,2,3,4,5,g+1,g+2\}=\emptyset$.
\item There exists an even positive integer $n$ such that $\Per(f) \cap \{6,7, \ldots, g\}=\{n\}$. Moreover, there is only one orbit of period $n$ and each  of its points has index $1$ for $f^{hn}$  for each positive integer $h$ such that $hn \le g$.
\end{numlist}
Then there exists a positive integer $m$ such that  $m \le 2g-6$ and $L(f^m) <0$.
\end{prop}

\begin{proof} By the same arguments as in the proof of \propc{cas-part-op} we find

\begin{eqnarray*}
L(f^i) & = & \cases{n & if\,  $n$ divides $i$,\cr
                    0 &   \mbox{otherwise.}}\cr
\end{eqnarray*}
Consequently, by  \refeq{numerar},
\begin{equation}
p_i = \left\{\begin{array}{ll}
        2-n & \mbox{if $n$ divides $i$},\\
        0 &  \mbox{if $i$ is odd,}\\
                    2 & \mbox{otherwise,} \\
                              \end{array}
        \right.
        \label{sapitos}
\end{equation}
for $i \le g+2$.

Let $r, j$ be the non-negative integers such that $g=nj+r$ and $0 \le r \le n-1$. Since
$n \le g$, $j \ge 1$. Observe that since $n$ and $g$ are even so is $r$. In particular,
$r \notin \{ 1,3\}$. We split the proof into  three cases.

\begin{case}{1}{ $r=0$. } \end{case}
Here, $g=nj$. Since $g+2 \le 2g$,  by \lemc{pol-or}, $$ p_{g+2}+\sum_{h=0}^{j-1}
p_{nh+2}-\sum_{h=1}^{j}p_{nh}+ \sum_{h=1}^{1}s_{nj+h}p_{2-h}+(g+2)s_{g+2}=0. $$ By
\refeq{sapitos},  $p_{g+2}+nj+(g+2)s_{g+2}=0$.  By Lemmas \refc{symmetry} and
\ref{pol-or}, $s_{g+2}=s_{2g-(g+2)}=s_{g-2}=s_{n(j-1)+n-2} =0$ because $2<n-2<n$.
Therefore, $p_{g+2}=-nj=-g$ and, by  \refeq{numerar},   $L(f^{g+2})=g+2$. This implies
the existence of an orbit whose period divides $g+2=nj+2$. But this contradicts our
hypothesis because the only orbit of period at most $g+2$ is the orbit of period $n$ and
since $n >2$, $n$ does not divide $g+2$.

\begin{case}{2}{ $r =2$.} \end{case}
By \lemc{pol-or}, $$ p_{g+2}+\sum_{h=0}^{j-1} p_{nh+4}-\sum_{h=1}^{j}p_{nh+2}+
\sum_{h=1}^{3}s_{nj+h}p_{4-h}+(g+2)s_{g+2}=0. $$ Since $n>4$  is even, $p_{nh+4}=2$ for
$h \in\{0,1,\ldots, j-1\}$ and $p_{nh+2}=2$ for $h \in \{0,1,\ldots,j \}$.  By
\lemc{pol-or}, $s_{g}=s_{m j+2}=-1$. Then $$ p_{g+2}+2j-2j-2+(g+2)s_{g+2}=0. $$ By Lemmas
\refc{symmetry} and \ref{pol-or},  $s_{g+2}=s_{2g-(g+2)}=s_{g-2}=s_{nj}=1$. So,
$p_{g+2}=-g$ and $L(f^{g+2})=g+2$. Since $g+2=nj+4$ and $n \ge 6$, by an  argument
analogous to that of the previous case, we obtain  a contradiction.

\begin{case}{3}{ $r >2$.}\end{case}
We will show by induction on $t$ that if $r+1 \le t <\min\{n,2r-2\}$ then
\begin{equation}
p_{nj+t} = \left\{\begin{array}{ll}
        2  &\mbox{if $t$ is even},\\
                    0 & \mbox{if $t$ is odd.} \\
                              \end{array}
        \right.
        \label{toros}
\end{equation}
Since $nj+r+1=g+1 <g+2$ and $4 \le r+1 \le t <n$, $n$ does not divide $g+1$, by
\refeq{sapitos},  the first step of the induction. Now fix $r+1 < t <\min\{n,2r-2\}$ and
assume that \refeq{toros} holds for every $r+1 \le s <t$.  Observe that $nj+t<nj+n<2g$.
Then, by  \lemc{pol-or}, (taking $v=g$), $$ p_{nj+t}+\sum_{h=0}^{j-1}
p_{nh+t}-\sum_{h=1}^{j}p_{nh+t-2}+ \sum_{h=1}^{t-1}s_{nj+h}p_{t-h}+(nj+t)s_{nj+t}=0. $$
Since $n$ is even, if $t$ is odd all the terms of the form $s_{nj+h}p_{t-h}$ for  $h \in
\{1,2, \ldots, t-1\}$ are equal to $0$, because $s_{nj+h}=0$ for each $h \in
\{1,3,\ldots, t-1\}$ and $p_{t-2}=0$.  By \lemc{symmetry},
$s_{nj+t}=s_{2g-(nj+t)}=s_{g+2r-t}=0$ because $2<2r-t < 2r-(r+1)<n$. Thus $p_{nj+t}=0$.
Assume now that $t$ is even. Since $4 \le r+1< t <n$ we have  $t-2, t \not\equiv_n  0$.
By Lemmas \refc{pol-or} and \ref{symmetry}  $s_{nj+2}=-1$ and $s_{nj+l}=0$ if $l \in
\{1,3,4 \ldots, t-1\}$. Then, substituting these values in the above equation, and using
\refeq{sapitos} and the induction hypothesis, we obtain $$
p_{nj+t}+2j-2j-2+s_{nj+t}(nj+t)=0. $$ Since $2<2r-t<n$, by \lemc{symmetry},
$s_{nj+t}=s_{2g-(nj+t)}=s_{nj+2r-t}=0$. So, $p_{nj+t}=2$ as desired. Now we split the
proof of this case into two subcases.

\begin{case}[Subcase]{3.1}{$2r-2 \ge n$.} \end{case}
By \lemc{pol-or}, $$ p_{nj+n}+\sum_{h=0}^{j-1} p_{nh+n}-\sum_{h=1}^{j}p_{nh+n-2}+
\sum_{h=1}^{n-1}s_{nj+h}p_{n-h}+(nj+n)s_{nj+n}=0. $$ Since $n \ge 3$, $p_{nh+n-2}=2$ and
$p_{nh+n}=2-n$ for $0 \le h <j$. Further, we have proved that $p_{nj+n-2}=2$. So, by
\refeq{sapitos}, $$ p_{nj+n}+(2-n)j-2j-2+(nj+n)s_{nj+n}=0. $$ Since $2 \le 2r-n <n$, by
\lemc{symmetry},  $s_{nj+n}=s_{2g-(nj+n)}=s_{nj+2r-n} \in \{0, -1\}$. So, $p_{nj+n} \in
\{2+nj,2+2nj+n \}$. In all cases, $L(f^{nj+n})<0$ which completes this subcase because
$nj+n<2g-2r \le 2g-6$.

\begin{case}[Subcase]{3.2}{$2r-2<n$.} \end{case}
By \lemc{pol-or}, $$ p_{nj+2r-2}+\sum_{h=0}^{j-1} p_{nh+2r-2}-\sum_{h=1}^{j}p_{nh+2r-4}$$
$$ +\sum_{h=1}^{2r-3}s_{nj+h}p_{2r-2-h}+(nj+2r-2)s_{nj+2r-2}=0. $$ By Lemmas
\refc{pol-or} and \ref{symmetry},   $s_{nj+l}=0$ if $l \in \{1,3, \ldots, 2r-3\}$. Since
$r \ge 3$, $2 \le 2r-4 < 2r-2<n$. Thus,   $p_{nh+2r-2}=p_{nh+2r-4}=2$ if $h \in \{0, 1,
\ldots, j-1\}$. Also, $p_{nj+2r-4}=2$. Further, $s_{nj+2}=-1$. Hence $$
p_{nj+2r-2}+2j-2j-2+(nj+2r-2)s_{nj+2r+2}=0. $$ By Lemmas \refc{symmetry} and
\ref{pol-or},   $s_{nj+2r-2}=s_{2g-(nj+2r-2)}=s_{nj+2}=-1$. So, $p_{nj+2r-2}=2+nj+2r-2$.
By  \refeq{numerar},   $L(f^{nj+2r-2})=-(nj+2r-2)<0$. Since $nj+2r-2 <  g-r+n \le g-r+g-r
\le 2g-6$, we are done.
\end{proof}

The following result will be used in the proof of \propc{pAgodd}

\start{prop}{cas-part-or-odd} Let $g$ be a positive odd integer such that $g \ge 3$ and
let $f \in \hgr$ satisfy the following conditions:
\begin{numlist}
\item $1, 2 \notin \Per(f)$.
\item There exists an even positive integer $n$ such that $\Per(f) \cap \{3,4, \ldots, g\}=\{n\}$. Moreover, there is only one orbit of period $n$ and each of its points has index $1$ for $f^{hn}$ for each positive integer $h$ such that $hn \le g$.
\end{numlist}
Then either there exists $m$ such that $1 \leq m \leq 3g-3$ and $L(f^m) <0$, or there
exists  $l \in \Per(f)$ such that $g+1 \le l \le \frac{4g-4}{3}$.
\end{prop}

\begin{proof} Using the same arguments as in the proof of  \propc{cas-part-op} we  find

\begin{eqnarray*}
L(f^i) & = & \cases{n & if\,  $n$ divides $i$,\cr
                    0 &   \mbox{otherwise,}}\cr
\end{eqnarray*}
for $i \le g$. Consequently, by  \refeq{numerar},
\begin{equation}
p_i = \left\{\begin{array}{ll}
        2-n & \mbox{if $n$ divides $i$},\\
        0 &  \mbox{if $i$ is odd,}\\
                    2 & \mbox{otherwise,} \\
                              \end{array}
        \right.
        \label{conejitos}
\end{equation}
for $i \le g$.

Let $r, j$ be the non-negative integers such that $g=nj+r$ and $0 \le r \le n-1$. Since
$n \le g$, $j \ge 1$. Observe that $n$ is even and $g$ is odd, so $r$ must be odd. In
particular, $r \notin \{0, 2\}$. We split the proof into three cases.

\begin{case}{1}{ $r =1$.}   \end{case}
By Lemmas  \refc{pol-or} and \ref{symmetry}
$$P(x)=(x^2-1)(x^{2nj}+x^{n(2j-1)}+\ldots+x^n+1)=\frac{x^{2nj+n}-1}{x^n-1}(x-1) ^2.$$ As
in the proof of  \propc{cas-part-op} we can see that this implies that the $(2nj+n)$-th
power of each of  the eigenvalues of $\fs$ is $1$ and since $2nj+n \le 2g-2+g-1=3g-3$, we
are done.

\begin{case}{2}{ $r \ge 3$.} \end{case}
By an argument analogous to that of the  Case $r \ge 3$ of  the proof of
 \propc{cas-part-or-even} (doing the first step of the induction with $r$ instead of $r+1$) we can show by induction on $t$ that if $r \le t <\min\{n,2r-2\}$ then
\begin{equation}
p_{nj+t} = \left\{\begin{array}{ll}
        2  &\mbox{if $t$ is even},\\
                    0 & \mbox{if $t$ is odd.} \\
                              \end{array}
        \right.
        \label{vacas}
\end{equation}

Now we split the proof of this case into three subcases.

\begin{case}[Subcase]{2.1}{$2r-2= n$.}\end{case}
By Lemmas  \refc{pol-or} and \ref{symmetry},
$$P(x)=(x^2-1)(x^{2nj+n}+x^{n(2j-1)}+\ldots+x^n+1)=\frac{x^{2nj+2n}-1}{x^n-1}(x-1) ^2.$$
As in the proof of \propc{cas-part-op}, we can see that this implies that the
$(2nj+2n)$-th power of each  eigenvalue of $\fs$ is $1$. Since $2nj+2n=2g-2r+n+n \le
2g-2+g-3=3g-5<3g-3$, we are done.

\begin{case}[Subcase]{2.2}{$2r-2 > n$.} \end{case}
Using the  same arguments as in the Case $2r-2 \ge n$  of the proof of
\propc{cas-part-or-even} we can prove that $p_{nj+n}+(2-n)j-2j-2+(nj+n)s_{nj+n}=0$. It
follows from \lemc{symmetry} that $s_{nj+n}=-s_{2g-(nj+n)}=-s_{nj+2r-n}=0$ because
$2<2r-n<n$. Thus $p_{nj+n}=2+nj$, so $L(f^{nj+n})=-nj<0$. Since $nj+n\le 2g-6 <3g-3$, we
are done.

\begin{case}[Subcase]{2.3}{$2r-2<n$.}\end{case}
By an argument analogous to that for the Case $2r-2<n$ of the proof of
\lemc{cas-part-op}, we obtain $$ p_{nj+2r-2}+2j-2j-2+(nj+2r-2)s_{nj+2r+2}=0. $$ By Lemmas
\refc{symmetry} and \ref{pol-or},   $s_{nj+2r-2}=-s_{2g-(nj+2r-2)}=-s_{nj+2}=1$, so
$p_{nj+2r-2}=2-(nj+2r-2)$. By  \refeq{numerar}, $L(f^{nj+2r-2})=nj+2r-2$. As  in  Case
$2r-2<n$ of the proof \propc{cas-part-op}, this assures the existence of a periodic orbit
of period $l=nj+2r-2 \ge g+1$ such that $3l \le 4g-4$.

\end{proof}

\chapter{Consequences of the Thurston-Nielsen theory}\label{S8}

This chapter is devoted to studying some results about homeomorphisms in  Thurston
canonical form and in standard form. In Section~\ref{2pa} we state some properties of
pseudo-Anosov maps. In Section~\ref{standarita} we determine sufficient conditions for
the existence of fixed-point classes of negative index of iterates of pseudo-Anosov maps
and reducible maps in standard form, and, in Section~\ref{sporm}, we study finite-order
maps and reducible maps.

\section{Properties of pseudo-Anosov maps}\label{2pa}
\markboth{THURSTON-NIELSEN THEORY}{PROPERTIES OF PSEUDO-ANOSOV MAPS}

The following  lemmas are  consequences of the Euler-Poincar\'{e} Formula \refeq{ep}, and
they will be used in Chapter\refc{S11}.

\start{lem}{b4} If there exists a pseudo-Anosov map $\Smap{f}{\Suf{0}{b}}$ then
$b~\ge~4$.
\end{lem}

\begin{proof} Suppose $f$ exists. Let $\widetilde{b}$ be the number of one-pronged boundary components of  the foliation on $\Suf{0}{b}$. Let $B'$ be a boundary component from which emanates more than one prong and $p_{B'}$ denote the number of such prongs. Then by \remc{pprong},
$\sum_{s \in \Sing(B')} =-p_{B'} \le -2$. Since there are $b-\widetilde{b}$  such
components, by the Euler-Poincar\'{e} Formula  \refeq{ep}, $$2(2-b)=2
\chi(\Suf{0}{b})=\sum_{s \in \Sing(\Suu)}(2-p_s)$$ $$ \le \sum_{s \in
\Sing(\partial(\Suu)) }(2-p_s) \le - \widetilde{b} - 2(
b-\widetilde{b})=\widetilde{b}-2b.$$ Thus $4 \le \widetilde{b}$ and, since $\widetilde{b}
\le b$, the result follows.
\end{proof}

\start{lem}{libros} Let $g \ne 1$ and let $\Smap{f}{\Su{g}}$ be a pseudo-Anosov map. If
$k$ is the number of $f$-orbits of singularities of the foliation on $\Su{g}$ then $k \ge
1$.
\end{lem}
\begin{proof}
If $k=0$ then the foliation has no singularities. Hence, by the Euler-Poincar\'{e} Formula
\refeq{ep},  $g=1$, which is impossible.
\end{proof}

\section{Fixed-point classes of iterates of maps in standard form}\label{standarita}
\markboth{THURSTON-NIELSEN THEORY}{FIXED-POINT CLASSES} This section is devoted to
stating some conditions which ensure the existence of  fixed-point classes of negative
index for iterates of pseudo-Anosov maps and reducible maps in standard form.  Most of
these results are consequences of the characterization of fixed-point classes given in
Chapter~\ref{S4}, mainly in \propc{types}.

\start{lem}{tango} Let $p$ and $n$  be  positive integers and  let $\Smap{f}{\Su{g}}$ be
a pseudo-Anosov map in standard form. If $x$ is a $p$-pronged periodic
point of period $n$ then the following holds.
\begin{numlist}
\item If $f$ preserves orientation then $\ind{f^{pn}}{x}=1-p$.
\item If $f$ reverses orientation and $n$ is odd
then $\ind{f^{2n}}{x}=1-p$.
\item If  $x$ is a regular point then $\ind{f^{2n}}{x}=1-2=-1$.
\end{numlist}
\end{lem}
\begin{proof}
Observe that $x$ is a fixed point of $f^n$. Moreover, $f^n$ is orientation-reversing if
and only if $f$ is orientation-reversing and $n$ is odd. Suppose first that $f$ preserves
orientation. By \propc{dash}, $x$ is of type $(p,0)^+$ for $f^{np}$ and, by \propc{pa},
$\ind{f^{np}}{x}=1-p$.  On the other hand, if $f$ reverses orientation and $n$ is odd, by
\propc{dash},  $x$ is of type $(p,0)^+$ for $f^{2n}$. Hence, as before,
$\ind{f^{2n}}{x}=1-p$. Finally, the third statement follows from the first two.
\end{proof}

\start{prop}{fix} Let $\Smap{f}{\Suu}$ be an orientation-preserving homeomorphism in
standard form. Assume that  $f$ has  a pseudo-Anosov component with either an interior
fixed point or an invariant boundary component. Moreover, suppose that every prong
emanating from the fixed point or invariant boundary component is fixed under the action
of $f$. Then $f$ has a fixed-point class of negative index.
\end{prop}
\begin{proof} By \propc{pa},  an interior fixed point $x$ of a pseudo-Anosov component is a fixed-point class. Moreover, if  $p$ is the number of prongs emanating from $x$, then $x$ is of type $(p,0)^+$ because,  by hypothesis,  all the prongs emanating from $x$ remain fixed under the action of $f$.  By \propc{pa},
the index of $x$ is $1-p$.  Since $p \ge 2$, this index is negative.

Now suppose that $f$ has an invariant boundary component $B$ such that all prongs remain
fixed under the action of $f$.  By \remc{rotation}, $f|_B$ is a rotation which implies
that $f|_B=\id$. Therefore, $B$ is contained in a fixed-point class $C$. This class $C$
must be as in  \propc{types}(C.2), (C.3) or (F). In each of the three cases,
$\ind{f}{C}\le -p<0$, where $p$ is the number of prongs emanating from $B$.
\end{proof}

\start{lem}{laberinto}  Let $m$ and  $g$ be positive integers such that $g \ge 2$, and
let $\Smap{f}{\Su{g}}$ be a pseudo-Anosov map in standard form.
\begin{numlist}
\item If there exists $x \in \Su{g}$ such that $x$ is an isolated fixed point of $f^m$, $m$ is even if $f$ is orientation-reversing, and $\ind{f^m}{x} \ne 1$, then  $f^m$ has a fixed-point class of negative index.
\item If $L(f^m)<0$ then  $f^m$ has a fixed-point class of negative index.
\item If there exists a singularity $x$ such that $f^m(x)=x$ and each prong emanating from $x$ remains fixed under the action of $f^m$, then  $f^m$ has a fixed-point class of negative index.
\item If there exists a $p$-pronged singularity $x$ of period $n$, $n$ is even if $f$ reverses orientation, and $m=np$,  then  $f^m$ has a fixed-point class of negative index.
\item If $m$ is even, $m/2$ is odd, $f$ is orientation-reversing and there exists a singularity of period $m/2$,  then  $f^m$ has a fixed-point class of negative index.
\item  If $m$ is even and  there exists a regular point of period  $m/2$  then  $f^m$ has a fixed-point class of negative index.
\end{numlist}
\end{lem}
\begin{proof}
Observe that by \propc{pa}, isolated fixed points of pseudo-Anosov maps on closed
surfaces are fixed-point classes. Therefore, it suffices  to show that, in each case,
$f^m$ has a fixed point of negative index.

Suppose that the hypotheses of (1) hold. Let us denote by $p$ the number of prongs
emanating from $x$. By \propc{pa}, if $m$ is even or $f$ is orientation-preserving,
$\ind{f^m}{x} \in \{1,1-p\}$. Since $p \ge 2$ and $\ind{f^m}{x} \ne 1$ we have that
$\ind{f^m}{x}=1-p \le -1$.

If the hypothesis of  (2) holds, then the conclusion follows from \theoc{essential}.

The statement  (3)  is a consequence of \propc{fix} and the statements (4), (5) and (6)
are consequences of \lemc{tango}.
\end{proof}

\start{lem}{espejos}  Let $g \ge 2$, let $f \in \hg$ be a pseudo-Anosov map in standard
form, and  let $k$ be the number of $f$-orbits of singularities of the foliation on
$\Su{g}$.
\begin{numlist}
\item If $k \ge 3$ then  there exists a positive integer $m$ such that $m \le 4g-4$ and  $f^m$ has a fixed-point class of negative index.
\item If $k=2$, $n_1 \le n_2$ and $p_2 \ge 4$ where $n_1, n_2$ are the periods of the  orbits of singularities and $p_2$ is the number of prongs emanating from each  point in the orbit of period $n_2$  then there exists a positive integer $m$ such that $m \le 4g-4$  and  $f^m$ has a fixed-point class of negative index.
\end{numlist}
\end{lem}
\begin{proof}
Let us denote by $O_1, O_2, \dots, O_k$ be the $f$-orbits of singularities. For each
$\ii{i}{k}$, let $n_i$ be the period of $O_i$, and let $p_i$ be the number of prongs
emanating from each element of $O_i$. Hence, we may write the Euler-Poincar\'{e} Formula
\refeq{ep}  in the form, $$ \sum_{i=1}^kn_i(p_i-2)=4(g-1). $$ Assume that  $k \geq 3$.
Since  $p_i \geq 3$, for each $i$,
\begin{equation}
\sum_{i=1}^k n_i  \leq \sum_{i=1}^k n_i(p_i-2)=4(g-1).\label{mirror}
\end{equation}
Let   $\ii{i_0}{k}$ be such that $n_{i_0}p_{i_0}=\min_{1 \leq i \leq k}\{n_i p_i\}$.
Then, $$ 3n_{i_0}p_{i_0} \leq  \sum_{i=1}^k n_ip_i=4(g-1)+2\sum_{i=1}^kn_i \leq 12(g-1),
$$ i.e., $n_{i_0}p_{i_0} \le 4(g-1).$ If $n_{i_0}$ is even or $f$ is
orientation-preserving then \lemc{laberinto}(4) holds for $m=n_{i_0}p_{i_0}$. If
$n_{i_0}$ is odd and $f$ is orientation-reversing then \lemc{laberinto}(5) holds for
$m=2n_{i_0}$. Hence, (1) is proved.

Let us prove (2). By \refeq{mirror}, $$n_1p_1 \le n_1(p_1-2)+2n_2 \le
n_1(p_1-2)+n_2(p_2-2)=4(g-1).$$ As in the preceding paragraph, this gives a proof of (2).
\end{proof}

\start{lem}{tiger} Let $m$, $n$, $p$, $g$ be positive integers such that $n \ge 2$ and $g
\ge 2$, and let $f \in \hg$ be a reducible map in standard form.
\begin{numlist}
\item Suppose  there exist a pseudo-Anosov $f$-component $C$ and a $p$-pronged boundary component $B$ of $C$ such that, $f^n(B)=B$, $pn=m$, and $n$ is even or $f$ is orientation-preserving. Then $f^m$ has a fixed-point class of negative index.

\item   Suppose there exist a pseudo-Anosov $f$-component $C$,  and a boundary component $B$ of $C$ such that $f^n(B)=B$, $n$ is odd, $2n=m$, and $f$ is orientation-reversing. Then $f^m$ has a fixed-point class of negative index.

\item  Suppose there exist a pseudo-Anosov $f$-component $C$,  and $x \in \Int(C)$, such that $f^n(x)=x$, $n$ is odd, $2n=m$, and $f$ is orientation-reversing. Then  $f^m$ has a fixed-point class of negative index.

\item Suppose there exist a pseudo-Anosov $f$-component $C$, and a $p$-pronged boundary singularity $x \in C$ such that,  $f^n(x)=x$, $pn=m$, and $n$ is even or $f$ is orientation-preserving. Then $f^m$ has a fixed-point class of negative index.

\item Suppose there exists a finite-order $f$-component $C$ such that $f^m|_C=\id_C$. Then $f^m$ has a fixed-point class of negative index.
\end{numlist}
\end{lem}
\begin{proof}
The statements (1), (2), (3) and (4) are consequences of \lemc{fix} and \propc{dash}.

If there exists a finite-order component $C$ such that $f^m|_C=\id_C$ then $C$ is
included in a fixed-point class $\widetilde{C}$ of $f^m$. Recall that by definition, a
component of a reducible map has negative Euler characteristic. Then, by
\propc{types}(F), $\ind{f^m}{\widetilde{C}} \le \chi(C)<0$.  Hence, (5) is proved.
\end{proof}

\begin{rem}
Despite the title of this section, these results (and their consequences) hold for maps
in Thurston canonical form. Indeed,  observe that if we have a fixed point, or an
invariant boundary component, of a pseudo-Anosov map in  Thurston canonical  form
$\Smap{f}{\Suu}$, we can define type in the same way as we did for maps in  standard
form: A $p$-pronged fixed point or a $p$-pronged invariant boundary component will be of
type $(p,k)^+$ (resp. $(p,k)^-$) if $f$ preserves (resp.  reverses) orientation and $f$
acts as the map  $r^+_{(p,k)}$ (resp. $r^-_{(p,k)}$) on the prongs emanating from it. By
\lemc{class}, there exists a pseudo-Anosov homeomorphism $\Smap{g}{\Suu}$ in standard
form isotopic to $f$.  The isotopy between these two maps preserves the types of fixed
points and invariant boundary components. Since, by \lemc{jiang}, indices of fixed-point
classes are preserved under isotopy, the index of a fixed point, or of an invariant
boundary  component, of any pseudo-Anosov map can be calculated as a function of its type
according to Tables \refc{indicespA} and \ref{indicesbdpA}. However,  the more
restrictive statements for homeomorphisms in standard form are sufficient for our
purposes.
\end{rem}

\section{Reducible maps}\label{sporm}

\markboth{THURSTON-NIELSEN THEORY}{REDUCIBLE MAPS}

This section is devoted to the study of reducible maps. In Subsection~\ref{corm} we state
some properties of a system of reducing curves. In Subsection~\ref{opc} (resp. \ref{orc})
we determine some properties of finite-order components of orientation-preserving (resp.
orientation-reversing) reducible maps.

\subsection{Components of reducible maps}\label{corm}

We begin with a basic property.

\start{rem}{ab} If $A$ and $B$ are surfaces then $$ \chi(A \cup B)=\chi(A)+\chi(B)-\chi(A
\cap B), $$ see \cite[Corollary V.4.6 and Proposition V.5.8]{Do}.
\end{rem}

By using the preceding remark, we shall prove the next result, which belongs to the class
of well-known facts whose proof is hard to find in the literature. Before stating it, let
us introduce the following notation: For each $\Suu$, denote by $\genus(\Suu)$ (resp.
$\bc(\Suu)$) its genus (resp. its number of boundary
components).\index{1g@$\genus(\Suu)$} \index{1bc@$\bc(\Suu)$}

\start{lem}{euler} Let $\Gamma=\{\Gamma_1,\Gamma_2,\ldots,\Gamma_n\}$ be a finite set of
pairwise disjoint simple closed curves  in a surface $\Suu$ and let $\Suu^1, \Suu^2,
\ldots, \Suu^k$ denote the closure of the connected components of $\Suu \setminus
\Gamma$. Then $\chi(\Suu)=\sum_{i=1}^k\chi(\Suu^i)$. Further,
$$\sum_{i=1}^k\genus(\Suu^i) \le \genus(\Suu).$$
\end{lem}

\begin{proof} We will prove  by induction that if $j \le k$, then
$$ \chi(\cup_{i=1}^j\Suu^i)=\sum_{i=1}^j\chi(\Suu^i). $$ That equality holds when $j =1$
is trivial. We now assume that the formula holds for some positive integer $j$ such that
$j <k$. The intersection of $\cup_{i=1}^j \Suu^i$ with $\Suu^{j+1}$ is a finite (possibly
empty) union of pairwise disjoint simple closed curves.  In any case, the Euler
characteristic of such an intersection is $0$, so by \remc{ab} and the inductive
hypothesis, $$ \chi\left(\bigcup\limits_{i=1}^{j+1}
\Suu^i\right)=\chi\left(\bigcup\limits_{i=1}^{j} \Suu^i\right)
+\chi\left(\Suu^{j+1}\right)-\chi\left(\left(\bigcup\limits_{i=1}^j \Suu^i\right) \cap
\Suu^{j+1}\right) $$ $$ =\sum_{i=1}^j\chi\left(\Suu^i\right)+\chi\left(\Suu^{j+1}\right)=
\sum_{i=1}^{j+1}\chi\left(\Suu^i\right) $$ as desired.

We prove  now that $\sum_{i=1}^k\genus(\Suu^i) \le \genus{\Suu}$. The surface
$\Suu=\cup_{i=1}^k\Suu^i$ is connected, therefore, without loss of generality, we can
assume that for each $\ii{i}{k-1}$, $$\Suu^{i} \cap \Suu^{i+1} \ne \emptyset.$$ Since
$\genus(\cup_{i=1}^k\Suu^{i})=\genus(\Suu)$, it suffices to show that
\begin{equation}
\sum_{i=1}^j\genus(\Suu^i) \le \genus(\bigcup\limits_{i=1}^j\Suu^{i})\label{field}
\end{equation}
for each $\ii{j}{k}$.

Clearly, \refeq{field}, holds for $j=1$. Let $j$ be a positive integer such that $j \le
k$ and assume \refeq{field} holds for each $\ii{i}{j-1}$. Denote by  $C_j$ the number of
connected components of $\left(\cup_{i=1}^j\Suu^{i}\right) \cap \Suu^{j+1}$. Observe that
our assumption implies that $C_j \ge 1$ for each $\ii{j}{k-1}$. By \remc{ab}, since
$\chi\left((\cup_{i=1}^j\Suu^{i}) \cap\Suu^{j+1}\right)=0$, $$
\chi\left(\bigcup\limits_{i=1}^j\Suu^{i}\right)+\chi(\Suu^{j+1})=\chi\left(\bigcup\limits_{i=1}^{j+1}\Suu^{i}\right).
$$ Hence,
\begin{equation}
2\genus(\bigcup\limits_{i=1}^j\Suu^{i})+\bc(\bigcup\limits_{i=1}^{j}\Suu^{i})-2+2g_{j+1}
+b_{j+1}-2\label{olas}
\end{equation}
$$
=2\genus(\bigcup\limits_{i=1}^{j+1}\Suu^{i})+\bc(\bigcup\limits_{i=1}^{j+1}\Suu^{i})-2.
$$ Each boundary component of $\cup_{i=1}^j\Suu^{j+1}$ is either a  boundary component of
$\cup_{i=1}^{j}\Suu^{i}$ or a boundary component of $\Suu^{j+1}$. Further,
$\cup_{i=1}^{j+1}\Suu^{i}$ has  exactly $\bc(\cup_{i=1}^{j}\Suu^{i}) - C_j$ boundary
components lying in $\cup_{i=1}^j\Suu^{i}$ and exactly $b_{j+1} -C_j$ boundary components
lying in $\Suu^{j+1}$. Thus, $\bc(\cup_{i=1}^{j+1}\Suu^{i})=\bc(\cup_{i=1}^{j}\Suu^{i}) -
C_j+b_j -C_j=\bc(\cup_{i=1}^{j}\Suu^{i})+b_j-2C_j$. Now, substituting this equality in
\refeq{olas} we obtain $$
2\genus(\bigcup\limits_{i=1}^j\Suu^{i})+2g_{j+1}+2C_{j}-2=2\genus(\bigcup\limits_{i=1}^{j+1}\Suu^{i}),
$$ so
$$\genus(\bigcup\limits_{i=1}^j\Suu^{i})+g_{j+1}+C_{j}-1=\genus(\bigcup\limits_{i=1}^{j+1}\Suu^{i}).
$$ Since $C_j \ge 1$ the result follows from the inductive hypothesis.
\end{proof}

We introduce the following notation which will be used frequently.

\start{nota}{dan} Let $g \ge 2$, let $\Smap{f}{\Su{g}}$ be a reducible homeomorphism, let
$\Gamma$ be a system of invariant curves, and let $N(\Gamma)$ be an invariant tubular 
neighborhood for $\Gamma$. An {\em   $f$-transversal} \index{transversal,
$f$-}\label{transversal} is a subset  $\{C_1,C_2,\ldots,C_k\}$ of $f$-components  such
that the set of all  $f$-components  is the disjoint union of the $f$-orbits of the
$C_i$'s. For each $1 \le i \le k$, we denote by $g_i$, $b_i$, and $n_i$ the genus, the
number of boundary components, and the period of $C_i$, respectively. If
$\Smap{f^{n_i}|_{C_i}}{C_i}$ is finite-order then, to simplify notation, we denote
$\sfi$ by $\sigma_i$.
\end{nota}

\start{lem}{corep} With \notc{dan} the following hold.
\begin{numlist}
\item $\sum_{i=1}^k (2g_i+b_i-2)n_i=2g-2.$
\item For each $\ii{i}{k}$, $2g_i+b_i-2 \ge 1$.
\item  For each $\ii{i}{k}$,  $n_i \le n_i(2g_i+b_i-2) \le 2g-2$.
\item For each, $\ii{i}{k}$, if $g_i=0$  then $b_i \ge 3$.
\end{numlist}
\end{lem}

\begin{proof} Observe that the set of all $f$-components is $\{f^j(C_i)\}_{1 \le j \le n_i}$, $1 \le i \le k$ and that, since $f$ is a homeomorphism, for each $\ii{i}{k}$ and each $\ii{j}{n_i}$, $\chi(C_i)=\chi(f^j(C_i))$. Thus  (1)  follows from  \lemc{euler}.

By the definition of standard form, $2-2g_i-b_i=\chi(C_i) \le -1$, hence (2) holds.

Clearly, (3) is a consequence of (1) and (2), and (4) a consequence of (2).
\end{proof}

Now we study transversals with only one component.

\start{lem}{ciclo} With \notc{dan}, if there exists an $f$-transversal containing only
one component, $C_1$, and $n_1\ge 2$,  and the  boundary components of $C_1$ form a
cycle under the action of $f^{n_1}$, then $n_1=2$ and $g=2g_1+b_1-1$.
\end{lem}
\begin{proof}
Let $B_1$ be a boundary component of $C_1$. By the hypothesis, the boundary components of
$C_1$ are $\{f^{kn_1}(B_1)\}_{k \in \{0,1,\dots, b_1-1\}}$.

Let $A \subset N(\Gamma)$ be the annulus which has $B_1$ as a boundary component and let
$C$ be an $f$-component such that $C \cap \Cl(A) \ne B_1$ and $C \cap \Cl(A)  \ne
\emptyset$. Since $f$ is a homeomorphism, and $C=f^j(C_1)$ for some $\ii{j}{n_1}$, and
the boundary components of $C_1$ form a cycle under the action of $f^{n_1}$, so do the
boundary components of $C$.  Therefore,  each boundary component of $C$ is of the form
$f^{n_1k}(\Cl(A) \cap C)$ for some $k \in \{0,1,\dots, b_1-1\}$.

Observe that $\bigcup \limits_{k=0}^{b_1-1}f^{n_1k}(A)  \cup C \cup C_1$ is a closed
subsurface of $\Su{g}$. Therefore, $ \cup_{k=0}^{b_1-1}f^{n_1k}(A)  \cup C \cup
C_1=\Su{g}$. If  $C=C_1$ then $n_1=1$. Since this is impossible, $n_1=2$, and, by
\lemc{corep}(1), $g=2g_1+b_1-1$ and $n_1=2$, as desired.
\end{proof}

If $f \in \fgb$ then $\sf=\sigma_{\tf}$. Then, by  \theoc{www} we have the following
result.

\start{cory}{wwwc} If $g \ge 2$ and $f \in \fgbp$ (resp. $\fgbr$)  then $\sf \le 4g+2$
(resp. $4g+(-1)^g4$).
\end{cory}

\subsection{The orientation-preserving case}\label{opc}

We begin with a basic property of finite-order maps of $\Su{0}$.

\start{lem}{s0} If $f \in \fgpx{0}$ is not the identity, then $f$ has exactly two fixed
points, $x_1$ and $x_2$. Therefore,  the $f$-period of every $x \in \Su{0} \setminus
\{x_1, x_2\}$ is $\sf$.
\end{lem}
\begin{proof} By \refeq{anne}, for every $m \in \N$, $L(f^m)=2- \tr \fs^m$. Since $H_1(\Su{0})$ is trivial, $\tr \fs^m=0$, so $L(f^m)=2$.  Now, the result follows by \theoc{suma} and \lemc{jiji}.
\end{proof}

As a consequence we have the following.

\start{cory}{s0b} Let $f \in \fgbpx{0}{b}$. If $\sf \ge 2$, and $b \ge 3$, then there
exists a positive integer $j$ and $a \in \{0,1,2\}$ such that $b=j \sf + a$. In
particular $b \ge \sf$. Moreover, if $\sf=b=3$ then the three boundary components of
$\Suf{0}{b}$ form a cycle under the action induced by $f$.
\end{cory}
\begin{proof} Consider the induced map $\Smap{\tf}{\Su{0}}$. By \lemc{s0}, each $\tf$-orbit has $1$ or $\sf$ elements and there are exactly two $\tf$-orbits of only one element.  The set of collapsed boundary components is $\tf$-invariant, so it is a disjoint union of orbits of $\tf$. Thus, $b=j\sf+a$,  where $j \in \N$ and $a \in \{0,1,2\}$. Clearly, if $b=\sf=3$, then $a=0$, $j=1$ and the three collapsed boundary components  form an $\tf$-periodic orbit. Therefore, the boundary components of $\Suf{0}{b}$ also form a cycle under the action of $f$.
\end{proof}

\start{lem}{numbertoro} Let $\Smap{f}{\Suf{1}{b}}$ be an orientation-preserving
finite-order map  such that  $\sf>2b$. Then $\sf \in \{3, 4,6\}$ and $b \in \{1,2\}$.
Moreover, if $b=2$, then the two boundary components of $\Suf{1}{b}$ are interchanged
under the action of $f$.
\end{lem}
\begin{proof}  Consider the induced map $\Smap{\tf}{\Su{1}}$. If $\tf$ is of type $\typeop{n}{\,\,}$  for some positive integer $n$, then $\Per(\tf)=\{\sigma_{\tf}\}$. With arguments analogous to the used in the proof of \coryc{s0b}, we see that $b=j\sf$ for some positive integer $j$. In particular, $b \ge \sf$. Since this contradicts our hypotheses, $\tf$ cannot be of type $\typeop{n}{\,\,}$ and the result follows from  \lemc{ordertoro}.
\end{proof}

\start{lem}{one fixed} If $g \ge 1$ and $f \in \fg$.
\begin{numlist}
\item If   $\sf >4g$ then $f$ has at most one fixed point. Consequently, if $h \in \hgbpx{g}{1}$ (resp. $h \in \hgbpx{g}{2}$) is such that $\widetilde{h} \in \fgp$ and $\sigma_h >4g$, then $h$ does not have fixed points (resp. both boundary components of $\Suf{g}{2}$ are interchanged under the action of $h$.).
\item If $\sf>2$ and $f$ has an isolated fixed point then $f$ is orientation preserving. Consequently, if $h \in \hgbrx{g}{1}$  and $\widetilde{h} \in \fgr$ then $\sigma_h=2$.
\end{numlist}
\end{lem}
\begin{proof}

The second statement is a consequence of \lemc{jijior}, so let us prove the first. Assume
first that  $g=1$. Since $\sf>4$, by \lemc{ordertoro} $f$ is either of type
$\typeop{6}{1,2,3}$ or $\typeop{n}{\,\,}$ for some $n \ge 5$ and the result follows
directly. Therefore, we can assume that $g \ge 2$. Let $\typeop{\sf}{p_1, p_2, \dots,
p_R}$ be the type of $f$.  If
 $f$ has at least two fixed points, we can assume that $R \ge 2$ and $p_{R-1}=p_R=1$. Set $R_0=R-2$ and $T=-R_0+\frac{2g+\sum_{i=1}^{R-2}p_i}{\sf}$. By  \theoc{signatures}, $T$ is a non-negative integer. On the other hand, for each $\ii{i}{R-2}$, $p_i \le \frac{\sf}{2}$. Hence, since by hypothesis $\sf >4g$,
$$ R_0 \le \frac{2g+\sum_{i=1}^{R_0}p_i}{\sf} \le \frac{2g+R_0 \frac{\sf}{2}}{\sf}
<\frac{1+R_0}{2}. $$ Then $R_0 <1$, so $R_0=0$ and $T=\frac{2g}{\sf}$. Since $\sf>4g$,
$T$ is not integer, which is impossible.
\end{proof}

\start{lem}{items0}   With \notc{dan},  let $\ii{i}{k}$ be such that $C_i$ is a
finite-order component.
\begin{numlist}
\item If $\sigma_i  \le 4g_i$  then $n_i \si \le 4g$.
\item If $\sigma_i \le 2$ then $n_i \si  \le 4g-4$.
\end{numlist}
\end{lem}
\begin{proof}
If $\sigma_i  \le 4g_i$, by \lemc{euler},  $\si n_i \le 4 g_in_i \le 4g$. Therefore, (1)
holds.

If $\si \le 2$, (2) follows directly from \lemc{corep}(3).
\end{proof}

\start{lem}{items}   With \notc{dan},  let $\ii{i}{k}$ be such that $C_i$ is a
finite-order component of $f$ and $f^{n_i}|_{C_i}$ is orientation preserving.
\begin{numlist}
\item  If $g_i=0$, and  $b_i  \ne 3$ or $\si \ne 3$ or $n_i \le \frac{4}{3}(g-1)$  then   $n_i \si \le 4g-4$.
\item  If $g_i=1$, and $b_i \ge 3$ or $\si  \ge 7$ then $n_i \si \le 4g-4$.
\item  If $g_i \ge 2$, and $b_i \ge 3$  then $n_i \si \le 4g-4$.
\end{numlist}
\end{lem}
\begin{proof}
Suppose now the hypotheses of (1) hold. By \lemc{items0}(2) we can assume that $\si \ge
3$. By  \lemc{s0b},  $3 \le \si  \le b_i$. If $b_i=3$,  then $\si=3$, $n_i \le
\frac{4}{3}(g-1)$ and the result follows directly. If $b_i \ge 4$, by \lemc{corep}(3),
$2n_i \le (b_i-2)n_i \le 2g-2$. Therefore, $n_i \si \le n_i(b_i-2)+2n_i \le 4g-4$, which
completes the proof of (1).

Now, we prove (2).  If  $\si \le 2b_i$ then by \lemc{corep}(3), $n_i\si \le 2n_ib_i \le
4g-4$ as desired. If $\si >2b_i$, the result holds by \lemc{numbertoro}.

Finally, we prove (3).  By \coryc{wwwc}, $\si \le 4g_i+2$. If $b_i \ge 3$,  by
\lemc{corep}(3), $$ n_i \si \le 2 n_i (2g_i+1) \le 2n_i(2g_i+b_i-2) \le 4g-4, $$ as
desired.
\end{proof}

\subsection{The orientation-reversing case}\label{orc}

The aim of the present subsection is to obtain analogous results to those of the previous
section for  finite-order components of orientation-reversing reducible maps.

\start{rem}{order even} If $\Smap{f}{\Suu}$ is a finite-order orientation-reversing map
then its order is even.
\end{rem}

\start{lem}{s0bro} If $b \ge 3$ and $f \in \fgbrx{0}{b}$ then there exist integers $j$
and $a$ such that $j \ge 1$, $a \in \{0,2\}$, and $b=j\sf/2+a$.  In particular, $b \ge
\sf/2$. Moreover, if $b=\sf/2$, then $\sf/2$ is odd and the $b$ boundary components of
$\Suf{0}{b}$ form a cycle under the action of $f$.
\end{lem}

\begin{proof} Consider the induced map $\Smap{\tf}{\Su{0}}$. Clearly, $\tf^m$ preserves orientation if  $m$ is even, and reverses orientation if $m$ is odd. Since $H_1(\Su{0})$ is trivial, by \refeq{anne},
\begin{eqnarray*}
L(\tf^m) & = &\left\{\begin{array}{ll}
        0  & \mbox{if $m$ is odd},\\
        2 & \mbox{if $m$ is even.}
    \end{array}
    \right.
\end{eqnarray*}

That the result holds if $\sf=2$ is trivial. Therefore,  we can assume without loss of
generality that $\sf>2$. If $\Fix(\tf) \ne \emptyset$ then, by \lemc{jijior},  $f$  has a
pointwise-fixed simple closed curve $\gamma$ and, $\tf$ acts as a reflection in a neighborhood of
$\gamma$. Then, $\tf^2=\id$ and $\sigma_{\tf}=\sf=2$, a contradiction. Then, we can
assume that $\Fix(\tf)=\emptyset$.

Now, by \remc{order even}, $h=\tf^2$ is an orientation-preserving map such that
$\sigma_{h}=\sf/2$ and $L(h^m)=2$ for every $m \in \N$. By \lemc{s0}  it follows that $h$
has two fixed points which are the only  points of $h$-period strictly less than $\sf/2$.
Since $\Fix(\tf)=\emptyset$,  we can conclude that $\tf$ has a two periodic orbit, which
is the only periodic orbit of even period strictly less than $\sf$.

By \lemc{parageo}, if there exists a simple closed curve  $\gamma \subset \Fix(f^i)$ then $i$ is odd
and $i=\sf/2$. In any case, $\Per(\tf) \subset \{2,\sigma/2, \sigma\}$ and there is
exactly one orbit of period $2$. Since the set of collapsed boundary components is
$\tf$-invariant, and has cardinal $b$, $b=k \sigma/2+a$ for some $a \in \{0,2\}$ and some
$k \in \N$. Since $b \ge 3$, $k \ge 1$.

Now, observe that if $b=\sf/2$, since $b \ge 3$, the set of collapsed boundary components
must consist of an $\tf$-periodic orbit of period $\sf/2$. Moreover, by \lemc{parageo},
$\sf/2$ is odd. Clearly, in this case, the boundary components of $\Smap{f}{\Suf{0}{b}}$
form an $f$-cycle.
\end{proof}

\start{lem}{s1bro} If $f \in \fgbrx{1}{b}$ then there exist a positive integer $j$ such
that $b=j\sf/2$. In particular, $b \ge \sf/2$.
\end{lem}
\begin{proof} Consider the induced map $\Smap{\tf}{\Su{1}}$. The characteristic polynomial of $\ffs$ is
 $x^2-\tr (\ffs) x-1$, whose roots are
 $$\lambda_1=\frac{\tr(\ffs)+\sqrt{\tr(\ffs)^2+4}}{2}$$
 and
$$ \lambda_2=\frac{\tr(\ffs)-\sqrt{\tr(\ffs)^2+4}}{2}.$$ Since $\tf$ is finite-order,
$\lambda_1$ and $\lambda_2$ are roots of unity. Clearly, $\lambda_1, \lambda_2 \in \R$,
and are not equal, so $\{\lambda_1, \lambda_2\}=\{-1,1\}$. Therefore,
\begin{eqnarray*}
\tr(\ffs^m)&= &\left\{\begin{array}{ll}
        2  & \mbox{if $m$ is even},\\
        0 & \mbox{if $m$ is odd.}
    \end{array}
    \right.
\end{eqnarray*}
Hence, by \refeq{anne}, $L(\tf^i)=0$ for every $i \in \N$. By \lemc{jiji}, $\tf$ does not
have periodic points of even period strictly less than $ \sf$.

By \lemc{jijior} and \remc{order even}, the periodic points of odd period are contained
in simple closed curves. With an argument similar to that used in the proof of \lemc{s0bro}, we can
show that if for some $i <\sf$ there exists a simple closed curve  $\gamma \subset \Fix(\tf^i)$, and
$i$ is minimal with this property, then $i$ is odd and $i=\sf/2$.

Then $\Per(\tf) \subset \{\sf, \sf/2\}$ and we can complete the proof as we did in
\lemc{s0bro}.
\end{proof}

\begin{rem} We could have proved \lemc{s1bro} in the same way we proved \lemc{numbertoro}, but the argument would have been much longer.
\end{rem}

\start{lem}{itemsor}   With \notc{dan}, let $f$ be orientation-reversing, and  let
$\ii{i}{k}$ be such that $f_{C_i}$ is a finite-order orientation-reversing component of $f$.
\begin{numlist}
\item If  $\sigma_i \le 3$ then $n_i \si \le 4g-4$.
\item  If $g_i=0$,  and  $b_i$ is even, or  $b_i \ne \sigma_i/2$ or $k =1$ then $n_i \si \le 4g-4$.
\item  If $g_i=1$ then $n_i \si \le 4g-4$.
\item If $g_i \ge 2$ and  $b_i \ne 2$  then $n_i \si \le 4g-4$.

\end{numlist}
\end{lem}
\begin{proof}
Assume that the hypotheses of (1) hold. By \remc{order even}, $\si$ is even, so $\si=2$.
Hence, by \lemc{corep}(3), $\si n_i=2n_i \le 4g-4$ and the proof of (1) is complete.

Now, suppose that $g_i=0$. By \lemc{corep}(4), $b_i \ge 3$ and, by (1) we can assume that
$\si \ge 4$. If $\si \le b_i$, then, $b_i \ge 4$ and, by \lemc{corep}(3) $2n_i  \le
(b_i-2)n_i \le 2g-2$. Here, by \lemc{corep}(3) $$ \si n_i \le (b_i-2)n_i +2n_i \le 4g-4.
$$ Therefore, we can assume that $\si >b_i$. By \lemc{s0bro}, $b_i=\si/2+a$, for some
integer $a$ such that $a \in \{0,2\}$. If $a=2$, then, by \lemc{corep}(3) $$ \si
n_i=2(b_i-2)n_i \le 4g-4, $$ and the result holds. Hence, we can restrict ourselves to
the case $b_i=\si/2$. By \lemc{parageo}, $\si/2$ is odd. Then, by hypothesis, $k=1$. By
\lemc{corep}(1), $(b_i-2)n_i=2g-2$ and $b_in_i$ is even. Since this is impossible, (2)
holds.

Let us prove (3). Suppose that $g_1=1$. We claim that we can assume that
$\sigma_i>2b_i$. Indeed, if $\sigma_i \le 2b_i$ then, by \lemc{corep}(1), $$ n_i\sigma_i
\le 2n_ib_i=2(2g_i+b_i-2)n_i \le 2(2g-2), $$ a contradiction. So, the claim is proved. By
\lemc{corep}(3) there exists  a positive integer $k$, such that $b_i=k\sigma_i/2$. Hence,
$2b_i=k\sigma_i \ge \sigma_i$, which contradicts $\sigma_i >2b_i$.  So, the proof of (3)
is complete.

Finally, we prove (4). Observe that, by \coryc{wwwc}, $\sigma_i \le 4g_i+4$. Hence, if
$b_i \ge 4$, by \lemc{corep}, $$ \sigma_i n_i \le 2(2g_i+2)n_i \le 2(2g_i+b_i-2)n_i \le
2(2g-2), $$ as required. Now, if $b_i \in \{1,3\}$ and $h $ denotes the map induced by
$f^{n_i}|_{C_i}$ on $\Su{g_i}$, $h$  is an orientation-reversing finite order map with a
fixed point or a three periodic orbit. By Lemmas~\ref{jijior} and \ref{parageo},
$\sigma_h \in \{2, 6\}$. Observe that $\sigma_h=\sigma_i$. In particular $\sigma_i \le
2g_i+2$ and we can complete the proof as in the preceding paragraph.
\end{proof}

\start{lem}{sandalo} Let $g \ge 2$ and let $f \in \hgr$ be a reducible map in standard
form such that each of its components is finite-order. If the $f$-period of each
connected component of $N(\Gamma)$ is odd then $f$ is finite-order.
\end{lem}
\begin{proof} Let $A$ be a connected component of $N(\Gamma)$ and  denote  its $f$-period by $n$. Since $n$ is odd, $f^n|_{A}$is orientation-reversing. Hence, from the description of the standard form for the tubular neighborhood $N(\Gamma)$ it follows that $f^n|_{A}$ is conjugate to a map  $\Smap{\psi}{\SI \times I}$ of one of the following forms.
\begin{numlist}
\item $(z,t) \mapsto (z e^{2a\pi i},1-t)$, where $a$ is a rational number.
\item  $(z,t) \mapsto (\overline{z},t)$.
\end{numlist}
In particular,  $f^n|_{A}$ is finite-order. Since $N(\Gamma)$ is $f$-invariant,
$f|_{N(\Gamma)}$ is finite-order. Since $f|_{\Su{g} \setminus N(\Gamma)}$ is finite-order
by hypothesis, the lemma is proved.
\end{proof}

\part{Proofs of the main results}\label{inmemoriales}

\chapter{Bounds for minimum periods}\label{S9}

The aim of this chapter is to obtain bounds for the minimum periods of maps of surfaces
with boundary. In Section~\ref{eucaliptos}  we apply some of the techniques developed in
Chapters~\ref{S2}  and \refc{S7} in order to find upper bounds. In
Section~\ref{esferas}, we construct examples of maps with ``large" minimum periods. and,
by means  of these examples, we prove the existence of some lower bounds.
Section~\ref{mp} is dedicated to the study of the minimum periods for classes of
finite-order maps of surfaces with boundary. Theorems~B, C, D and E are proved there.

\section{Upper bounds for minimum periods}\label{eucaliptos}
\markboth{SURFACES OF LOW GENUS}{UPPER BOUNDS FOR MINIMUM PERIODS} The two main results
of this section are Proposition~A, which was stated in the introduction
(Chapter~\ref{S1}), and \propc{razon}. We begin with the proof of  Proposition~A which is
based in some of the ideas developed by Nielsen in \cite{N}, and uses Lefschetz's
Fixed-Point Theorem and some elementary algebra.

\begin{proclama}{Proposition A}\index{Proposition A}
If $2g+b\geq 4$ then $\mmm(\hgb)\leq 2g+b-2$.
\end{proclama}

\begin{proof} Suppose that  $f \in \hgb$ is a  homeomorphism such that $f, f^2 \ldots, f^{2g+b-3}$ are fixed-point free. In particular,  $2g+b-3 \ge 1$, so $f$ is fixed-point free. Since $\Fix(f)=\Fix(f^{-1})$, $f^{-1}$ is also fixed-point free. Then, by \theoc{lef},
$$L(f^{-1})=L(f)=L(f^2)=\ldots=L(f^{2g+b-3})=0,$$ so, by \refeq{anne},
$$\tr(\fs^{-1})=\tr(\fs)=\tr (\fs^2)=\ldots=\tr (\fs^{2g+b-3})=1.$$ (Recall that $\fs$
denotes the linear map induced by $f$ on the first homology group.) The dimension of
$H_1(\Suf{g}{b}, \Q)$, the first Betti number, is $2g+b-1 \ge 3$.  So, by \propc{paratraza}, $\tr (\fs^{2g+b-2}) \neq 1$ and, by
\refeq{anne}, $L(f^{2g+b-2}) \neq 0$. Thus, by \theoc{lef}, $f^{2g+b-2}$ has a fixed
point.
\end{proof}

In order to prove \propc{razon} we need to introduce some notation and prove an
elementary lemma.

If $\Fix(f^m)$ is a finite set, we denote by $P_m(f)$ \index{1p@$P_m(f)$} the number of
periodic points of $f$-period $m$, that is, $$\Card(\left\{x \in \Fix(f^m) \,\,: \,\,x
\notin \Fix(f^k)\,\, \mbox{ for $k \in \{1,2,\ldots,m-1\}$}\right\}).$$ Then $$
\Card(\Fix(f^m))=\sum_{d \vert m}P_d(f), $$ where $\sum_{d \vert m}$ denotes the sum over
all positive divisors of $m$. Consider the {\em Mobius function} \index{M\"{o}bius function}
$\map{\mu}{\N \setminus \{0\}}{\{-1,0,1\}}$ defined by

\begin{eqnarray*}
\mu(m) & = &\left\{\begin{array}{ll}
                1  & \mbox{if $m=1$},\\
                0 & \mbox{if $k^2 \vert m$ for some $k \ge 2$},\\
                        (-1)^r & \mbox{if $m=p_1p_2\ldots p_r$ distinct prime factors.}
        \end{array}
        \right.
\end{eqnarray*}
By the M\"{o}bius Inversion Formula, see \cite[Proposition 13.B.3]{Ch}
\begin{equation}
P_m(f)=\sum_{d \vert m}\mu(d)\Card(\Fix(f^{m/d})).\label{elena}
\end{equation}

Let $f \in \hg$. For each positive integer $i$ set \index{1l@$l(f^i)$} $$ l(f^i) =
\sum_{d \vert i}\mu(d)L(f^{m/d}). $$

Also, for each positive integer $n$, we define $l_n(f)=(l(f),l(f^2),\dots,l(f^n))$.
\index{1l@$l_n(f)$}

\start{rem}{do} Dold \cite{Do1} proves that if  $Y$ is an ENR,  and $V \subset Y$ is an
open set, and $\map{f}{V}{Y}$ is a map and $n$ is a positive integer such that
$\Fix(f^n)$ is compact, then $n$ divides $l(f^n)$. This property is easily checked for
the class of maps we are going to consider. As Llibre remarks in \cite{Ll}, these
``modified'' Lefschetz numbers are interesting because, for many classes of maps, if
$l(f^m) \ne 0$ then $m \in \Per(f)$. In particular, this holds for finite-order maps,
provided that $m$ is strictly smaller than the order of the map (by \lemc{musica}), and
for maps $\tf$ induced by maps $f \in \hgb$ with ``large'' minimum period (by
\lemc{diversidad}).
\end{rem}

Recall that if $f \in \hgb$, then $\tf$ denotes the homeomorphism induced by $f$ on
$\Su{g}$. If the minimum period of a map $f$ is large enough, then $l(\tf^i)$ counts the
number of points of $\tf$-period $i$.

\start{lem}{diversidad} Let $f \in \hgb$. For each $i <\mmm(f)$, $l(\tf^i)=P_i(\tf)$.
\end{lem}
\begin{proof} Observe that,  for each $i<\mmm(f)$, by  \propc{index},  $\tf^i$ is a homeomorphism with finitely many fixed points, each having index one. So, by \theoc{suma}, $L(\tf^i)=\Card(\Fix(\tf^i))$.  Hence, by  \refeq{elena},
$$ P_i(\tf)=\sum_{d \vert i}\mu(d)\Card(\Fix(\tf^d))=\sum_{d  \vert
i}\mu(d)L(\tf^d)=l(\tf^i), $$ as desired.
\end{proof}

\start{rem}{isla} Let $i$ be a positive integer and let $f \in \hgb$. If  $\mmm(f)>i $
then, by \lemc{diversidad}, $l(\tf^i) \ge 0$.
\end{rem}

\start{conj}{xx} Let $g \ge 2$. By \remc{isla}, if $f \in \hgbp$, then $l(\tf^i) \ge 0$
for each $i <\mmm(f)$. Hence, $\mmm(\hgbp)$ is less than or equal to the natural number
$m_g^+$ defined as $$\max\{n \ge 1\,\,:\,\ \mbox{there exists $f \in \hgp$ such that
$l(f^i) \ge 0$ for each $i < n$}\}. $$ We conjecture that $m_g^+=4g+2$. Analogously, if
$m_g^-$ equals $$ \max\{n \ge 1\,\,:\,\ \mbox{there exists $f \in \hgr$ such that $l(f^i)
\ge 0$ for each $i < n$ }\}, $$ then $m(\hgbr) \le m_g^-$. We conjecture that
$m_g^-=4g+(-1)^g4$. The case $g=2$ of these conjectures follows from some results of
Chapter~\ref{S10}.
\end{conj}

Let $f \in \hgb$  and for each $i$, let $k_i$ denote the number of $f$-cycles of period
$i$ of boundary components of $\Suf{g}{b}$. By \lemc{diversidad} and \propc{index}, if
$i<\mmm(f)$, then $\tf^i$ must have exactly $i k_i$ fixed points, each  with index one;
hence, $l(\tf^i)=ik_i$. Now, consider the  vector $(k_1,2k_2,\dots,bk_b)$. Clearly,
$\sum_{i=1}^bik_i=b$. If $i$ is a positive integer such that $l(\tf^i) \ne i k_i$, then
$\mmm(f) \le i$, i.e., $i$ is an upper bound for $\mmm(f)$. This idea is formalized  in
\propc{razon}. To state this result, we need to define a map which gives the minimum of
these $i$.

Let $n$ be a positive integer. Denote by  $C_{b,n}$  the subset of  $\N^{n}$ consisting
of all the $n$-tuples of the form $(k_1,2k_2,\dots,n k_{n})$ such that  $\sum_{i=1}^n
ik_i=b$. We define  $\map{\alpha}{C_{b,n} \times \Z^n}{\N \cup \{\infty\}}$  in the
following way. For  each $(\overline{x},\overline{y}) \in C_{b,n} \times \Z^n$,
\index{1a@$\alpha(\overline{x},\overline{y})$}
\begin{eqnarray*}
\alpha(\overline{x},\overline{y})  &= &
                   \cases{\infty  & \mbox{if $\overline{x}=\overline{y}$,} \cr
                                    \min\{\ii{i}{n}\,\,:\,\, x_i \ne y_i\}     & \mbox{otherwise.}}\cr
\end{eqnarray*}

For  $f \in \hg$ we set $$\gamma(f,b) =\max_{\overline{k} \in C_{b,b}}\{
\alpha(\overline{k},l_b(\tf))\}.$$ \index{1g@$\gamma(f,b)$}

Now, for each map $f \in \hgb$, we can  give the promised upper bound for $\mmm(f)$,
depending  on $l_b(\tf)$ and on the action of $f$ on the boundary components of
$\Suf{g}{b}$. This proposition will be very useful in Chapter~\ref{S10}, where we will
need to bound the minimum periods of maps $f \in \hgb$ such that the sequence $l_b(\tf)$
is in a given finite subset of $\Z^b$.

\start{prop}{razon} Let $f \in \hgb$. For each $\ii{j}{b}$, let $k_j$ denote the number
of $f$-cycles of boundary components of period $j$. Then $$\mmm(f) \le
\alpha((k_1,2k_2,\dots,bk_b),l_b(\tf)).$$ Consequently, for each $\Smap{f}{\Suf{g}{b}}$,
$\mmm(f) \le \gamma(\tf,b)$.
\end{prop}
\begin{proof}  Let $\ii{j}{b}$. If $\mmm(f) > j$, by \lemc{diversidad}, $l(\tf^j)=P_j(\tf)$. On the other hand, by \propc{index}, every point of $\tf$-period $j$ must be a collapsed boundary component. Thus $P_j(\tf)=jk_j$. Since  $(k_1,2k_2,\dots,bk_b) \in C_{b,b}$, and
$$ \alpha((k_1,2k_2,\dots,bk_b),l_b(\tf)) \ge \mmm(f), $$ the result  holds.
\end{proof}

\section{Lower bounds for minimum periods}\label{esferas}
\markboth{SURFACES OF LOW GENUS}{LOWER BOUNDS FOR MINIMUM PERIODS} In this section we
construct examples of maps with large minimum periods on surfaces with boundary. The
basic idea of these constructions consists in considering a map $f$ of a closed surface
$\Suu$ with an $f$-invariant set $D$ which is a union of finitely many pairwise disjoint
open disks and consider the restriction $f|_{\Suu \setminus D}$. Clearly, $f \in \hgb$,
where $b$ is the number of connected components of $D$.

We begin with a lemma which asserts the existence of such an $f$-invariant set $D$ as in
the preceding paragraph, when $f$ is finite-order and there exists an $f$-invariant
finite set $F$ of cardinal $b$.

\start{lem}{diamante} Let $\Smap{f}{\Suu}$ be a finite-order map. If
$F=\{x_1,x_2\dots,x_b\}$ is a subset of $\Int{\Suu}$ which is  $f$-invariant  then there
exists a subset $D$ of  $\Suu$ such that
\begin{numlist}
\item $D$ is $f$-invariant and $\Cl(D) \subset \Int(\Suu)$.
\item There exist $b$ pairwise disjoint open disks $D_1, D_2, \dots D_b$ such that $D=\union{i=1}{b}D_i$ and  the center of $D_i$ is $x_i$ for each $\ii{i}{b}$.
\end{numlist}
Moreover, if  $C$ is a closed subset of $\Suu$ and $F \cap C=\emptyset$ we can assume
that $C \cap \Cl(D)=\emptyset$.
\end{lem}
\begin{proof} It is a consequence of Lemmas~\ref{jiji} and \refc{jijior}.
\end{proof}

\start{rem}{donde} Let $S$ be a (possibly non-connected) compact orientable surface and
let $\Smap{f}{S}$ be a finite-order map, that is, $f^n=\id_{S}$ for some positive integer
$n$. Then \lemc{diamante} holds replacing $\Suu$ by $S$.
\end{rem}

As mentioned before, the examples we shall construct are based on finite-order maps of
closed surfaces. In the following, we state a property of such maps.

\start{lem}{musica} If $f \in \fgp$  then, for each $i < \sf$, $P_i(f)=l(f^i).$ In
particular, $i$ divides $l(f^i)$.
\end{lem}
\begin{proof} By \lemc{jiji}. $\Fix(f^i)$ is a finite set consisting of points of index one. By
\theoc{suma},  $\Card(\Fix(f^i))=L(f^i)$, and we can now argue as in the proof of
\lemc{diversidad}.
\end{proof}

Consider a map $f \in \fgp$ and let $F \subset \Su{g}$ be the $f$-invariant set
consisting of  all the points whose $f$-period is strictly less than $\sf$. Let
$b=\sum_{i=1}^{\sf-1}l(f^i)$. By \lemc{musica},  $\Card(F)=b$. So, if $D$ is  as in
\lemc{diamante}, then $h=f|_{\Su{g}\setminus D} \in \fgbp$ and  $\mmm(h)=\sf$.

This construction applies only in the special case  $b=\sum_{i=1}^{\sf-1}l(f^i)$.
Nevertheless, in several cases, by modifying $f$ by means of an isotopy we can obtain
another map $h \in \hgp$ with an $h$-invariant set $F$ of a certain cardinal $b$ to which
we can apply \lemc{diamante} in order to find an $h$-invariant set consisting of $b$
disks,  and $h|_{\Su{g}\setminus D}$ will be the  desired map. This is done in
\propc{fulgor}.

Our next step will be to define an auxiliary map which will used to calculate the minimum
period of the   map obtained by the above procedure. Let $\sigma$ be a positive integer.
We define $\map{\beta}{C_{b,\sigma} \times \Z^\sigma}{\N \cup \{\infty\}}$  in the
following way. For  each $(\overline{x},\overline{y}) \in C_{b,\sigma} \times \Z^\sigma$,
\index{1b@$\beta(\overline{x},\overline{y})$}
 \begin{eqnarray*}
\beta(\overline{x},\overline{y})  &= &
                   \cases{\infty  & \mbox{if $\overline{x}=\overline{y}$,} \cr
                                    \min(\{\frac{x_i}{i}\,:\,\, x_i < y_i\}\cup \{x_i-y_i\,\,:\,\,x_i>y_i\}) & \mbox{otherwise.}}\cr
\end{eqnarray*}

The following result, in which we formalize the above procedure, will be applied in
\propc{idioma} to construct examples giving certain lower bounds for minimum periods.

\start{prop}{fulgor} Let $f \in \fgp$. Assume there exists $$\overline{k}=(k_1,2 k_2,
\dots, \sf k_\sf) \in C_{b,\sf}$$ such that  $k_j=0$ if $l(f^j)=0$.  Then there exists $h
\in \hgbp$ such that $\widetilde{h}$ is isotopic to $f$ and $$\mmm(h)= \min \{ \sf,
\beta(\overline{k}, l_\sf(f))     \}.$$ Furthermore, if, for each $\ii{j}{\sf}$,
$jk_j-l(f^j)$ is a multiple of  $\sf$ then $\widetilde{h}=f$, so  $h$ is finite-order.
\end{prop}
\begin{proof} By \lemc{musica}, for each $\ii{j}{\sf}$,
$$P_j(f)=l(f^j) \ge \min\{l(f^j), jk_j\}.$$ Then there exists an $f$-invariant subset of
$P_j(f)$ of cardinal $\min\{l(f^j), jk_j\}$.  Let us denote this subset by $F_j$. Set $F
= \cup_{j=1}^\sf F_j$, an $f$-invariant set.  Let $D_1$ be the $f$-invariant set as in
\lemc{diamante} corresponding to $F$.

For each  $\ii{j}{\sf}$ we  will  inductively define an $f$-invariant subset of $\Su{g}$,
$B_j$, and a map $\Smap{\vv_j}{B_j}$ in the following way: If $l(f^j) \ge jk_j$, we set
$B_j=\emptyset$. Otherwise, by \remc{isla}, $0 \le l(f^j)<jk_j$. So $k_j>0$ and, by
hypothesis, $l(f^j)>0$. By \lemc{musica},  $P_j(f) \ne \emptyset$.  Consider $x \in
P_j(f)$. By \lemc{jiji}, there exists a disk  with center $x$ where $f^j$ is conjugate to
a rotation of order $\sf/j$. Therefore,  there exists an open annulus $A \subset \Su{g}$
centered at  $x$ such that
\begin{numlist}
\item $f^j(A)=A$.
\item $\union{i=0}{j-1}f^i(A) \cap (D_1 \cup \union{i=1}{j-1} B_j)=\emptyset$.
\item $f^j|_A$ is conjugate to a rotation of order $\sf/j$.
\end{numlist}
Set $B_j=\union{i=0}{j-1}f^i(A)$.

Now, we have $jk_j-l(f^j)>0$, and we consider two cases: $jk_j-l(f^j)$ is a multiple of
$\sf$ or not.

In the former case  $jk_j-l(f^j)=a_j \sf$ for some positive integer $a_j$, and we proceed
as follows. Choose $a_j$ disjoint $f$-orbits, $O_1, O_2, \dots, O_{a_j} \subset B_j$ and
set $H_j =\union{i=1}{a_j}O_i$. Let $D_2^j \subset B_j$ be a set as given by
\lemc{diamante} taking $\Suu=B_j$, $F=H_j$ and $C=\Cl(D_1)$. Denote $\vv_j=f|_{B_j}$.

Finally, if $jk_j-l(f^j)$ is not a multiple of $\sf$, set $m_j=\frac{jk_j-l(f^j)}{j}$.
Since $m_j >0$, by \lemc{musica}, $m_j$ is a positive integer. Let $C \subset A$ be a
closed annulus with non-empty interior.  Then there exists a homeomorphism
$\Smap{\vv_j}{B_j}$ such that
\begin{numlist}
\item $\vv_j|_{\partial B_j}=f|_{\partial B_j}$ and $\vv_j$ is isotopic to $f|_{B_j}$ relative to $\partial B_j$,
\item $\vv_j^j|_{C}$ is conjugate to a rotation of order $m_j$,
\item $\mmm(\vv_j) = \min\{jm_j, \sf\}$.
\end{numlist}

Choose $x \in \Int(C)$ and denote its $\vv_j$-orbit by $H_j$.  Then,
$$\Card(H_j)=jm_j=jk_j-l(f^j).$$ (Observe that the step $j=1$ can done in the same way as
the step $j>1$). Let $D_2^j$ be the set defined by \lemc{diamante} and \remc{donde},
taking $S$ as $\cup_{i=1}^jf^i(C) $ and $f$ as $\vv_j|_{\cup_{i=0}^{j-1}f^i(C)}$ It is
easy to check that \lemc{diamante} holds for $\vv_j|_{\cup_{i=1}^jf^i(C)}$,  so there
exists a set $D_2^j \subset \cup_{i=1}^jf^i(C)$ satisfying $(1)$ and $(2)$ of
\lemc{diamante}.

Clearly,  $D_1 \cup \cup_{j=1}^{\sf}D_2^j$ is a union of pairwise disjoint disks.
Moreover, the number of connected components of $D_1 \cup \cup_{j=1}^{\sf}D_2^j$ is equal
to $\Card(F \cup \cup_{j=1}^{\sf}H_j)$. Since

$$ \Card(F \cup \union{j=1}{\sf}H_j)= \sum_{j=1}^\sf
\min\{l(f^j),jk_j\}+\sum_{jk_j-l(f^j)>0}jk_j-l(f^j)=\sum_{j=1}^\sf jk_j=b, $$ $S=\Su{g}
\setminus (\union{j=1}{\sf} D_2^j \cup D_1)$ is a surface of genus $g$ with $b$ boundary
components. Define $\Smap{h}{S}$ in the following way:

\begin{eqnarray*}
h(x) & = & \cases{f(x) & \mbox{if  $x \notin \union{j=1}{\sf} B_j$},\cr
                                 \vv_j(x)& \mbox{if $x \in B_j$ for some $\ii{j}{\sf}.$}}\cr
\end{eqnarray*}

Let us prove that $\mmm(h)= \min \{ \sf, \beta(\overline{k}, l_\sf(f))     \}.$ Observe
that $ \union{j=1}{\sf}B_j$ and $S \setminus  \union{j=1}{\sf}B_j$ are $h$-invariant.
Hence, $$ \Per(h)= \Per(h|_{\cup_{j=1}^{\sf}B_j}) \cup  \Per(h|_{S \setminus
\cup_{j=1}^{\sf}B_j}). $$ Then $$\min \Per(h)=\min \{\min\Per(h|_{\cup_{j=1}^{\sf}B_j}),
\min \Per(h|_{S \setminus  \cup_{j=1}^{\sf}B_j})\}.$$ By definition of $h$, $$ \min
\Per(h|_{\cup_{j=1}^{\sf}B_j})=\min\left(\{ \sf \} \cup
\left\{jk_j-l(f^j)\,\,:\,\,jk_j>l(f^j), 1 \le j < \sigma_j \right\}\right). $$ On the other hand,
$h|_{S \setminus \cup_{j=1}^{\sf}B_j}=f|_{S \setminus  \cup_{j=1}^{\sf}B_j}$. Hence, $j
\in  \Per(h|_{S \setminus \cup_{j=1}^{\sf}B_j})$ if and only if there exists $x \in
P_j(f)$  such that $x \notin F_j$. By \lemc{musica}, this occurs if and only if $$
l(f^j)=P_j(f)>\Card(F_j)=\min\{jk_j, l(f^j)\}. $$ Since this is equivalent to
$l(f^j)>jk_j$, $$ \min\Per(h|_{S \setminus \cup_{j=1}^{\sf}B_j})=\min\{\sf\}
\cup\{j\,\,:\,\, l(f^j)>jk_j\}. $$ So  $\mmm(h)= \min \{ \sf, \beta(\overline{k},
l_\sf(f))  \}$, as desired.
\end{proof}

Let $f \in \fgp$. It is not always possible to find $\overline{k}$ satisfying the
hypotheses of \propc{fulgor}. However, if $f$  has a fixed point, we can take
$\overline{k}=(b,0,0,\dots)$. Also,  $k_1$ can take any value in  $\{0,1,\dots,b\}$
provided that the sum of all the $k_i$'s is $b$. This allows us, in several cases, to
choose $b$ and $\overline{k}$ in such a way that the map $h \in \hgb$ given by
\propc{fulgor} achieves a ``large'' $\mmm(h)$ and hence, a ``large'' lower bound for
$\mmm(\hgbp)$ is achieved. This result is stated in the following proposition which
includes one of the inequalities of Theorem~I.

\start{prop}{idioma} For $g \ge 2$, \tabc{singular} shows certain lower bounds for the
values of $\mmm(\hgbp)$.
\end{prop}

\begin{table}
\begin{center}
\begin{tabular}{|l|l|}
\hline $b$ & $\mmm(\hgbp) \ge$\\ \hline $b \le 2g+2$ & $b-2$\\ $b=2g+3$ & $2g$\\ $2g+4
\le b \le 3g+3$ & $2g+1$\\ $3g+3 \le b \le 4g+2$ & $b-g-2$\\ $4g+2 \le b \le 5g+3$ &
$3g$\\ $5g+2 \le b \le 6g+2$ & $b-2g-2$\\ $6g+2 \le b \le 6g+4$ & $4g$\\ $b=6g+5$ &
$4g+1$\\ $b \ge 6g+6$ & $4g+2$\\ \hline
\end{tabular}
\end{center}
\caption{Lower bounds for $\mmm(\hgbp)$.}\label{singular}
\end{table}

\begin{proof} By \coryc{triangle}, there exist maps $\vv_1, \vv_2, \vv_3 \in \fgp$ of types $$\typeop{4g}{1,1,2g}, \typeop{2g}{1,1,g} \mbox{ and } \typeop{4g+2}{1,2,2g+1}$$ respectively. Set $f_b=\vv_1$, $\vv_2$ or $\vv_3$  as indicated in \tabc{mapa}. Denote by  $\sigma_b$ the order of $f_b$ and  define $\overline{k}_b=(k_1,k_2,\dots,k_{\sigma_b})$ in the following way. If $i \notin \{1,2,g,2g,2g+1\}$ then $k_i=0$; otherwise, $k_i$ is defined as in \tabc{mapa}.
\begin{table}
\begin{center}
\begin{tabular}{|l|l|l|l|l|l|l|l|}
\hline $b$ &  $f_b$ & $k_{1}$ &$k_{2}$ &$k_{g}$ &$k_{2g}$ & $k_{2g+1}$&$\mmm(h_b)$\\
\hline $b \le 2g+2$  &$\vv_1$&   $b$ &   $0$ & $0$ &  $0$ & $0$&$b-2$ \\ $b=2g+3$
&$\vv_1$&  $b$ &    $0$ & $0$ &  $0$ & $0$&$2g$  \\ $2g+4 \le b \le 3g+3$  &$\vv_3$&
$b-2$ &  $2$ & $0$ &  $0$ & $0$&$2g+1$  \\ $3g+3 \le b \le 4g+2$  &$\vv_2$&  $b-g$ & $0$
& $g$ &  $0$ &  -  &$b-g-2$   \\ $4g+2 \le b \le 5g+2$  &$\vv_2$&   $b-g$ & $0$ & $g$ &
$0$ &  - &$3g$   \\ $5g+2 \le b \le 6g+2$   &$\vv_1$& $b-2g$ &$0$ & $0$ & $2g$ &
$0$&$b-2g-2$\\ $6g+2 \le b \le  6g+4$   &$\vv_1$& $b-2g$ &$0$ & $0$ & $2g$ &  $0$&$4g$\\
$b=6g+5$      &$\vv_3$&   $b-2g-3$ &  $2$ & $0$ & $0$ &  $2g+1$&$4g+1$\\ $b \ge 6g+6$
&$\vv_3$&   $b-2g-3$ &  $2$ & $0$ & $0$ &  $2g+1$ &$4g+2$\\ \hline
\end{tabular}
\end{center}
\caption{Values of $\overline{k}$  and $f$.}\label{mapa}
\end{table}
Denote by $\beta_b$ the value of the lower bound of $\mmm(\hgbp)$ of \tabc{singular}.
Applying  \propc{fulgor} to $\overline{k}_b$ and $f_b$, we obtain a map $h_b \in \hgbp$
such that $\mmm(h_b)=\beta_b$. Hence, $\mmm(\hgbp) \ge \mmm(h_b)=\beta_b$, as desired.
\end{proof}

\start{rem}{cos} As we shall see in Chapter~\ref{S10}, the lower bounds given in
\propc{idioma} are the best possible when $g=2$ and $b \notin \{1,2,3,4,6,8\}.$
\end{rem}

\start{rem}{ppp}Let $g \ge 2$. The results of \tabc{singular} are not best possible in
general. For example,  assume that there exist positive integers $p_1,p_2,p_3$ pairwise
coprime, such that each of them divides $2g+p_1+p_2+p_3-2$. Moreover, assume that there
exist non-negative integers $n_1, n_2, n_3$ such that $$ b=n_1p_1+n_2p_2+n_3p_3. $$ Let
$$n= \min (\{2g+p_1+p_2+p_3-2\}\cup \{p_i \,\,: n_i=0\} \cup \{(n_i-1)p_i\,\,:\,\ n_i \ge
2\}).$$ By Harvey's \theoc{Harvey} and \propc{fulgor}, there exists $f \in \hgbp$ such
that $\mmm(f)=n$.

Hence, $\mmm(\hgbp)$ is bounded below by the maximum of all $n$ obtained as above.
\end{rem}

By \lemc{jijior}, the fixed-point set of an iterate of an orientation-reversing map of a
closed surface can contain not only isolated fixed points, but also pointwise-fixed
simple closed curves.  Taking this into account, we can prove the following proposition,
in the same way as we proved \propc{fulgor}.

\start{prop}{cedro} Let $f \in \fgr$.  Assume that there exists $$\overline{k}=(k_1,2
k_2, \dots, \sf k_{\sf}) \in C_{b,\sf}$$ such that  $k_j=0$ if $l(f^j)=0$.  Then there
exists $h \in \hgbr$ such that $$\mmm(h)= \min \{ a_f, \beta(\overline{k}, l_{\sf}(\tf))
\},$$ where $a_f=\sf$  if $\dim(\Fix(f^{\sf/2}))=0$ and $a_f=\sf/2$, otherwise.
Furthermore, if for each $\ii{j}{\sf}$, if $jk_j-l(f^j)$ is a multiple of $\sf$ then
$\widetilde{h}=f$, so  $h$ is finite-order.
\end{prop}

The next result is a corollary of Propositions~\ref{fulgor} and\refc{cedro}

\start{cory}{cifra} Let $f \in \fgp$ (resp. $\fgr$) be a map of type
$$\typeop{n}{p_1,p_2,\dots,p_{_R}}.$$ Assume that $p_i \le p_j$ if $i \le j$ and that
$\sum_{i=1}^sp_i=b$ for some $s \le R$. Then there exists a map $h \in \fgbp$ (resp.
$\fgbr$) such that
\begin{eqnarray*}
\mmm(h)  &= &  \cases{p_{s+1}  & \mbox{if $s <R$,} \cr
                           n  & \mbox{if $s=R$.}}\cr
\end{eqnarray*}
\end{cory}

By \lemc{jijior}, orientation-reversing finite-order maps have no isolated  fixed points.
Therefore, if we want to apply \propc{cedro} to a given map $f \in \fgr$ and
$\overline{k}=(k_1,2k_2,\dots,bk_b)$,    $k_1$ must be equal to $0$. This restricts
considerably the  values of $\overline{k}$ satisfying the hypotheses of \propc{cedro}.
(Observe that the $i$-th entry  of $\overline{k}$ must be a multiple of $i$.) This is
solved by imposing the additional requirement on $f$ of having an invariant annulus which
plays the role of the fixed point in the orientation-preserving case, allowing us to
consider $k_1 \ne 0$.  This idea is formalized in the following.

\start{prop}{caoba} Let $f \in \fgr$. Suppose that the following hold.
\begin{numlist}
\item $\dim(\Fix(f^{\sf/2})=0$
\item There exists an $f$-invariant annulus $A \subset \Su{g}$ such that $f|_A$ is conjugate to the map $(z, t) \mapsto (z e^{\frac{2\pi i}{\sf}},1-t)$ on $\SI \times [0,1]$.
\item There exists $$(k_1,2k_2,\dots,\sf k_\sf) \in C_{b,\sf}$$
such that for each $j \in \{2,3,\dots, b\}$, $l(f^j)=0$ implies $k_j=0$.
\end{numlist}
Then there exists a map $h \in \hgbr$ such that $$\mmm(h)= \min \{ \sf,
\beta(\overline{k}, l_{\sigma_f}(\tf)) \}.$$
\end{prop}
\begin{proof} Define $D_1$ as in the proof of \propc{fulgor}. For each $j \le \sf$ we  define an $f$-invariant set $B_j \subset \Su{g}$ and a map $\Smap{\vv_j}{B_j}$ as follows. If $j=1$, we set $B_j=A$ and let $\vv_1$ be a map satisfying the following conditions.
\begin{numlist}
\item $\vv_1|_{\partial B_1}=f|_{\partial B_1}$ and $\vv_1$ is isotopic to $f|_{B_1}$ relative to $\partial B_1$.
\item $\vv |_{B_1}$ is conjugate to the map $(z, t) \mapsto (z e^{\frac{2\pi i}{k_1}},1-t)$ of $\SI \times [0,1]$.
\item $\mmm(\vv_1) = \min\{k_1, \sf\}$.
\end{numlist}
Now, for each   $i \in \{2,3,\dots,\sf\}$ we define $B_j$ and $\vv_j$ as in the proof of
\propc{fulgor}, and complete the proof with arguments analogous to those used there.
\end{proof}

Combining \lemc{8g} and \propc{caoba}, we obtain a bound for the minimum period of the
class of orientation reversing homeomorphisms of surfaces of odd genus.

\start{cory}{prisma} If $g \ge 2$ is odd and $b \ge 6g-6$, then there exists  $f \in
\hgr$ such that $h(f)=4g-4$. Consequently, $\mmm(\hgbr) \ge 4g-4$ if $b \ge 6g-6$.
\end{cory}
\begin{proof} Let $f \in \fg$ be as in \lemc{8g}. Define
$$\overline{k}=(k_1,k_2,\dots,k_{4g-4}) \in C_{b, 4g-4}$$ by

\begin{eqnarray*}
k_i & = & \cases{ b-2g+2 & \mbox{if  $i=1,$}\cr
                2g-2 & \mbox{if  $i=2g-2$},\cr
                               0 & \mbox{otherwise.}}\cr
\end{eqnarray*}

By \propc{caoba}, there exists $h \in \hgbr$ such that $\mmm(h)=\min\{4g-4,
b-2g+2\}=4g-4$.
\end{proof}

Another application of \propc{caoba} yields the following.

\start{cory}{felicidad} Let $g \ge 2$. There exists $f \in \hgbr$ such that $\mmm(f)
=\min\{b, 2g-2\}.$ Consequently, if $b \le 2g-2$, then $\mmm(\hgbr) \ge 2g-2$.
\end{cory}
\begin{proof}
Let $\overline{k}=(b,0,0,\dots) \in C_{b, 2g-2}$. If $g$ is even (resp. odd), applying
\propc{caoba} to $\overline{k}$ and the map of \lemc{rosa} (resp. \lemc{8g}), we obtain a
map $f \in \hgbr$, such that  $\mmm(f) =\min\{b, 2g-2\}$, as desired.
\end{proof}

The following lemma is the analogue of \lemc{idioma} for the orientation-reversing case
with the additional requirement that the genus $g$ be even.

\start{prop}{sabiduria} If $g \ge 2$ is even then \tabc{universo} shows certain lower
bounds for $\mmm(\hgbr)$.
\end{prop}
\begin{table}[t]
\begin{center}
\begin{tabular}{|l|l|}
\hline $b$ & $\mmm(\hgbr) \ge$\\ \hline $2g \le b \le 2g+2$ & $b-2$\\
 $2g+2 \le b \le 2g+4$  &  $2g$\\
 $2g+4 \le b \le 2g+6$  &  $b-4$\\
 $2g+6 \le b \le 4g+4$  &  $2g+2$\\
 $4g+4 \le b \le 6g+2$  &  $b-2g-2$\\
 $6g+2 \le b \le 6g+6$  &  $4g$\\
 $6g+6 \le b \le 6g+10$  &  $b-2g-6$\\
 $b \ge 6g+10$  &  $4g+4$\\
\hline
\end{tabular}
\end{center}
\caption{Lower bounds for $\mmm(\hgbr)$.}\label{universo}
\end{table}

\begin{proof} We prove this result by an argument analogous to that used in the proof of \propc{idioma}. In this case, we apply  \propc{fulgor} to $\vv_1, \vv_2 \in \fgr$, the maps defined in Lemmas~\ref{intimos dones} and \refc{8g+4} respectively.
Set $f_b=\vv_i$ as indicated in \tabc{criaturas}. Define
$\overline{k}_b=(k_1,k_2,\dots,k_{\sigma_{f_b}})$ in the following way. If $i \notin
\{1,2,4,2g,2g+2\}$ then $k_i=0$. Otherwise, $k_i$ is defined as in \tabc{criaturas}.
\begin{table}
\begin{center}
\begin{tabular}{|l|l|l|l|l|l|l|l|}
\hline $b$ &  $f_b$ & $k_{1}$ &$k_{2}$ &$k_{4}$ &$k_{2g}$ & $k_{2g+2}$&$\mmm(h_b)$\\
\hline $2g \le b \le 2g+2$  & $\vv_1$ & $b-2$   & $2$ & $0$ & $0$ & $0$&$b-2$ \\ $2g+2
\le b \le 2g+4$  & $\vv_1$ &$b-2$   & $2$ & $0$ & $0$ & $0$&$2g$ \\ $2g+4 \le b \le 2g+6$
& $\vv_2$ &$b-4$   & $0$ & $4$ & $0$ & $0$&$b-4$ \\ $2g+6 \le b \le 4g+4$  & $\vv_2$
&$b-4$   & $0$ & $4$ & $0$ & $0$ &$2g+2$\\ $4g+4 \le b \le 6g+2$  &  $\vv_1$ &$b-2g-2$
& $2$ & $0$ & $2g$ & $0$&$b-2g-2$ \\ $6g+2 \le b \le 6g+6$  & $\vv_1$ &$b-2g-2$   & $2$ &
$0$ & $2g$ & $0$&$4g$ \\ $6g+6 \le b \le 6g+10$  &$\vv_2$ &$b-2g-6$   & $0$ & $4$ & $0$ &
$2g+2$&$b-2g-6$ \\ $b \ge 6g+10$  &  $\vv_2$& $b-2g-6$   & $0$ & $4$ & $0$ & $2g+2$&$4g+4$ \\
\hline
\end{tabular}
\end{center}
\caption{Values of $\overline{k}$  and $f$.}\label{criaturas}
\end{table}

Denote by $\beta_b$ the value of the lower bound of $\mmm(\hgbr)$ of \tabc{universo}.
Applying  \propc{cedro} to $\overline{k}_b$ and $f_b$ we obtain a map $h_b \in \hgbr$
such that $\mmm(h)=\beta_b$. Hence, $\mmm(\hgbr) \ge \mmm(h_b)=\beta_b$ and the proof is
complete.
\end{proof}

\start{rem}{cosi} As we shall see in Chapter~\ref{S10} the lower bounds given in
\propc{sabiduria} are the best possible when $g=2$ and $b \notin \{1,2,3,4,6,10\}$.
\end{rem}

\section{Minimum periods of finite-order maps}\label{mp}
\markboth{Surfaces of low genus}{Minimum periods for finite-order maps}
\markboth{SURFACES OF LOW GENUS}{MINIMUM PERIODS OF FINITE-ORDER MAPS} In this section we
prove Theorems~B, C, D and E. Almost all the statements of Theorems~D and E follow from
Theorems~B and C. The basic idea for the proof Theorem~B (resp. Theorem~C)  is that for
each $g \ge 2$, $b \ge 3$,  a pair $(g,b)$ satisfies certain algebraic condition  if and
only if there exists a finite-order map $f \in \fgp$ (resp. $\fgr$) with exactly $b$
periodic points of $f$-period strictly less than $2g+b-2$. In such a case, \coryc{cifra}
gives the desired result.

\start{lem}{about b}Let $g \ge 2$ and $b \ge 5$ and suppose that for some  $f \in \fgbp$
(resp. $\fgbr$),  $m(f)=2g+b-2$. Then there exists $k \in \fgp$   (resp.$\fgr$)  of type
$$\typeop{2g+b-2}{p_1,p_2, \dots, p_{_R}}$$ where $R$ is a positive integer and
$\sum_{i=1}^Rp_i=b$.
\end{lem}

\begin{proof}
Let $f \in \fgbp$ (resp. $f \in \fgbr$) be such that $m(f)=2g+b-2$ and let $k=\tf$. Then
$k \in \fgp$ (resp. $k \in \fgr$). We  will prove that  the type of $k$  is
$\typeop{2g+b-2}{p_1,p_2, \dots, p_R}$ where $R \ge 1$ and $\sum_{i=1}^Rp_i=b$. We claim
that $\sigma_k=2g+b-2$. Clearly, $\sigma_k=\sf$. Since $f$ is finite-order and has a
periodic orbit of period $2g+b-2$, $\sigma_k=j(2g+b-2)$ for some positive integer $j$.
By \theoc{www}, $$ j=\frac{\sigma_k}{2g+b-2} \le \frac{4g+2}{2g+3}<2. $$ Thus, $j=1$ and
the claim is proved. Hence, the type of $k$ is $$\typeop{2g+b-2}{p_1,p_2,\dots,p_{_R}}$$
and it only remains to check that $R \ge 1$ and $\sum_{i=1}^Rp_i=b$.

Notice that  every point in $\Suf{g}{b}$ has $f$-period $2g+b-2$. Thus, the points of
$\Su{g}$ of $k$-period $p_1, p_2, \dots, p_R$ are collapsed boundary components. Let $m$
be the number of $k$-orbits which are collapsed boundary components of $k$-period
$2g+b-2$. Then $m$ is a non-negative integer and $b=\sum_{i=1}^Rp_i+m(2g+b-2)$. If $m \ge
1$ then $g \le 1$, which is impossible. So, $m=0$ and $b=\sum_{i=1}^Rp_i$. Since $b \ge
5$, $R \ge 1$, so the proof is complete.
\end{proof}

The main tool for the proof of the next result is \coryc{har}.

\begin{proclama}{Theorem B}\index{Theorem B} Let $g \ge 2$. Then $\mmm(\fgbp)=2g+b-2$ if and only if  $b \in \{2,3,4\}$ or there exist positive integers $p_1, p_2, p_3$ such that they are pairwise coprime, each of them divides $2g+b-2$, and  $p_1+p_2+p_3=b$.
\end{proclama}

\begin{proof}
We begin  with the ``only if'' direction. Assume that $b=1$. If there exists $f \in
\fgbpx{g}{1}$ such that $\mmm(f)=2g+b-2=2g-1$ then, since $g \le 2g-1$, $f,
f^2,\dots,f^g$ are fixed-point free. The single boundary component of $\Suf{g}{1}$ is
$f^i$-invariant for each $\ii{i}{g}$, so the iterates of the induced map $\tf^i$ for
$\ii{i}{g}$ have only one fixed point, the collapsed boundary component. By
\propc{index}, this fixed point has index one for $\tf^i$ for each $\ii{i}{g}$.  Then, it
follows from \theoc{suma}, that $L(\tf^i)=1$ for each $\ii{i}{g}$. By \lemc{nortumbria},
$L(\tf^{g+1})=-g<0$.  Hence, by \lemc{nw}, $$ -g=L(\tf^{g+1})=\chi(\Su{g})=2-2g $$ and
$\tf^{g+1}=\id$. Therefore, $g=2$ and $f^3=\id$.  Thus, $\sf$ divides $3$. Since $f,
f^2=f^g$ are fixed-point free, $\sf=3$. Then $\Per(\tf)=\{1,3\}$. On the other hand,
$\tf$ has only one fixed point. So, $f$ has type $\typeop{3}{1}$. Since this  contradicts
condition (5) of \coryc{har},   $b \ge 2$.

If $b \in \{2,3,4\}$ the conclusion holds. So we can assume that $b \ge 5$. By
hypothesis,  there exists $f \in \fgbp$ such that $\mmm(f)=2g+b-2$. By \lemc{about b},
there exists $h \in \fgp$ of  type $\typeop{2g+b-2}{p_1,p_2,\dots,p_{_R}}$, for some
positive integer $R$ and $p_1, p_2, \dots p_{_R} \ge 1$ such that $\sum_{i=1}^Rp_i=b$.
Set $$ T=\frac{2g-2+\sum_{i=1}^Rp_i}{2g+b-2}-R+2=3-R.$$

Since the condition (1) of  \coryc{har} holds,  $T$ must be even and non-negative. Then
$R \in \{1,3\}$. By  the condition (5) of \coryc{har}, $R \ne 1$. Therefore,  $R=3$ and
the desired conclusion holds by \coryc{triangle}.

Let us see the ``if'' direction. By Proposition~A, $\mmm(\fgbp) \le 2g+b-2$. To see
equality, it only remains to show that if the hypotheses hold then  there exists a map $f
\in \fgbp$ such that $\mmm(f)=2g+b-2$.

By \coryc{triangle}, there exist maps $f_1, f_2, f_3 \in \fgp$ of  type $$
\typeop{4g}{1,1,2g}, \typeop{2g+1}{1,1,1}, \typeop{2g+2}{1,1,2}.$$ respectively. For each
$b \in \{2,3,4\}$, applying  \coryc{cifra} to $f_{b-1}$ we can  see that there exists
$h_b \in \fgb$ such that $\mmm(h_b)=2g+b-2$.

By \coryc{triangle},  if $b \ge 5$, there exists a map of  type
$$\typeop{2g-2+p_1+p_2+p_3}{p_1,p_2,p_3}.$$ Since $p_1+p_2+p_3=b$  the result  holds by
\coryc{cifra}.
\end{proof}

\begin{rem} Combining a result of Gilman \cite{G} and \coryc{har} we obtain the following: suppose that $f \in \fgp$, then the isotopy class of $f$ is irreducible if and only if there exists positive integers $p_1, p_2$ and $p_3$ such that they are pairwise coprime and the type of $f$ is $\typeop{2g+p_1+p_2+p_3}{p_1, p_2, p_3}$. Hence, if $b \ge 3$, by Theorem~B, the maps $f \in \fgbp$ such that $\mmm(f)=2g+b-2$ are the ones for which the isotopy class of $\tf$ is irreducible.
\end{rem}

\begin{proclama}{Theorem D} \index{Theorem D} Let $g \ge 2$. Then $\mmm(\hgbp)=2g+b-2$ if one of the following conditions holds.
\begin{numlist}
\item There exist positive integers $p_1, p_2, p_3$ such that they are pairwise coprime, each of them divides $2g+b-2$, and  $p_1+p_2+p_3=b$.
\item $b-2$ divides $2g$.
\item $b-3$ divides $2g+1$.
\item $b \in \{1, 2, 3, 4, g+2, 2g+2,2g+4\}$.
\end{numlist}
\end{proclama}

\begin{proof} If (1) holds then the desired result is a consequence of Theorem~B.

If $b-2$ (resp. $b-3$) divides $2g$, we can apply Theorem~B to $p_1=p_2=1$ and $p_3=b-2$
(resp. $p_1=1$, $p_2=2$ and $p_3=b-3$). So, if (2) holds, we are done.

Now, assume that  (4) holds. If $b \in \{2,3,4\}$, the conclusion follows from Theorem~B.
If $b \in \{g+2, 2g+2\}$ (resp. $b \in \{2g+4\}$) then (2) (resp. (3)) holds, so the
desired result holds.

Finally, assume that $b=1$. By \lemc{ordertoro} and \lemc{triangle}, for each $g \ge 2$
there exists a map $f \in \fgpx{g-1}$ of type $\typeop{2g-1}{1,1,1}$. Let $\{x_1,
x_2,x_3\} \subset  \Su{g}$  be the $f$-fixed points. Set $F=\{x_1,x_2,x_3\}$ and let
$D=D_1 \cup D_2 \cup D_3$ be as in \lemc{diamante}. Clearly, $f(D_i)=D_i$ for each $i \in
\{1,2,3\}$. Now, we  glue to $\Su{g} \setminus D$ an annulus $A$ whose boundary
components are $\partial D_1$ and $\partial D_2$. Moreover, we extend $f|_{\Su{g}
\setminus D}$ to a map $$ \Smap{h}{(\Su{g} \setminus D) \cup A} $$ such that
$\mmm(h)=2g-1$. Since $(\Su{g} \setminus D) \cup A$ is a surface of genus $g$ with one
boundary component, the proof is complete.
\end{proof}

\begin{rem} The map constructed in the above proof for the case $b=1$ is not finite-order because, even if $f$ is conjugate to a rotation of order $2g-1$ around each fixed point, it can be proved that the angles of these rotations cannot sum to $0$, so $h|_A$ cannot be finite-order.
\end{rem}

\begin{proclama}{Theorem C} \index{Theorem C}Let $g \ge 2$. Then   $\mmm(\fgbr)=2g+b-2$ if and only if $b \in \{2,4\}$ or one of the following conditions holds.
\begin{numlist}
\item $g$ is even and there exist positive integers $p_1,  p_2$ such that $\gcd(p_1,p_2)=2$, each of them divides $2g+b-2$, and $p_1+p_2=b$.
\item $g$ is odd, $b$ is even, and $b$ divides $2g-2$.
\end{numlist}
\end{proclama}

\begin{proof}
We begin with the ``only if'' direction. Suppose $\mmm(\fgbr)=2g+b-2$.  Then there exists
$f \in \fgbr$ such that $\mmm(f)=2g+b-2$. By \lemc{suenio},  $b$ is even. If $b \in
\{2,4\}$ we are done. Hence, we can assume that $b \ge 5$. By \lemc{about b}, there
exists $h \in \fgr$ of  type $\typeop{2g+b-2}{p_1,p_2,\dots,p_{_R}}$, for some positive
integer $R$, $p_1, p_2, \dots p_{_R} \ge 1$. Set $$
T=\frac{2g-2+\sum_{i=1}^Rp_i}{2g+b-2}-R+2=3-R.$$ Since condition (4) of \coryc{anterior}
holds, $T \ge 1$. Hence, $R \in \{1,2\}$ and the desired conclusion follows from
\lemc{arpa}.

Let us see the ``if'' direction. By Proposition~A, $\mmm(\fgbr) \le 2g+b-2$. Suppose that
$g$ is odd and consider $b$ such that either $b \in \{2,4\}$ or $b$ divides $2g-2$. By
\lemc{arpa} there exists a map of type $\typeop{2g+b-2}{b}$. By  \lemc{cifra}, there
exists $f \in \fgbr$ such that $\mmm(f)=2g+b-2$, which proves that $\mmm(\fgbr) \le
2g+b-2$ in this case.

Now, assume that $g$ is even. By  \lemc{arpa}, there exist maps of type
$\typeop{4g}{2,2g}$ and $\typeop{4g+4}{4,2g+2}$. Applying \coryc{cifra} to these maps we
obtain the desired equality for $b \in \{2,4\}$. If $b \ge 5$, by \lemc{arpa}, there
exists a  map of type $\typeop{2g+b-2}{p_1,p_2}$. Since $p_1+p_2=b$ we can apply
\lemc{cifra} to this map in order to obtain $f \in \fgbr$ such that $\mmm(f)=2g+b-2$.
\end{proof}

Now we use Theorem~C to prove the following result.

\begin{proclama}{Theorem E} \index{Theorem E}Let $g \ge 2$.
\begin{numlist}
\item If $b$ is odd then $\mmm(\hgbr) \le b$, and equality holds if $b \le 2g-2$.
\item $\mmm(\hgbr)=2g+b-2$ if one of the following conditions holds.
\begin{romlist}
\item $b \in \{2,4\}$.
\item $g$ is odd, $b$ is even and $b$ divides $2g-2$.
\item $g$ is odd, $b \in \{g-1, 2g-2\}$.
\item $g$ is even, and  there exists positive integers $p_1,  p_2$ such that each of them divides $2g+b-2$, $\gcd(p_1,p_2)=2$, and $p_1+p_2=b$.
\item  $g$ is even, and $b-2$ divides $2g$.
\item  $g$ is even, and $b-4$ divides $2g+2$.
\item $g$ is even, and  $b \in \{g+2,2g+2,2g+6\}$.
\end{romlist}
\end{numlist}
\end{proclama}
\begin{proof}  We begin by proving (1). Let $f \in \hgbr$. If $b$ is odd, there is a boundary component $B$ of $\Suf{g}{b}$ such that $f^i(B)=B$ for some  odd $i \le b$. Since $f^i|_B$ is conjugate to an orientation-reversing homeomorphism of the circle, then, by \remc{grado}, $f^i|_B$ has a fixed point. Therefore, $\mmm(f) \le i \le b$. That equality holds for $b \le 2g-2$ follows from \coryc{felicidad}.

Let us prove (2). If (i), (ii), (iii) or (iv) hold, then the result is a consequence of
Theorem~C.  If (v)  (resp. (vi)) holds, we apply Theorem~C, to $p_1=2$ and $p_2=b-2$
(resp. $p_1=4$ and $p_2=b-4$) and obtain the desired conclusion.

Finally, we prove that   (vii) implies (2). If $b \in \{g+2, 2g+2\}$ (resp. $b=2g+6$),
then $b-2$ divides $2g$ ($b-4$ divides $2g+2$), and  (v) (resp. (vi)) holds, so the proof
is complete.
\end{proof}

\chapter{Homeomorphisms of surfaces of low genus}\label{S10}
The purpose of this chapter is to  study  the minimum periods of maps of $\Suf{g}{b}$ for
$g \in \{0,1,2\}$. As we will see, the cases $g=0$ and $g=1$ are not hard to solve,
whereas the case $g=2$ requires more effort. In the latter case, one of the inequalities
can be proved by means of the examples constructed in Chapter~\ref{S9}. The
main tools for studying the other inequality  will be  Newton's equations,
\lemc{symmetry}, and \propc{razon}. Indeed, by means of simple calculations we will show
that if a homeomorphism  $\Smap{h}{\Su{2}}$ is in a ``large'' set then $L(h^n)<0$ for
some $n \in \{1, 2, 3, 4\}$. So, \propc{razon} implies that the minimum period of  maps
$\Smap{f}{\Suf{2}{b}}$ such that $\tf$ is in this ``large'' set of homeomorphisms is less
than or equal to $4$.  For the rest of the maps  $\Smap{h}{\Su{2}}$ we will show that
the pair $(L(h), L(h^2))$ can take only finitely many values. For each of these possible
values $(c,d)$, \propc{razon}  will give  an upper bound for the minimum period of the
class of $f \in \hgbp$ (resp. $\hgbr$) inducing a map $\tf \in \hgp$ (resp. $\hgr$) such
that $(L(\tf),L(\tf^2))=(c,d)$.

\section{The orientation-preserving case}
\markboth{SURFACES OF LOW GENUS}{THE ORIENTATION-PRESERVING CASE} This section is devoted
to prove Theorem~F. We begin by introducing some notation.  For each $f \in \hg$ and each
positive integer $n$, let  $L_n(f)$ denote the $n$-uple of integers, $$(L(f),
L(f^2),\dots,L(f^n))$$ Also, denote  the infinite sequence $$(L(f), L(f^2), L(f^3),
\ldots)$$ by $L_\infty(f))$. For  $v=(v_1,\ldots,v_r) \in \Z^r$, we write $v^\infty$ for
the sequence $(v_1,\ldots,v_r,v_1,\ldots,v_r,\ldots) \in \Z^\N$.

\start{lem}{poema} Let $f \in \hgbp$. If $L(\tf) \ge 4$ then $\mmm(f) \le 2g$.
\end{lem}
\begin{proof} If $\mmm(f)>2g$, by \propc{index}, for each $\ii{i}{2g}$, the fixed points of $\tf^i$ are isolated and have index one. By \theoc{suma}, $\tf$ has at least $4$ fixed points, i.e., $\Card(\tf) \ge 4$. Since $\Fix(\tf) \subset \Fix(\tf^i)$, for every positive integer $i$, $\Card(\Fix(\tf^i)) \ge 4$, if $i \le 2g$.
Thus, by \theoc{suma},  $L(\tf^i) \ge 4$ for each $\ii{i}{2g}$. This contradicts
\lemc{amistad}, so the lemma is proved.
\end{proof}

\start{lem}{epitafio} If $f \in \hgbpx{2}{b}$ is such that $\mmm(f) \ge 5$ then either
$L_5(\tf)=(0,6,12,6,-20)$ or $$L_\infty(\tf) \in \{(0,4,6,4,0,-2)^\infty,
(1,3,1,3,6,1,3,1,3,-2)^\infty,$$ $$ (2,2,2,6,2,2,2,-2)^\infty,  (3,3,3,3,-2)^\infty,
(2,4,2,4,2,-2)^\infty\}.$$
\end{lem}
\begin{proof} Fix $f \in \hgbpx{2}{b}$ such that $\mmm(f) \ge 5$. For each positive integer $i$, let $\gg_i$ denote $l(\tf^i)$.
Combining Newton's equations  (page \pageref{newton}) with \lemc{symmetry}, and
\refeq{numerar}, we obtain the following system of equations,

\begin{eqnarray*}
p_1+s_1 &=& 0, \\ p_2+s_1 p_1+2s_2&=& 0,  \\ p_3+s_1 p_2 +s_2 p_1 +3 s_3 &=& 0, \\
p_4+s_1 p_3 +s_2 p_2 + s_3 p_1 +4 s_4 &=& 0,\\ s_1&=&s_3,\\ s_4&=&1,\\ L(\tf^i)&=&2-p_i,
\mbox{ for $i =1,2,3$ and $4,$}\\ \gg_1&=&L(\tf),\\ \gg_i&=&L(\tf^i)-L(\tf),\mbox{ for $i
=2$ and $3,$}\\ \gg_4&=&L(\tf^4)-L(\tf^2).
\end{eqnarray*}
Solving the system for $\gg_3$ and $\gg_4$ we get
\begin{eqnarray*}
\gg_3&=&  \frac{1}{2} (-12 +4 \gg_1 +3 \gg_1^2-\gg_1^3+6\gg_2-3\gg_1\gg_2),\\
\gg_4&=&\frac{1}{2} (-24 + 26 \gg_1 + 3 \gg_1^2  - 6 \gg_1^3  + \gg_1^4  + 10 \gg_2 - 10
\gg_1 \gg_2 +  2 \gg_1^2  \gg_2 - \gg_2^2 ).
\end{eqnarray*}
By  \lemc{poema}, $\gg_1 \in \{0,1,2,3\}$. For each of these values of $\gg_1$, the
values of $\gg_3$ and $\gg_4$ as a function of $\gg_2$ are given in  the second and third
columns of \tabc{vikings}, respectively. By \lemc{poema}, $\gg_i=l(\tf^i) \ge 0$ for each
$i \in \{1,2,3,4\}$. The last column of \tabc{vikings} gives us the values of $\gg_2$ for
which $\gg_3 \ge 0$, and $\gg_4 \ge 0$, for $\gg_1 \in \{0,1,2,3\}$.
\begin{table}
\begin{center}
\begin{tabular}{|l|l|l|l|}
\hline $\gg_1$ & $\gg_3$ & $\gg_4$ &$\gg_2$\\ \hline $0 $&$  3( \gg_2-2)$&$ (-24 + 10
\gg_2 - \gg_2^2)/2$&$ 4,6$\\ $1 $&$ 3( \gg_2-2)/2 $&$  \gg_2(2- \gg_2)/2$&$2$\\ $2 $&$
0 $&$            (8 -2 \gg_2 - \gg_2^2)/2$&$0,2$\\ $3 $&$ -3\gg_2/2 $&$
-\gg_2(\gg_2+2)/2$&$0$\\ \hline
\end{tabular}
\end{center}
\caption{Values of $\gg_3$ and $\gg_4$ for $\gg_1\in \{0,1,2,3\}$.}\label{vikings}
\end{table}
In other words, $L_2(\tf) \in \{(0,4),(0,6),(1,3),(2,2),(2,4),(3,3)\}$.

If $L_2(\tf)=(0,6)$, then by  Newton's equations  (page \pageref{newton}) and
\lemc{symmetry}, $L_5(\tf)=(0,6,12,6,-20)$.

By \coryc{har}, there exists $f_1 \in \fgpx{2}$ of type $$\typeop{6}{2,2,3,3},$$ and by
\coryc{triangle}, there exist maps $f_2,f_3,f_4,f_5 \in \fgpx{2}$ of types
$$\typeop{10}{1,2,5}, \typeop{8}{1,1,4}, \typeop{6}{1,1,2}, \typeop{5}{1,1,1}$$
 respectively. By  \lemc{nw}, for $i=1,2,3,4,5$, $L_\infty(f_i)$ equals
$$(0,4,6,4,0,-2)^\infty, (1,3,1,3,6,1,3,1,3,-2)^\infty, (2,2,2,6,2,2,2,-2)^\infty, $$
$$(2,4,2,4,2,-2)^\infty, (3,3,3,3,-2)^\infty,$$ respectively. Hence, the result follows
from  \remc{torre}.
\end{proof}

\start{lem}{25} If $f \in \hgbpx{2}{5}$ then $\mmm(f) \le 3$.
\end{lem}
\begin{proof}
Assume that there exists  $f \in \hgbpx{2}{5}$ such that $\mmm(f) >3$. Observe that $$
C_{5,5}=\{(0,0,0,0,5),(0,2,3,0,0),(1,4,0,0,0),$$
$$(1,0,0,4,0),(2,0,3,0,0),(3,2,0,0,0),(5,0,0,0,0)\}.$$ Hence, by \propc{razon}, $l(\tf)
\in \{0,1,2,3,5\}$. If $l(\tf)=5$ then, by \propc{razon},  $l_2(\tf)=(5,0)$, and by
Newton's equations  (page \pageref{newton}), and \lemc{symmetry},  $l(\tf^3)=-21 <0$, so,
by \remc{isla}, $\mmm(f) \le 3$, a contradiction.  So, $l(\tf) \in \{0,1,2,3\}$. By
\propc{razon}, $l_3(\tf) \in B$, where $$
B=\{(0,0,0),(0,2,3),(1,4,0),(1,0,0),(2,0,3),(3,2,0)\}. $$ On the other hand, in
\tabc{vikings},  for each $l(\tf) \in \{0,1,2,3\}$, the value of $l(\tf^3)$ is given in
terms of $l(\tf^2)$. Using this table, a simple calculation shows that for each $h \in
\hgpx{2}$, if $l(h) \in \{0,1,2,3\}$ then $l_3(h) \notin B$. Since this is a
contradiction, the lemma is proved.
\end{proof}

\start{lem}{ajedrez} The values of $\mmm(\hgbpx{2}{b})$ are given in \tabc{tortugi}.
\end{lem}

\begin{table}
\begin{center}
\begin{tabular}{|l|llllllllllllll|}
\hline b &1&2&3&4&5&6&7&8&9&10&11&12&13&14\\ \hline $\mmm(\hgbpx{2}{b})$
&3&4&5&6&3&8&4&10&5&6&6&6&7&8\\ \hline
\end{tabular}
\end{center}
\begin{center}
\begin{tabular}{|l|llll|}
\hline b&15 &16&17&$b \ge 18$\\ \hline $\mmm(\hgbpx{2}{b})$&8 &8&9&10\\ \hline
\end{tabular}
\end{center}
\caption{Values of $\mmm(\hgbpx{2}{b})$.}\label{tortugi}
\end{table}

\begin{table}
\begin{center}
\begin{tabular}{|l|lllllllllllllll|}
\hline $l_2(\tf) / b$ & 5 & 6 & 7 & 8 & 9 & 10 & 11 & 12 & 13 & 14 & 15 & 16 & 17 & 18 &
$\cdots$\\ \hline (0,4) & 2 & 2 & 3 & 3 & 3 & 6 & 3 & 3 & 3 & 4 & 5 & 6 & $\cdots$ & & \\
(1,3) & 2 & 3 & 4 & 10 & 5 & 5 & 5 & 5 & 5 & 6 & 7 & 8 & 9 & 10 &$\cdots$ \\ (2,2) & 3 &
8 & 4 & 4 & 4 & 4 & 5 & 6&7&8&$\cdots$& & & & \\ (2,4) & 2 & 2 & 3 & 4 & 5 & 6 & $\cdots$
& & && & & & & \\ (3,3) & 2 & 3 & 4 & 5 & $\cdots$&& & & & && & & & \\ (0,6) & 2 & 3 & 2
& 2 & 3 & 3 & 3 & 3&3 & 3&3 & 3&3&5 &$\cdots$\\ \hline
\end{tabular}
\caption{Values of $\gg(\tf,b)$ for $b \ge 5$ and $\tf \in B$.}\label{f2}
\end{center}
\end{table}

\begin{proof} For each $b$, denote by $m_b$ the value claimed for $\mmm(\hgbpx{2}{b})$ in \tabc{tortugi}. If $b \in \{1,2,3,4,6,8\}$, by Theorem~C, $\mmm(\hgbpx{2}{b})=m_b$.

By \propc{idioma}, for each $b \notin  \{1,2,3,4,6,8\}$ there exists a map $f \in
\hgbpx{2}{b}$ such that $\mmm(\hgbpx{2}{b}) \ge m_b$. In particular, by  \lemc{25},
$\mmm(\hgbpx{2}{5})=3$.

In \tabc{f2}, we list the values of $\gamma(\tf,b)$ for each $f$ such that $l_2(\tf) \in
B$, where $$B=\{(0,4), (1,3),(2,2),(2,4),(3,3),(0,6)\}.$$ By \propc{razon}, if $f \in
\hgbp$ is such that $l_2(\tf) \in B$ then $\mmm(f)$ is less than or equal to the
corresponding entry of \tabc{f2}.

Consider $f \in  \hgbx{2}{7}$. We claim that $\mmm(f) \le 4$.  Indeed, if $\mmm(f)>4$, by
\lemc{epitafio}, \propc{razon} and \tabc{f2}, $\mmm(f) \le 4$, a contradiction. Hence the
claim is proved. Therefore, $4 \ge \mmm(\hgbx{2}{7}) \ge m_7=4$.

If $b \notin \{1,2,3,4,5,7\}$ then  $\mmm(\hgbpx{2}{b}) \ge m_b >4$. So, by
\lemc{epitafio}, \propc{razon} and \tabc{f2}, $m_b \le \mmm(\hgbpx{2}{b}) \le m_b$, as
desired.
\end{proof}

\start{rem}{costumbre} Observe that $\mmm(\hgbpx{2}{b})=2.2+b-2=b+2$ if and only if $b=1$
or $\mmm(\fgbpx{2}{b})=b+2$.
\end{rem}

\begin{proclama}{Theorem F}\index{Theorem F}
 \begin{numlist}
\item
\begin{eqnarray*}
\mmm(\hgbpx{0}{b}) & = &  \left\{\begin{array}{ll}
                                           1 & \mbox {if  $b=1$,}\\
                                           \infty & \mbox{if $b=2$,}\\
                                          b-2  & \mbox{if  $b \ge 3$.}
                           \end{array}
                           \right.\\
\end{eqnarray*}
\item
\begin{eqnarray*}
\mmm(\hgbpx{1}{b})& = &   \left\{\begin{array}{ll}
                        2 & \mbox{if $b=1$,}\\
                        b & \mbox{if  $b \ge 2$.}
                     \end{array}
                              \right.
\end{eqnarray*}
\item \tabc{tortugi} shows the values of $\mmm(\hgbpx{2}{b})$.
\end{numlist}
\end{proclama}

\begin{proof}We will split the proof into various cases.
\begin{case}{1}{$g=0$, $b=1$.}
\end{case}
In this case, the result can be deduced from Brouwer's Fixed-Point Theorem and also, from
Fuller's \theoc{fuller}.

\begin{case}{2}{$g=0$, $b=2$.}
\end{case}

View the annulus $\Suf{0}{2}$ as $\SI \times I$. Consider the homeomorphism $(z, \rho)
\mapsto (R_{\alpha}(z), \rho)$, where $\alpha \in \R$ is irrational. It is clear that it
preserves orientation and that it has no periodic points.

\begin{case}{3}{$g=0$, $b=3$.}
\end{case}
Observe that if $f \in \hgbpx{0}{3}$ then $\tf$ is an orientation-preserving
homeomorphism of the sphere $\Su{0}$.  Since the first homology group of the sphere
$\Su{0}$ is trivial, by  \refeq{anne},  $L(\tf^i)=2$ for each positive integer $i$.
Therefore, $l_3(\tf)=(2,0,0)$. Since $C_{3,3}=\{(3,0,0),(0,0,3),(1,2,0)\}$, the result
follows from \propc{razon}.

\begin{case}{4}{$g=0$, $b \ge 4$.}
\end{case}
By Proposition~A, $h_p(\Suf{0}{b}) \leq b-2$. The following example gives the reverse
inequality. Consider a sphere with $b-2$ holes symmetrically distributed on the equator
and two more at the poles (see \figc{grupoop}) and take $f$ to be rotation through an
angle of $2\pi /(b-2)$ with respect to the axis $R$. Clearly, $f$ has minimum period
$b-2$.

\placedrawing{cg.lp}{Examples for the proofs of Theorem F}{grupoop}

\begin{case}{5}{$g=1$, $b=1$.}
\end{case}
By Fuller's \theoc{fuller}, $h_p(\Suf{1}{1}) \le 2$. To complete the proof of this case
we shall exhibit an example of a map on $\Suf{1}{1}$ without fixed points.

By  \lemc{ordertoro}, there exists a map $f \in \fgpx{1}$ of type $\typeop{6}{1,2,3}$.
Let $x$ be the fixed point of $f$. Let $D$ be a set as in \lemc{diamante} for $F=\{x\}$.
Then $\Smap{f|_{\Su{1}\setminus D}}{\Su{1} \setminus D}$ is a homeomorphism  without
fixed points, as desired.

\begin{case}{6}{$g=1$, $b \ge 2$.}
\end{case}
By Proposition~A, $\mmm(\hgbpx{1}{b}) \leq b$. To see that  equality holds, consider a
torus with $b$ holes distributed as in \figc{grupoop}. Rotation through an angle of $2\pi
/b$ with respect to the axis $R$ has minimum period $b$, so the proof of this case is
complete.

\begin{case}{7}{$g=2$.}
\end{case}
See \lemc{ajedrez}.
\end{proof}

\section{The orientation-reversing case}
\markboth{SURFACES OF LOW GENUS}{THE ORIENTATION-REVERSING CASE} As in the
orientation-preserving case, most of this section will be devoted to studying
$\mmm(\hgbrx{2}{b})$. One difference between the two cases is given by the following
result.

\start{lem}{monedas} Let $i$ be a positive odd integer and $f \in \hgbr$. If  $L(\tf^i)
\ne 0$ then $\mmm(f) \le i$.
\end{lem}
\begin{proof} Assume that $\mmm(f)>i$. Since  $L(\tf^i) \ne 0$, by \theoc{lef}, $\tf$ has a fixed point. Since $\mmm(f)>i \ge 1$, this fixed point must be a collapsed boundary component $B$. Clearly, $B$ is $f^i$-invariant. Since $f^i$ is orientation-reversing, $f^i|_B$ is conjugate to an orientation-reversing homeomorphism of the circle. Thus, by \remc{grado}, $f^i|_B$ has a fixed point. Then, $m(f) \le i$, a contradiction.
\end{proof}

\start{lem}{inagotable} Let $f \in \hgbr$. If $L(\tf^2)>4$   then $\mmm(f) \le 2g$.
\end{lem}
\begin{proof} Assume that $\mmm(f) >2g$ and $L(\tf^2)>4$. By \lemc{monedas},  $L(\tf^i)= 0$ for each  $i \in \{1,3,5,\dots,2g-1\}$.

By \propc{index},  for each $\ii{i}{2g}$, the fixed points of $\tf^i$ are isolated and
have index one. Since $\Fix(\tf^2) \subset \Fix(\tf^{2i})$, for every positive integer
$i$, $\Card(\Fix(\tf^{2i})>4$ for each $\ii{i}{g}$. By \theoc{suma}, $4<L(\tf^{2i})$ for
each $\ii{i}{g}$. This contradicts \lemc{verlaine}, so the proof  is complete.
\end{proof}

\begin{table}[t]
\begin{center}
\begin{tabular}{|l|llllllllllllll|}
\hline b &1&2&3&4&5&6&7&8&9&10&11&12&13&14\\ \hline $\mmm(\hgbrx{2}{b})$
&1&4&3&6&4&8&4&4&5&12&6&6&7&8\\ \hline
\end{tabular}
\end{center}
\begin{center}
\begin{tabular}{|l|llllllll|}
\hline b&15 &16&17&18&19&20&21&$b \ge 22$\\ \hline $\mmm(\hgbrx{2}{b})$&8 &8&8&  8&  9&
10&11&12\\ \hline
\end{tabular}
\end{center}
\caption{Values of  $\mmm(\hgbrx{2}{b})$.}\label{tortugo}
\end{table}

\start{lem}{noche} The values of $\mmm(\hgbrx{2}{b})$ are as given in \tabc{tortugo}.
\end{lem}

\begin{table}[t]
\begin{center}
\begin{tabular}{|l|lll|}
\hline $b \setminus l_2(\tf)$ & (0,0) &  (0,2) & (0,4) \\ \hline 5  &   4   &    3   &
2   \\ 6  &   4   &    8   &   2   \\ 7  &   4   &    4   &   3   \\ 8  &   4   &    4
&   4   \\ 9  &   5   &    4   &   5   \\ 10 &   12   &    4   &   6   \\ 11 &   6  &
5   &   6   \\ 12 &   6   &    6   &$\vdots$\\ 13 &   6   &    7   &       \\ 14 &   6
&    8   &       \\ 15 &   6   &    8   &       \\ 16 &   6   &$\vdots$&       \\ 17 &
7   &        &       \\ 18 &   8   &        &       \\ 19 &   9   &        &       \\ 20
&  10   &        &       \\ 21 &  11   &        &       \\ 22 &  12   &        &       \\
b  &    $\vdots$   &        & \\ \hline
\end{tabular}
\caption{Values of $\gg(\tf,b)$ for $b \ge 5$ and $\tf \in B$.}\label{f2r}
\end{center}
\end{table}

\begin{proof} For each $b$, denote by $m_b$ the value claimed for $\mmm(\hgbrx{2}{b})$ in \tabc{tortugo}.

If $b \in \{1,2,4\}$ then conditions (1) and (2i) of Theorem~E hold, so
$\mmm(\hgbrx{2}{b})=m_b$.

By Theorem~E, $\mmm(\hgbrx{2}{3}) \le 3$. By \lemc{arpa}(2), there exists $f \in
\fgrx{2}$ of type $\typeop{12}{4,6}$. Let $\overline{k} \in C_{3,12}$ be such that
$k_1=3$ and $k_i=0$ for each $i \in \{2,3,\dots,12\}$. By \propc{fulgor}, there exists $h
\in \hgbrx{2}{3}$ such that $\mmm(h)=3$. Thus, $\mmm(\hgbrx{2}{3})=3=m_3$.

If $b \in \{6,10\}$, condition (vii) of Theorem~E holds, so $\mmm(\hgbrx{2}{b})=m_b$.

Hence, we can assume $b \ge 5$ and $b \notin \{6,10\}$. Here,  by \propc{sabiduria},
$\mmm(\hgbrx{2}{b}) \ge m_b$. Suppose there exists $f \in \hgbrx{2}{b}$ such that
$\mmm(f)>m_b$.  Since $m_b \ge 4$,  by Lemmas~\ref{monedas} and \ref{inagotable},
$\l_2(\tf) \in B$, where $$ B=\{(0,0),(0,2),(0,4)\}. $$ In \tabc{f2r}, we list the values
of $\gamma(\tf,b)$ for each such $f$. By \ \propc{razon},  $\mmm(\hgbrx{2}{b}) \le m_b$,
which contradicts our assumption and completes the proof.
\end{proof}

\start{rem}{costumbrita} As in  the orientation-preserving case,
$\mmm(\hgbrx{2}{b})=2.2+b-2$ if and only if $\mmm(\fgbrx{2}{b})=4+b-2$.
\end{rem}

\begin{proclama}{Theorem G} \index{Theorem G}  \begin{numlist}
\item
\begin{eqnarray*}
\mmm(\hgbrx{0}{b}) & = &  \left\{\begin{array}{ll}
                                           1 & \mbox {if  $b=1$,}\\
                                           \infty & \mbox{if $b=2$,}\\
                                           2 & \mbox{if $b=3$,}\\
                                          b-2  & \mbox{if  $b \ge 4$.}
                           \end{array}
                           \right.\\
\end{eqnarray*}
\item $\mmm(\hgbrx{1}{b})  =  b-2$.
\item \tabc{tortugo} shows the values of $\mmm(\hgbrx{2}{b})$.
\end{numlist}
\end{proclama}
\begin{proof} We split the proof into various cases.

\begin{case}{1}{$g=0$, $b=1$.}
\end{case} See Fuller's \theoc{fuller}.

\begin{case}{2}{$g=0$, $b=2$.}
\end{case}
View the annulus $\Suf{0}{2}$ as $\SI \times [0,1]$ and consider the homeomorphism
 $(z, \rho) \mapsto (R_{\alpha}(z), -\rho)$, where $\alpha \in \R$ is irrational. It is clear that this homeomorphism reverses orientation and has no periodic points.

\begin{case}{3}{$g=0$, $b=3$.}
\end{case}
By Fuller's  \theoc{fuller}, $h_r(\Suf{0}{3}) \le 2$. We  prove equality by means of an
example.  Consider a sphere $\Suf{0}{3}$ with three holes symmetrically distributed on
the equator, see \figc{grupor}. Let $\Smap{r}{\Suf{0}{3}}$ be rotation through an angle
of $2 \pi /3$ with respect to the axis $R$, and  let $\Smap{s}{\Suf{0}{3}}$ be reflection
in the plane containing the equator. Define $f=s \circ r$. Clearly,  $f$ has no fixed
points.

\begin{case}{4}{$g=0$, $b \ge 4$.}\end{case}
By Proposition~A, $h_r(\hgbrx{0}{b}) \le b-2$. The following example shows that equality
holds. Let $\Smap{f}{\Suf{0}{b}}$ be the map constructed in the proof of  Case $g=0$, $b
\ge 4$ of the proof of Theorem~B, and let $s$ be reflection in  the plane containing the
equator. Then $m(s \circ f)=b-2$.

\begin{case}{5}{$g=1$.}
\end{case}
By \lemc{felicidad} and Proposition~A,  $h_p(\Suf{1}{b}) \le b$. The following example
shows that equality holds. Consider a torus with $b$ holes distributed as  in
\figc{grupor}. Let $\Smap{r}{\Suf{1}{b}}$ be rotation through an angle of $2\pi /b$ with
respect to $R$ and let $\Smap{s}{\Suf{1}{b}}$ be a reflection in the plane $P$. Then
$h(s \circ r)=b$.

\placedrawing{grupor.lp}{Examples for the proofs of Theorem~G}{grupor}

\begin{case}{6}{$g=2$.}
\end{case}
See \lemc{noche}.
\end{proof}

\chapter{Proof of Theorems~H and I}\label{S11}

In this chapter we complete the proofs of Theorems H and I. Section~\ref{H} is dedicated
to the former, and Section~\ref{I}, to the latter.

\section{Proof of Theorem~H}\label{H}

This section is practically entirely devoted  to proving the following theorem, which
will be used to deduce the difficult inequality of  Theorem~H.

\start{theo}{resumenop} Let $g \geq 2$ and let $\Smap{f}{\Su{g}}$ be an orientation
preserving homeomorphism. Then there exists a positive integer $m$ such that $m \le 4g+2$
and  $f^m$ has a non-empty fixed-point class of  non-positive index.
\end{theo}

The strategy of the proof of \theoc{resumenop} consists in studying  fixed-point classes
of non-positive index of iterates of maps in standard form. The finite-order case is
trivial, so we shall concentrate on the pseudo-Anosov and reducible cases, in
Propositions~\ref{pA}, \refc{reducibleop}, and \refc{reduop}.

\start{prop}{pA} Let $g \geq 2$ and let \Smap{f}{\Su{g}} be an orientation-preserving
pseudo-Anosov homeomorphism in standard form. Then there exists a positive integer $m$
such that $m \le 4g$ and $f^m$ has a fixed-point class of negative index.
 \end{prop}

\begin{proof}
Let $$\{x_1,x_2, \ldots, x_k\} \subset \Su{g}$$ be such that the set  of singularities of
the foliation on $\Su{g}$ is the disjoint union of the $f$-orbits of the $x_i$'s. For
each $\ii{i}{k}$, let $n_i$ and $p_i$ be the period of $x_i$, and the number of prongs
emanating from $x_i$,  respectively. The Euler-Poincar\'{e} Formula  \refeq{ep} may be
written in the form,
\begin{equation}
\sum_{i=1}^kn_i(p_i-2)=4(g-1). \label{epp1}
\end{equation}
By  \lemc{libros} and \lemc{espejos}(1) we can assume that  $k \in \{1,2\}$. If $k=2$,
without loss of generality, we can assume that $n_1 \le n_2$ and, by \lemc{espejos}(2),
we can assume that $p_2=3$.

We now prove that for some positive integer $m$ such that $m \le 4g$,  one of the parts
of \lemc{laberinto} applies.

Notice that if there exists a regular point of period $n$ and $n \le 2g$ then
\lemc{laberinto}(6) applies for $m=2n$. Hence, we may assume that there are no regular
points of period  less than or equal to $2g$. In particular,
\begin{equation}
\Per(f) \cap \{1,2,\dots,2g\} \subset \left\{\begin{array}{ll}
                \{n_1\} & \mbox{if $k=1$,}\\
                \{n_1, n_2\} & \mbox{if $k=2$.}\\
                \end{array}
                \right.
                \label{sky}
\end{equation}

For every positive integer $h$, $x_1$ is a fixed point of $f^{n_1h}$. If
$\ind{f^{n_1h}}{x_1} \ne 1$ then \lemc{laberinto}(1) applies. Thus, we can assume that
$\ind{f^{n_1h}}{x_1}=1$ for every positive integer $h$ such that $n_1h \le g$. Now, we
split the proof into four cases.

\begin{case}{1}{$n_1 \in \{1,2\}$.} \end{case}
By \refeq{epp1}, $n_1p_1 \le 4(g-1)+2n_1 \le 4g$. So, \lemc{laberinto}(4) applies for
$m=n_1p_1$.

\begin{case}{2}{$k=2$ and $n_1, n_2 \le g+1$.} \end{case}
By \refeq{epp1}, $$ n_1 p_1+n_2 p_2 \leq 4(g-1)+2(n_1+n_2) \leq 4(g-1)+4g+4=8g. $$ Hence,
$n_ip_i \le 4g$ for some $i \in \{1,2\}$. So,  \lemc{laberinto}(4) applies for $m=n_ip_i$
for some $i \in \{1,2\}$, as desired.

\begin{case}{3}{$n_1 \ge g+1$.}\end{case}
By \refeq{sky}, since $g+1 \le 2g$, $\mmm(f) \ge g+1$. By \theoc{lef},
$$L(f)=L(f^2)=\ldots=L(f^{g})=0.$$ By Lemmas~\ref{polinor} and\refc{symmetry},  the
characteristic polynomial of $\fs$ is $(x-1)^2(x^{2g-2}+1)$. Thus, the $(4g-4)$-th power
of each eigenvalue is $1$. By \lemc{potency}, \lemc{laberinto}(2) applies with $m=4g-4$.

\begin{case}{4}{$3 \le n_1 \le g$ and,  if $k=2$, $n_2 >g+1$.}  \end{case}
Since $3 \le n_1 \le g \le 2g$, by \refeq{sky}, $\Per(f) \cap \{1,2,\dots,g\}=\{n_1\}$.
Moreover, the orbit of $x_1$ is the only periodic orbit whose $f$-period is smaller than
$g$. Then the hypotheses of  \propc{cas-part-op} hold with $n=n_1$. If for some
$\ii{m}{3g-3}$, $L(f^m)<0$ then \lemc{laberinto}(2) applies, as desired. Otherwise, there
exists a periodic point $y$ of period $l$ where $g+2 \le l \le 4(g-1)/3$.  Since
$4(g-1)/3 \le 2g$, $y$ is not a regular point. Therefore, $k=2$, $y$ is in the orbit of
$x_2$ and $l=n_2$. Recall that  $p_2=3$. Then $p_2n_2 =3l \le 4(g-1)$. Thus,
\lemc{laberinto}(4)  applies with $m=n_2p_2$, and the proof is complete.
\end{proof}

\start{prop}{reducibleop} Let \Smap{f}{\Su{g}} be an orientation-preserving reducible
homeomorphism in standard form which has a pseudo-Anosov component. Then there exists a
positive integer $m$ such that $m \leq 4g-4$ and $f^m$ has a fixed-point class  of
negative index.
\end{prop}
\begin{proof} Let $C$ be a pseudo-Anosov $f$-component. Denote by $n_1$, $g_1$, and $b_1$ the period, genus, and number of boundary components of $C$, respectively. Let  $B \subset C$ be  a $p$-pronged boundary component of $C$. Let $r$ be the least positive integer such that $f^{rn_1}(B)=B$.

By \remc{pprong}, $\sum (2-p_s) =-pr,$ where the sum is taken over all the prongs
emanating from  singularities $s$ lying in $\cup_{i=1}^r f^{n_1i}(B)$. Thus, since $p_s
\ge 3$,

\begin{equation}
pr \le \sum (p_s-2), \label{cheese}
\end{equation}
where the sum is taken over all singularities $s$ lying in $C$. By
the Euler-Poincar\'{e} Formula  \refeq{ep} applied to $\Smap{f^{n_1}|_{C}}{C}$, $$ \sum
(p_s-2) =-2 \chi(C)=2(2g_1+b_1-2), $$ where the sum is taken over all the singularities
$s$ of the foliation of $C$. Thus, by  \lemc{corep}(1)  and \refeq{cheese}, $$ prn_1 \le
2n_1(2g_1+b_1-2)  \le 2(2g-2)=4g-4. $$ Now, \lemc{tiger}(1)  applies and we are done.
\end{proof}

\start{prop}{reduop} Let $g \ge 2$ and let $\Smap{f}{\Su{g}}$ be an orientation
preserving reducible homeomorphism which satisfies the following conditions.
\begin{numlist}
\item $L(f) \ne 2$ or $L(f^2) \ne 0$.
\item Each  of its components  is finite-order.
\item $f$ is not finite-order.
\end{numlist}
Then there exist an $f$-component $C$ and a positive integer $n$ such that $n \le 4g$,
$f^n|_C=\id|_C$.
\end{prop}
\begin{proof} With \notc{dan}, for each $\ii{i}{k}$, $f^{n_i\si}|_{C_i}=\id_{C_i}$. Hence, it suffices to show that $n_i\si \le 4g$ for some $\ii{i}{k}$. We may assume that $g_1 \le g_2 \dots \le g_k$. Let us split the proof into various cases.

\begin{case}{1} $k \ge 2$, $g_1=g_2=0$.
\end{case}
By \lemc{corep}(1),  $n_1(b_1-2)+n_2(b_2-2) \le 2g-2$. Changing subindices if necessary,
by \lemc{corep}(4), we can assume that $n_1 \le n_1(b_1-2) \le g-1$. By \lemc{items}(1),
it suffices to prove the result for  $\sigma_1=3$.  Now, $\sigma_1 n_1=3n_1 \le 3g-3,$ as
required.

\begin{case}{2} $k \ge 2$, $g_1 \ge 1$, $g_2 \ge 1$.
\end{case}
By \lemc{euler}, $n_1g_1+n_2g_2 \le g$. Changing subindices if necessary, we can assume
that,  $n_1 g_1 \le g/2$.  If $g_1=1$, by \lemc{items}(2) we can assume that $\sigma_1
\le 6$. If $g_1 \ge 2$, by \coryc{wwwc}, $\sigma_1 \le 4g_1+2$. Since $n_1 \le g/2$, $$
n_1\sigma_1 \le n_1(4g_1+2) \le 2g+g \le 4g, $$ and the desired conclusion holds.

\begin{case}{3} $k \ge 2$, $g_1=0$, $g_2 \ge 1$.
\end{case}

By \lemc{items}(1), it suffices to prove the result for $b_1=\sigma_1=3$. Now, if $3n_1
\le 4g$  the result holds. Otherwise, $n_1>\frac{4g}{3}$. By \lemc{corep}(1),
$n_1+n_2(2g_2+b_2-2) \le 2g-2$. Then $$ n_2(2g_2-1) \le n_2(2g_2+b_2-2) <
2g-2-\frac{4g}{3} \le \frac{2g}{3} . $$ If $g_2=1$, by \lemc{items}(2), we can assume
that $\sigma_1 \le 6$. Hence, by
 \coryc{wwwc}, $\sigma_2 \le 4g_2+2$, so
$$ n_2 \le n_2\sigma_2 \le 2n_2(2g_2+1)=2n_2(2g_2-1)+4n_2 \le
\frac{4g}{3}+\frac{8g}{3}=4g, $$ which completes the proof of this case.

\begin{case}{4} $k=1$, $g_1=0$.
\end{case}
By \lemc{items}(1), we can assume that  $b_1=\sigma_1=3$. By  \coryc{s0b},  the three
boundary components of $C_1$ form a cycle under the action of $f^{n_1}$. By \lemc{ciclo},
$n_1=2$, and $g=2g_1+b_1-1=2$. So $\sigma_1 n_1=6 \le 8=4g$, and the desired conclusion
holds.

\begin{case}{5} $k=1$, $g_1 \ge 1$, $b_1 \ge 2$.
\end{case}

By \lemc{items}(2) and (3), we can assume that  $b_1=2$. By \coryc{wwwc} and
\lemc{items0}(1),  $\sigma_1\in  \{4g_1+1, 4g_1+2\}$. If $n_1=1$,   by \lemc{corep}(1),
$2g_1=n_1(2g_1+b_1-2)=2g-2$. Here, $$ n_1\sigma_1 \le 4g_1+2=4g-2, $$ as desired. If $n_1
\ge 2$, by \lemc{one fixed}(1), both boundary components are interchanged under the
action of $f^{n-1}$. Then \lemc{ciclo} gives $n_1=2$ and $g=2g_1+b_1-1=2g_1+1$. Hence, $$
n_1\sigma_1 \le 2(4g_1+2) =4g. $$ and we are done.

\begin{case}{6} $k=1$, $g_1 \ge 1$, $b_1=1$.
\end{case}

If  $n_1=1$, then by \lemc{corep}(1), $2g_1-1=2g-2$, which is impossible. Then we can
assume that  $n_1 \ge 2$. By \lemc{items0}(1) and \coryc{wwwc}, we can assume that
$\sigma_1\in  \{4g_1+1, 4g_1+2\}$. By \lemc{ciclo}, $n_1=2$.

Let $A$ be the closed annulus connecting $C_1$ and $f(C_1)$. Since $f|_{A}$ is an
orientation-preserving map which interchanges the boundary components of $A$, the
description of standard form (page \pageref{S4}) shows that there exists $a \in \Q$
such that  $f|_{A}$ is conjugate to the map  $\Smap{\psi}{\SI \times [0,1]}$ defined by
$$ (z,t) \mapsto  (\overline{z}e^{a(1-2t)\pi i},1-t). $$ Clearly, $\SI \times
\{\frac{1}{2}\}$ is $\psi$-invariant . Moreover, since $\psi|_{\SI \times
\{\frac{1}{2}\}}$ acts as the map $z \mapsto \overline{z}$, $\tr(\psi_{*1})=-1$.
Therefore $\tr((f|_{A})_{*1})=-1$, so, by \refeq{anne},  $L(f|_A)=2$. Observe that
$\Su{g}=C_1 \cup f(C_1) \cup A$ and, since $C_1 \cap f(C_1)=\emptyset$ and
$f^2(C_1)=C_1$, $f|_{C_1 \cup f(C_1)}$ does not have fixed points. Hence, there exist two
open subsets of $\Su{g}$, $U$ and $V$ such that $\Fix(f) \subset U \subset \Int(A)$, $U
\cup V=\Su{g}$ and $V \cap U \cap \Fix(f)=\emptyset$. Applying  \theoc{suma} twice we
obtain $$ L(f)=\inn(f|_{U})+\inn(f|_V)=L(f|_A)=2. $$ Now, observe that $f^2|_A$ is
conjugate to  the map $\Smap{\phi}{\SI \times [0,1]}$ defined by $$ (z,t) \mapsto  (z
e^{a(4t-2)\pi i},t). $$ Also, by \lemc{one fixed}(1),  $f^2|_{C_1}$ has no fixed points.
(Observe that the restriction of $f^2$ to the boundary component of $C_1$ is a rotation
of order $\sigma_1 \ge 4g_1>4)$. Similarly,  $f^2|_{f(C_1)}$ has no fixed points.

By arguments analogous to those used in the preceding paragraph, we can prove that
$L(f^2)=0$. Since  this contradicts the  hypotheses of the proposition, the proof of this
case is complete.
\end{proof}

\begin{rem} Although in the preceding proof we use the fact that the map on the annulus $A$
is the restriction of a reducible map in standard form, it is a simple matter to check
that any orientation-preserving map of an annulus which leaves invariant each boundary
component (resp. interchanges both boundary components) has Lefschetz number equal to $2$
(resp. $0$).
\end{rem}

\begin{prooftext}{Proof of \theoc{resumenop}} Assume first that $L(f)=2$ and $L(f^2)=0$. In this case, by \theoc{lef}, $\Fix(f) \ne \emptyset$. Since $\Fix(f) \subset \Fix(f^2)$,  $\Fix(f^2) \ne \emptyset$. Thus, if $f^2$ does not have a fixed-point class of index $0$, by  \lemc{suma},  $f^2$ has a fixed-point class of negative index. Since $2 \le 4g+2$, the result holds for this case. Now, assume that $L(f) \ne 2$ or $L(f^2) \ne 0$ and let us prove that $f^m$ has a fixed-point class of negative index for some positive integer $m$ such that $m \le 4g+2$.
By \lemc{class} and \theoc{jiang} we can also assume that $f$ is in standard form.

We know that there are three possibilities for $f$, namely, it can be of finite-order,
pseudo-Anosov or reducible. If $f$ is finite-order, then there exists a positive integer
$n$ such that $f^n=\id$. By \theoc{www}, we can take $n \leq 4g+2$. So,
$L(f^n)=L(\id)=2-\tr (\id)=2-2g <0$. Hence, the fixed-point class is all of $\Su{g}$  and
its index is $L(f^n)$.

The remaining cases follow from Propositions~\ref{pA}, \refc{reducibleop} and
\refc{reduop}.
\end{prooftext}

\index{Theorem H}
\begin{proclama}{Theorem H}
If $g \ge 2$ then $\mmm(\hgbp) \le 4g+2$. Moreover, if $b \geq 6g+6$, then equality
holds.
\end{proclama}
\begin{proof}
Let $g \ge 2$.  Observe that, by \coryc{idioma},  $\mmm(\hgbp) \ge 4g+2$ if $b \ge 6g+6$.
To complete the proof the theorem, it suffices to show that $\mmm( \hgbp) \le 4g+2$.

Let $f \in \hgbp$ and consider the induced map $\Smap{\tf}{\Su{g}}$. Let $m$ be as in
\theoc{resumenop} for $\tf$. Since $m \le 4g+2$, it is enough to prove that $\mmm(f) \le
m$. Consider a fixed-point class  $C \subset \Su{g}$  of $\tf^m$ of  non-positive index.
If $C$ is  finite,  by \lemc{elolvidado}, the index of $C$ with respect to $\tf^m$ is the
sum of the indices of each of its elements with respect to $\tf^m$. Since this sum is
non-positive, at least one of its terms must be non-positive. Thus  $\tf^m$ has a fixed
point of non-positive index. By \propc{index}, $f^m $ has a fixed point, so $\mmm(f) \le
m$

If $C$ is  infinite, it contains points which are not collapsed boundary components of
$\Suf{g}{b}$. Since the existence of these fixed points of $\tf^m$ implies the existence
of fixed points of $f^m$, the proof of Theorem~H is complete.
\end{proof}

\start{rem}{nada} Observe that if $f \in \hgbp$ is such that $\mmm(f) >4g$ then $\tf$ is
isotopic to a finite-order map.
\end{rem}

\start{conj}{xxx} If $g \ge 2$, there exists $f \in \hgbp$ such that $\tf$ is isotopic to
a finite-order map and $\mmm(f)=\mmm(\hgbp)$.

This holds if $b \ge 6g+6$ or if $g=2$.
\end{conj}

\section{Proof of Theorem~I}\label{I}

As in the previous section, our main objective is to prove  the following.

\start{theo}{resumenor} Let $g \geq 2$ and let $f \in \hgr$. Then there exists a positive
integer $m$ such that  $m \le 4g+(-1)^g4$ and  $f^m$ has a fixed-point class of  negative
index.
\end{theo}

To prove \theoc{resumenor} we will use the following results  which study the
pseudo-Anosov case, and the reducible case both with and without pseudo-Anosov
components.

\start{prop}{pAgodd} Let $g \geq 2$ be odd. If $\Smap{f}{\Su{g}}$ is an
orientation-reversing pseudo-Anosov map in standard form then there exists a positive
integer $m$ such that $m \le  4g-4$ and  $f^m$ has a fixed-point class of negative index.
 \end{prop}

\begin{proof} As in  the proof of \propc{pA}, we will show that one of the statements of \lemc{laberinto} applies for some $m$ such that $1 \le m \le 4g-4$.

We can repeat the first part of the proof of \propc{pA}. Combining the results obtained
there with \lemc{laberinto}(1), (5) and (6) we can assume
\begin{romlist}
\item $k \in \{1,2\}$.
\item If $k=2$, then $n_1 \le n_2$ and $p_2=3$.
\item $\ind{f^{n_1h}}{x_1}=1$ for every $h$ such that $n_1h \le g$ and $n_1h$ is even.
\item There are no regular points of period less than or equal to $2g-2$.
\item If for some $i \in \{1,2\}$, $n_i \le 2g-2$ then $n_i$ is even.
\end{romlist}

Now we split the proof into four cases in order to see that in each of them, one of the
statements of \lemc{laberinto} applies.

\begin{case}{1}{$k=2$, and $n_1, n_2 \le g-1$.} \end{case}
Here $$ n_1 p_1+n_2 p_2 \leq 4(g-1)+2(n_1+n_2) \leq 4(g-1)+4g-4=8g-8. $$ Therefore, we
may assume $n_1p_1 \le 4g-4$. By (v), \lemc{laberinto}(4) applies with $m=n_1p_1$.

\begin{case}{2}{ $n_1>g$.} \end{case}
By (iv) and \theoc{lef}, $$L(f)=L(f^2)=\ldots=L(f^g)=0.$$ By Lemmas~\ref{pol-or}
and\refc{symmetry}, $P(x)=(x^2-1)(x^{2g-2}+1)$, so the $(4g-4)$-th power of the
eigenvalues is $1$.  By \lemc{potency},  $L(f^{4g-4})<0$, so \lemc{laberinto}(2) applies.

\begin{case}{3}{$n_1=2$, and  $n_2 > g$ if $k=2$.} \end{case}
By \theoc{lef}, \theoc{suma} and (iv),
\begin{eqnarray*}
L(f^i) & = & \cases{0 &  \mbox{if $i$ is odd,}\cr
                    2 &   \mbox{otherwise}.}\cr
\end{eqnarray*}
By  \refeq{numerar},  $p_i=0$ for $1 \le i \le g$. Now,  Newton's equations (page
\pageref{newton}) show  that $P(x)=x^{2g}-1$, so the $(2g)$-th power of each eigenvalue
is $1$. By \lemc{potency}, $L(f^{2g})<0$,  so \lemc{laberinto}(2) applies with  $m=2g \le
4g-4$.

\begin{case}{4}{$3 \le n_1 \le g-1$ and $n_2 >g$  if $k=2$.} \end{case}

By (v), $n_1$ is even and by (iv), there are no periodic regular points of period less
than or equal to $g$ because $g \le 2g-2$. Thus by  (iii), \propc{cas-part-or-odd}
applies. Now, we can complete the proof of the proposition in this case as in  Case 4 of
the proof of \propc{pA}.
\end{proof}

\start{prop}{pAgeven} Let $g \ge 2$ be even. If  $\Smap{f}{\Su{g}}$ is an
orientation-reversing pseudo-Anosov map in standard form, then there exists a positive
integer $m$ such that $m \le 4g+4$ and $f^m$ has a fixed-point class of negative index.
 \end{prop}

\begin{proof} As in  the proof of \propc{pA}, we will show that one of the statements of \lemc{laberinto} applies for some positive integer $m$ such that $m \le 4g+4$.

We can repeat the first part of the proof of \propc{pA}. Combining the results obtained
there with \lemc{laberinto}(5) and (6) we obtain
\begin{romlist}
\item $k \in \{1,2\}$.
\item If $k=2$, then $n_1 \le n_2$.
\item $\ind{f^{n_1h}}{x_1}=1$ for every $h$ such that $n_1h \le g$.
\item There are no regular points of period less than or equal to $2g+2$.
\item If for some $i \in \{1,2\}$, $n_i \le 2g+2$ then $n_i$ is even.
\end{romlist}
We split the proof into four cases.

\begin{case}{1}{ $n_1 \in \{2,4\}$.} \end{case}
Here $n_1p_1 \le 4g-4+2n_1 \le 4g-4+8=4g+4$. Then \lemc{laberinto}(4) applies.

\begin{case}{2}{ $k=2$, and $n_1, n_2 \le g+2$.} \end{case}
Here $$ n_1 p_1+n_2 p_2 \leq 4(g-1)+2(n_1+n_2) \leq 4(g-1)+4g+8=8g+4. $$ We can assume
$n_1p_1 \le 4g+2$. Since $g+2 \le 2g+2$, by (v), $n_1$ is  even. Then \lemc{laberinto}(4)
applies.

\begin{case}{3}{ $n_1 \ge g+1$.}{}\end{case}
Since $g \le 2g+2$, by (iv) and \theoc{lef}, $$L(f)=L(f^2)=\ldots=L(f^{g})=0.$$ Now, by
Lemmas~\ref{pol-or} and\refc{symmetry}, $P(x)=(x^2-1)(x^{2g-2}-1)$. Therefore, the
$(2g-2)$-th power of each eigenvalue is equal to $1$. By \lemc{potency}, $L(f^{2g-2})<0$,
so \lemc{laberinto}(2) applies.

\begin{case}{4}{$5 \le n_1 \le g$, and $n_2 >g+2$, if $k=2$.} \end{case}
By (v), $n_1$ is even. By (iv) there are no periodic regular points of period less than
or equal to $g+2$, because  $g+2 \le 2g+2$. By (iii),  the hypotheses of
\propc{cas-part-or-even} hold. Therefore there exists a positive integer $m$ such that $m
\le 2g-6$ and $L(f^m)<0$. Since $2g-6<4g-4$, \lemc{laberinto}(2) applies in this case.
\end{proof}

\start{prop}{reducibleor} Let $\Smap{f}{\Su{g}}$ be an orientation-reversing reducible map in standard form
which has a pseudo-Anosov component. Then there exists a positive integer $m$ such that
$m  \le 4g-4$  and $f^m$ has a fixed-point class  of negative index.
\end{prop}
\begin{proof} We shall prove that one of the statements of \lemc{tiger} applies for some positive integer $m \le 4g-4$.

Let $C \subset \Su{g}$ be a pseudo-Anosov $f$-component of period $n_1$, genus $g_1$, and
$b_1$ boundary components. Let  $B$ be a  boundary component of $C$. We can repeat the
first part of the proof of \propc{reducibleop} to prove
\begin{equation}
p_{_B}n_1r_{_B} \le 4g-4,\label{pillow}
\end{equation}
where $r_{_B}$ is the least positive integer such that $f^{n_1r_{_B}}(B)=B$, and $p_{_B}$
is the number of prongs emanating from $B$. If $n_1r_{_B}$ is even, then, by
\refeq{pillow}, \lemc{tiger}(1) applies. Hence we can assume that
\paragrafetiq{(i)}{\em For each boundary component $B$ of $C$, $n_1r_{_B}$ is odd.}

If $p_{_B} \ge 2$ for some  boundary component $B$ of $C$ then by \refeq{pillow},
$2n_1r_{_B} \le  p_{_B}n_1r_{_B} \le 4g-4$. Hence, \lemc{tiger}(2) applies, so we can
assume that
\paragrafetiq{(ii)}{\em For every boundary component $B$ of $C$, $p_{_B}=1$.}

Clearly, $r_{_B} \le b_1$. Suppose now that $g_1 \ge 1$. Then by \lemc{corep}(1), $$
2r_{_B}n_1 \le 2b_1n_1 \le 2n_1(2g_1+b_1-2) \le 4g-4. $$ Now, by (i), \lemc{tiger}(3)
applies. On the other hand, if $g_1=0$, by \lemc{b4}. $b_1 \ge 4$. Hence, by
\lemc{corep}(1),  $2n_1 \le n_1(b_1-2) \le 2g-2$. Thus we can assume
\paragrafetiq{(iii)}{\em $g_1=0$,  $2n_1 \le n_1(b_1-2) \le 2g-2$.}

Now, observe that by (ii) and \remc{pprong}, $\sum_{s \in \Sing(B)}(2-p_s)=-1$ for each
boundary component $B$ of $C$. Thus $\sum_{s \in \Sing(\partial C)}(2-p_s)=-b_1$. By the
Euler-Poincar\'{e} Formula \refeq{ep} for $f^{n_1}|_{C}$ $$ \sum_{s \in \Sing(\Int
C)}(2-p_s)-b_1=\sum_{s \in \Sing(\Int C)}(2-p_s)+\sum_{s \in \Sing(\partial
C)}(2-p_s)=2(2-b_1). $$ So
\paragrafetiq{(iv)}{\em $\sum_{s \in \Sing(\Int C)}(p_s-2)=b_1-4$.}

By (iii), $f^{2n_1}|_{C_1}$ induces an orientation-preserving map of the sphere $\Su{0}$.
By Brouwer's theorem (or by \theoc{lef}, since $H_1(\Su{0})$ so $L(\tf^{2n_1})=2$), this
map has a fixed point $x$.  If $x$ is a collapsed boundary component $B$, then
$f^{2n_1}(B)=B$. Then \lemc{tiger}(1)  applies for $m=4n_1$, and,  by (iii),  $4n_1 \le
4g-4$. Thus, we can assume that $x$ is not a collapsed boundary component. In this case,
$f^{2n_1}|_C$ has  a fixed point $y \in \Int(C)$. If $y $ is a regular point,
\lemc{tiger}(1) applies with $m=4n_1$. As before, by (iii), we see that $m \le 4g-4$.
Hence, we can assume that  $y$ is a singularity. Denote by $p_1$ the number of prongs
emanating from $y$. By (iv), $$ 2(p_1-2) \le \sum_{s \in \Sing(\Int C)}(p_s-2)=b_1-4. $$
So, $2p_1 \le b_1$. Thus, $2p_1n_1\le b_1n_1=n_1(b_1-2)+2n_1 \le 4g-4.$ Since
$f^{2n_1}(y)=y$, \lemc{tiger}(4) applies and the proof of the proposition is complete.
\end{proof}

The following proposition is the analogue of  \propc{reduop} for the
orientation-reversing case.

\start{prop}{reduor} Let $g \ge 2$ and let $f \in \hgr$ be a reducible homeomorphism in
standard form such that each of its components  is finite-order.  Then there exist an
$f$-component $C$ and a positive  integer $m$ such that $m \le 4g+(-1)^g4$ and
$f^m|_C=\id_C$.
\end{prop}
\begin{proof} With \notc{dan}, as in the proof of \propc{reduop}, it suffices to show that there exists $\ii{i}{k}$ such that $n_i \si \le 4g+(-1)^g4$. We split the proof into various cases.

\begin{case}{1}{$k \ge 3.$}
\end{case}
By \lemc{corep}(1), $\sum_{i=1}^k (2g_i+b_i-2)n_i=2g-2.$ Changing subindices if
necessary,  by \lemc{corep}(2), we can assume
\begin{equation}
n_1 \le n_1(2g_1+b_1-2) \le \frac{2}{3}(g-1). \label{dream}
\end{equation}
Now we split the proof of this case into  two subcases.

\begin{case}[Subcase]{1.1}{$g_1 \in\{0,1\}$}\end{case}
We claim that the result holds if  $\sigma_1 \le 2b_1$. Indeed, by \refeq{dream} $$
n_1\sigma_1 \le2b_1n_1=2((2g_1+b_1-2)n_1+2n_1) \le 2(\frac{2}{3}(g-1)+
\frac{4}{3}(g-1))=4(g-1). $$ Hence, we can assume that $\sigma_1 >2b_1$. Now, consider
$f^{n_1}|_{C_1}$. By Lemmas~\ref{s0b},  \refc{s0bro} and \refc{s1bro}, $f^{n_1}$ must be
orientation preserving and $g_1=1$. So $n_1$ is even. By \lemc{numbertoro}, $\sigma_1 \le
6$. Hence, by \refeq{dream} $$ \sigma_1 n_1 \le 6 n_1 \le 4g-4, $$ which completes the
proof of this subcase.

\begin{case}[Subcase]{1.2}{$g_1 \ge 2$}\end{case}
Since $b_1 \ge 1$, $$ 3n_1 \le (2g_1-1)n_1 \le (2g_1+b_1-2)n_1 \le \frac{2}{3}(g-1). $$
Thus, by \coryc{wwwc}, $\sigma_1 \le 4g_1+4$ and, since $b_1 \ge 1$, by \refeq{dream} $$
n_1 \sigma_1 \le  n_1(4g_1+4)=2(n_1(2g_1-1)+3n_1) \le 2(\frac{2}{3}(g-1)
+\frac{2}{3}(g-1)) \le 4(g-1), $$ as desired.

\begin{case}{2}{$k=2$.}
\end{case}
We consider three subcases.

\begin{case}[Subcase]{2.1}{$g_1, g_2 \ge 1$.}\end{case}
Since, by \lemc{euler},  $n_1g_1+n_2 g_2 \le g$,  we can assume without loss of
generality that $n_1g_1 \le  g/2$.

If $g_1=1$ then by Lemmas~\ref{itemsor}(3) and \ref{items}(2),  the result holds if $n_1$
is odd or $\sigma_1 \ge 7$. Now, suppose that $n_1$ is even and $\sigma_1 \le 6$. In
particular, $2 \le n_1\sigma_1 \le g/2$. Then $g \ge 4$ and $$ n_1\sigma_1 \le 6n_1 \le
3g \le 4g-4, $$ as desired. If $g_1 \ge 2$, then  $2  \le n_1g_1 \le g/2$. Hence, $g \ge
4$, and by \coryc{wwwc}, $$ n_1\sigma_1 \le n_1 (4g_1+4) \le 2g+g \le 4g-4, $$ and the
proof of this subcase is complete.

\begin{case}[Subcase]{2.2}{$g_1= g_2 =0$.}\end{case}
By \lemc{corep}(1),
\begin{equation}
n_1(b_1-2)+n_2(b_2-2)=2g-2.  \label{dream1}
\end{equation}
If  $n_1$ and $n_2$ are  even,  then, by \refeq{dream1}, $n_i(b_i-2) \le g-1$ for some $i
\in \{1,2\}$ and the result follows from \lemc{items}(1).

If $n_1$ and $n_2$ are odd, by \lemc{itemsor}(2), we can assume that, for each $i \in
\{1,2\}$,  $b_i=\si/2$ and $b_i$ is odd. By \lemc{s0bro}, the boundary components if
$C_i$ form a cycle under the action of $f^{n_i}$. Therefore, they remain fixed under the
action of $f^{n_ib_i}$. Hence, the $f$-period of each connected component of $N(\Gamma)$
is odd. By \lemc{sandalo}, $f$ is finite-order. Then we take the whole $\Su{g}$ as $C_1$
and the result holds by \theoc{www}.

To complete the proof, changing subindices if necessary, we can assume that $n_1$ is odd
and $n_2$ is even. Then, by \refeq{dream1},  $b_1$ must be even and the result follows
from \lemc{itemsor}(2).

\begin{case}[Subcase]{2.3}{$g_1=0$, $g_2 \ge 1$.}\end{case}

By  \lemc{itemsor}(2) and (4), we can assume that $n_2(2g_2+b_2-2)$ is even. Then,
$n_1(b_1-2)$ is even. If $n_1$ is odd,  the result hold by \lemc{itemsor}(2). If $n_1$ is
even, by \lemc{items}(1), we can assume that $b_1=3$ and $n_1>\frac{4}{3}(g-1)$. By
\lemc{corep}(1), $$ n_2(2g_2+b_2-2)=2g-2-n_1(b_1-2) < \frac{2}{3}(g-1). $$ and we can
complete the proof as we did in Case 1.

\begin{case}{3}{$k=1$, $g_1=0$.}\end{case}
If $n_1=1$ the desired conclusion follows from \coryc{wwwc}. Therefore, we can assume
that $n_1>1$. By \lemc{corep}(1), $n_1(b_1-2)=2g-2$.

Assume that $n_1$ is even. By \lemc{items}(1), we can assume that $\sigma_1=b_1=3$ and
the three boundary components of $C_1$ form a cycle under the action of $f^{n_1}$. Then
we can apply \lemc{ciclo} to obtain  $n_1=2$ and $g=2g_1+b_1-1=2$. Hence, $\sigma_1n_1=6
\le 8=4g+(-1)^g4$.

If $n_1$ is odd, since $(b_1-2)n_1=2g-2$, $b_1$ is even and the result follows from
\lemc{itemsor}(2).

\begin{case}{4}{$k=1$, $g_1 \ge 1$.}\end{case}
Again we consider three subcases.
\begin{case}[Subcase]{4.1}{$b_1=1$.}\end{case}
By \lemc{ciclo}, $n_1=2$ and $g=2g_1+b_1-1=2g_1$. In particular, $g$ and $n_1$ are even.
Then, by \coryc{wwwc}, $\sigma_1 n_1 \le (4g_1+2)n_1=4g+(-1)^g4$.

\begin{case}[Subcase]{4.2}{$b_1=2$.}\end{case}
If the boundary components of $C_1$ form a cycle under the action of $f^{n_1}$ then, by
\lemc{ciclo}, $n_1=2$, so $C_1$ and $f(C_1)$ are connected by  two annuli $A_1$ and
$A_2$. Since $f(A_1)=A_1$ or $A_2$, $f^2(A_1)=A_1$ and  $f^2(A_2)=A_2$. Then the two
boundary components of $C_1$ cannot form a cycle form a cycle under the action of
$f^{n_1}$, which contradicts our assumption. Therefore, both boundary components of $C_1$
are mapped to themselves under the action of $f^{n_1}$. By \lemc{items0}(3), we can
assume that $\sigma_1 \ge 3$. By \lemc{one fixed}(2), $f^{n_1}|_{C_1}$ must be
orientation preserving, so $n_1$ is even and, by \lemc{one fixed}(1), $\sigma_1 \le
4g_1$. Then $n_1\sigma_1 <4g_1n_1=4g-4$.

\begin{case}[Subcase]{4.3}{$b_1=3$.}\end{case}
If $g_1=1$ the result follows from \lemc{items}(2). Now, let $g_1 \ge 2$. By
\lemc{corep}(1), $$ n_1(2g_1+1)=n_1(2g_1+b_1-2)=2g-2. $$ This implies that $n_1$ is even,
and \lemc{items}(1) and (2) completes the proof.
\end{proof}

\begin{prooftext}{Proof of \theoc{resumenor}}
By \lemc{class} and \theoc{jiang} we can assume that $f$ is in standard form.

We know that there are three possibilities for $f$, namely, it can be of finite-order,
pseudo-Anosov or reducible. If $f$ is finite-order, then there exists a positive integer
$n$ such that $f^n=\id$. By \theoc{www}, we can take $n \leq 4g+(-1)^g4$, so  by
\refeq{anne}, $L(f^n)=L(\id)=2-\tr (\id)=2-2g <0$. In this case, the fixed-point class is
all of $\Su{g}$  and its index is $L(f^n)$.

The remaining cases are consequences of Propositions~\ref{pAgodd}, \refc{pAgeven},
\refc{reducibleor} and \refc{reduor}.
\end{prooftext}

\index{Theorem I}
\begin{proclama}{Theorem I} Let $g \geq 2$. Then  $\mmm(\hgbr)  \le 4g+(-1)^g4$ and equality holds if  $b \ge 6g+2+(-1)^g8$.
\end{proclama}

\begin{proof}
By \propc{sabiduria}, $\mmm(\hgbp) \ge 4g+4$ if $b \ge 6g+10$ and $g$ is even. By
\coryc{prisma},  $\mmm(\hgbp) \ge 4g-4$ and $b \ge 6g-6$ if $g$ is odd. Hence, to
complete the proof of the theorem, it suffices to show that $\mmm( \hgbp) \le
4g+(-1)^g4$.  Now, we can complete the proof as we did for Theorem~H.

 \end{proof}

\printindex

\enddocument